\documentclass[10pt,reqno]{amsart}
\usepackage{mathrsfs}
\usepackage{amsmath,amsfonts,amssymb,amsthm}
\usepackage{pst-plot}
\usepackage{graphicx}
\usepackage{times}
\usepackage{verbatim}
\usepackage{hyperref}
\usepackage{cleveref}
\usepackage{nameref}
\usepackage{csquotes}
\usepackage{a4wide}
\usepackage{subfig}
\usepackage{resmes}
\usepackage[english]{babel}

\allowdisplaybreaks
\newcommand{\I}{{\bf 1}}

\newtheorem{proposition}{Proposition}[section]
\newtheorem{theorem}[proposition]{Theorem}

\newtheorem{lemma}[proposition]{Lemma}
\newtheorem{remark}[proposition]{Remark}

\newtheorem{example}[proposition]{Example}
\newtheorem{examples}[proposition]{Examples}

\newcommand{\nc}{\newcommand}
\nc{\FN}{{\mathbf N}}
\nc{\cB}{{\mathcal B}}
\nc{\R}{{\mathbb R}}
\nc{\N}{{\mathbb N}}
\nc{\Z}{{\mathbb Z}}

\DeclareMathOperator{\relint}{relint}

\DeclareMathOperator{\interior}{int}
\DeclareMathOperator{\conv}{conv_s}
\DeclareMathOperator{\cl}{cl}

\DeclareMathOperator{\ext}{ext}
\DeclareMathOperator{\dom}{\mathfrak{I}}

\DeclareMathOperator{\proj}{\Pi}
\DeclareMathOperator{\fracc}{\mathfrak{c}}
\newcommand{\sdn}{\sigma_d^\circ}
\newcommand{\sden}{\sigma_{d-1}^\circ}
\newcommand{\Ca}{{\sf C}_a}
\newcommand{\ZV}{Z_V}
\nc{\dint}{\mathrm{d}}
\nc{\mrG}{\mathrm{G}}
\nc{\SO}{\mathrm{SO}}
\nc{\He}{H^+(e)}
\nc{\Hp}{H^*(e)}
\nc{\bo}{\mathbf{o}}
\nc{\gas}{\gamma_S}
\nc{\BP}{\mathbb{P}}
\nc{\BE}{\mathbb{E}}
\nc{\BQ}{\mathbb{Q}}
\nc{\BS}{\mathbb{S}}
\nc{\cH}{{\mathcal H}}
\nc{\cS}{{\mathcal S}}
\nc{\cK}{{\mathcal K}}
\nc{\oo}{\overline{0}}
\nc{\sfp}{\mathsf{p}}
\nc{\indi}{\mathbf{1}}
\nc{\Ih}{{\mathcal{I}_d}}
\newcommand{\defeq}{\mathrel{\mathop:}=}
\nc{\bI}{\mathbf{1}}
\nc{\1}{\bI}
\nc{\cC}{\mathcal{C}}
\DeclareMathOperator{\TG}{G}
\begin{document}
\title[Geometric inequalities, stability  and Kendall's problem]{Geometric inequalities, stability results and Kendall's problem in spherical space}
\author[Daniel Hug and Andreas Reichenbacher]
{Daniel Hug and Andreas Reichenbacher}
\address{Karlsruhe Institute of Technology, Department of Mathematics,
D-76128 Karlsruhe, Germany}
\email{daniel.hug@kit.edu}\email{reichenbacher\underline{\phantom{x}}andreas@web.de}
\thanks{Authors supported in part by DFG grant FOR  1548}
%\maketitle
%\date{\today}
\date{\today}
\pagestyle{myheadings}
\markboth{Daniel Hug and Andreas Reichenbacher}{Kendall's problem in spherical space}

\begin{abstract}
In Euclidean space, the asymptotic shape of large cells in various types of Poisson driven random tessellations
has been the subject of a famous conjecture due to David Kendall. Since shape is a geometric concept and large cells
are identified by means of geometric size functionals, the resolution of the conjecture is inevitably connected with geometric
inequalities of isoperimetric type and their  improvements in the form of geometric stability results, relating geometric size functionals and hitting functionals. The latter  are deterministic  characteristics of the underlying random tessellation. The current
work explores specific and typical cells of random tessellations in spherical space. A key ingredient of our approach are new geometric
inequalities and quantitative strengthenings in terms of stability results for  general and also for some specific size and hitting functionals of spherically  convex bodies. As a consequence we obtain probabilistic deviation inequalities and asymptotic distributions of quite
general size functionals. In contrast to the Euclidean setting, where naturally the asymptotic regime concerns large size, in the spherical framework the asymp\-totic analysis is primarily concerned with high intensities.
\end{abstract}

\maketitle
%

%
%
%

%
%%%%%%%%%%%%%%%%%%%%%%%%%%%%%%%%%%%%%%%%%%%%%%%%%%%%%%%%%%%%%%%%%%%%%%%
%%%%%%%%%%%%%%%%%%%%%%%%%%%%%%%
%\tableofcontents
%%%%%%%%%%%%%%%%%%%%%%%%%%%%%%%%%%%%%%%%%%%%%%%%%%%%%%%%%%%%%%%%%%%%%%
%%%%%%%%%%%%%%%%%%%%%%%%%%%%%%%
\section{Introduction}\label{intro}

Deterministic tessellations, or mosaics, have been a subject of interest for a very long time. Even ancient cultures, like the Sumerians or
the Romans, used colored tiles to decorate floors and walls. Formal, mathematical definitions and deterministic tilings of the
plane or higher dimensional spaces were considered much later. A tessellation of Euclidean space $\R^d$ is usually understood to be a system of closed sets (often these are compact, convex
polytopes with nonempty interiors) in $\R^d$ which cover the whole space, have pairwise no common interior points, and do not accumulate locally.

Random tessellations of Euclidean space are a classical topic in stochastic geometry. They have been extensively studied in the
literature; see, e.g., \cite[Chap.~9]{SKM}, \cite[Chap.~10]{SW2008} and \cite{M89,MS,ReJu2025} for an overview and results for
general tessellations. In a large number of contributions,  various properties of special models have been explored, in a deterministic or random framework, including hyperplane
tessellations, Voronoi tessellations, Delaunay tessellations, Laguerre tessellations and generalizations of these models such as $\beta$-Voronoi tessellations
 \cite{GKT2024} and generalized balanced Voronoi tessellations (see \cite{JuRe2024} and the survey \cite{ReJu2025}); see also \cite{BL1, BL2, HS24, Mdiss, Miles74, M, Muche05, RedenbachLiebscher2015, Voss2013} and \cite[Chap.~6]{matheron}.

\begin{figure}
\begin{center}
\includegraphics[width=8cm,angle=0]{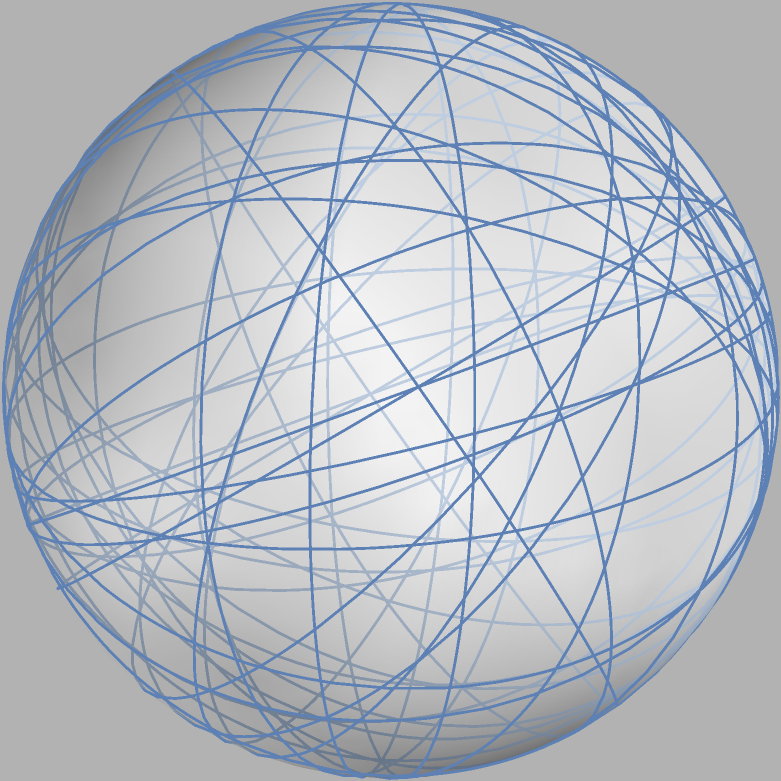} 
\end{center}
\caption{Tessellation of the unit sphere by 30 uniformly sampled hyperspheres (generated by code based on Wolfram Mathematica)}\label{Fig1}
\end{figure}

%
%----------------------------------
%
In the present work, we consider random tessellations of the unit
sphere $\BS^d$ in Euclidean space $\R^{d+1}$. This setting has not been studied as extensively in the literature as the Euclidean framework. The
intersection of the unit sphere with a $d$-dimensional linear subspace is the unit sphere
in the intersecting subspace and thus a
hypersphere (great subsphere) of $\BS^d$ with unit radius.
At the same time, $d$-dimensional linear subspaces
partition the Euclidean space $\R^{d+1}$ into polyhedral cones. This relation plays an important role in spherical random geometry,
see, e.g., \cite{AM,CE,KT21b,KTT21,Schneiderneu} and \cite{S22cones} (and the literature cited there).
Tessellations of the sphere generated by intersecting the unit sphere with $d$-dimensional linear subspaces are called spherical hyperplane (or hypersphere)  tessellations. Random  hypersphere tessellations, where the subspaces are selected randomly, are studied in
\cite{AZ, CE}, \cite[Sec.~6]{Miles}, and in recent work on conical tessellations \cite{HScone, HS21, HS22} (see Fig.~\ref{Fig1} for an illustration). A cell splitting scheme on $\mathbb{S}^2$, related to Poisson processes of hyperspheres, has been considered
in \cite{DHT}, a systematic study of splitting tessellations in spherical space (in analogy to iteration stable tessellations in Euclidean space) is carried out in \cite{HT+}.

\begin{figure}
\begin{center}
\includegraphics[width=8cm,angle=0]{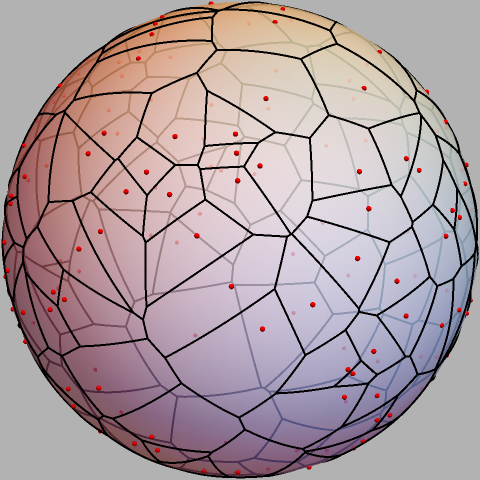} 
\end{center}
\caption{Voronoi tessellation of the unit sphere by a sample of 150 uniform random points on the unit sphere (generated by code based on Wolfram Mathematica)}\label{Fig2}
\end{figure}

Voronoi tessellations in spherical space can be defined as in the Euclidean case, using the geodesic distance on $\BS^d$ to determine the individual cells of the tessellation. Random Voronoi tessellations on the sphere and their applications have been investigated, for instance,  in \cite{Chaidee2018,Chaidee2020,KT21, Renk,Sugi,Tan,Yoshino2012}, \cite[Sec.~7]{Miles}, and  \cite[Sec.~3.7.6, Sec.~5.10]{Okabe}. Random Voronoi tessellations have also been studied in hyperbolic space \cite{Isokawa20a,Isokawa20b} and in more general Riemannian spaces. Fig.~\ref{Fig2} provides an illustration for a Voronoi tessellation of the unit sphere. From a foundational viewpoint, Voronoi diagrams have been explored in Riemannian manifolds in \cite{Innami15,Leibon02,Leibon,Nielsen}, connections to
biomedical imaging are discussed in \cite{Jones}, aspects of computational geometry and information theory are treated in \cite{Boissonat2010,Boissonat2019,Nass}, random Poisson--Voronoi tessellations on surfaces and Riemannian manifolds are studied in \cite{Calka19,Calka21}. Further motivation to consider Poisson--Voronoi tessellations in Riemannian spaces arises from intriguing problems in percolation theory (see \cite{Benjamini01,Benjamini18,Hansen22,Hansen24,Buehler}).

In the following, we focus on what became known as `Kendall's Problem' or `Kendall's Conjecture' and in particular
on its geometric foundations (see \cite[Chap.~12]{HS24} for surveys and background information, \cite[Note 9 for Sec.~10.4]{SW2008}, and  \cite{BonnetDiss2017,BonnetCalkaReitzner2016,Calka2013,Hug2013,HugReitzner2016,Reitzner2010,Schnsurvey2}).
Prior to our work, this line of investigation
has been exclusively explored in the Euclidean setting. In the present investigation, we formulate and study a spherical analogue.
An analysis of the problem in hyperbolic space has been initiated in \cite{Herold, HH+}. To recall
the problem in Euclidean space, consider a
stationary and isotropic Poisson line process in the Euclidean plane and denote the almost surely unique cell containing the
origin by $Z_0$. This cell is called the {\em Crofton cell} or {\em zero cell}. In the foreword
of the first edition of \cite{SKM} from 1987, David G. Kendall stated the following conjecture: The conditional law for the shape of $Z_0$,
given the area $A(Z_0)$, converges weakly, as $A(Z_0)\rightarrow\infty$, to the degenerate law concentrated at the circular
shape. This conjecture was strongly supported by heuristic arguments from R. Miles \cite{Milesheuristic}.
Two years later,  Kovalenko gave a proof in \cite{Kov1}. Kovalenko also provided a simplified proof in \cite{Kov2} and
an extension to the typical cell of a Poisson--Voronoi tessellation in the plane in \cite{Kov3}.
Further extensions to arbitrary dimensions and not necessarily isotropic Poisson
hyperplane tessellations were made in \cite{HRS}, where the size of the Crofton cell was measured by the volume. In \cite{HRS2}
the problem was extended and solved for typical cells of stationary Poisson--Voronoi tessellations in arbitrary dimensions and
the size was measured by an intrinsic volume. In \cite{LC} a very general setting with a very general class of size functionals
was considered, containing the aforementioned results as special cases, streamlined statements and arguments can be found in \cite[Chap.~12]{HS24}. In \cite{HSf}, Kendall's problem was
extended to the typical $k$-faces of a Poisson hyperplane tessellation ($k\in2,\ldots,d-1\}$) and in \cite{HSpv} to the typical
$k$-faces of a Poisson--Voronoi tessellation. In \cite{HSdel} typical cells of Poisson--Delaunay tessellations were considered.
In all these previous works, geometric inequalities, stability results and polytopal approximation have been crucial geometric ingredients in
addition to the probabilistic analysis. References to various applications, including extreme value theory and statistical learning theory, of the results obtained in these works are highlighted in \cite[Notes for Sec.~12.3]{HS24}.

We continue with some notation in order to present selected new results. Since on the unit sphere $\BS^d$ there is no naturally
distinguished point similar to the Euclidean origin in $\R^{d+1}$,  we choose an arbitrary fixed point $\sfp\in\BS^d$
as the {\em spherical origin}.
Let $d_s$ denote the geodesic metric on $\BS^d$, and let $B_s(x,r)$ denote the closed spherical (geodesic) ball with centre $x$ and radius $r\le \pi$. If
$r\le \pi/2$, then we call $B_s(x,r)$ a spherical cap.
A {\em proper (spherically) convex body} in $\BS^d$ is the intersection of the unit sphere with some nonempty
line-free closed convex cone in $\R^{d+1}$ that does not only consist of the Euclidean origin $o\in\R^{d+1}$. We denote the set of proper convex bodies by $\cK_s^d$.
If we do not require the cone to be line-free but only that it is not equal to some linear subspace of
$\R^{d+1}$, then the resulting set will be denoted by $\overline{\cK}_s^d$ and its elements are called {\em (spherically) convex bodies}. A {\em spherical polytope} is the intersection of $\BS^d$ with a polyhedral cone (a finite intersection of half-spaces whose bounding hyperplanes contain the origin) which is also a spherically convex body. For more details on spherical geometry, we refer to \cite{Glas} and \cite[Sec.~6.5]{SW2008}.

By a tessellation of $\BS^d$ we mean a finite collection of spherical polytopes that have nonempty interiors,
cover $\BS^d$ and have pairwise disjoint interiors, or the trivial tessellation consisting of $\BS^d$ only. Prominent examples are hypersphere and Voronoi tessellations which are based on a
finite number of hyperspheres and point sets, respectively.  Random tessellations are obtained, for instance, by selecting the underlying hyperspheres and points randomly, thus giving rise to random hypersphere and random Voronoi tessellations, respectively.
In general, a spherical random tessellation is said to be isotropic, if its distribution is invariant under any rotation in $\SO_{d+1}$.
In this case, almost surely there exists a unique cell containing $\sfp$ in its interior. We call this cell the {\em spherical Crofton cell} (or spherical zero cell) of the given random tessellation and denote it by $Z_0$.

The spherical Lebesgue measure on $\BS^d$ will be denoted by $\sigma_d$ and the surface area of the unit sphere by
$\omega_{d+1}\defeq\sigma_d(\BS^d)$. It is often useful to work with $\sdn\defeq\omega_{d+1}^{-1}\sigma_d$, the normalized spherical Lebesgue measure. The {\em intensity measure} of a point process $X$ on $\BS^d$ is  $\BE X(\cdot)$. (If $X$ is a simple point process, we can view $X$ as a
random collection of finitely many points in $\BS^d$.) If $X$ is an isotropic point process and $\BE X(\BS^d)<\infty$, then $\BE X(\cdot)=\gamma_s \sdn(\cdot)$ with some constant $\gamma_s\ge 0$, which is referred to as the {\em intensity} of $X$. Moreover, if $X$ is a (not necessarily isotropic) Poisson process, then for any Borel set $A\subset \BS^d$ the integer-valued random variable $X(A)$ is   Poisson distributed with Poisson parameter $\BE X(A)<\infty$. Henceforth, we always assume that the intensity measure  of a point process  on $\BS^d$ is finite and not the zero measure. Hence, in the isotropic case, the intensity $\gamma_s$ of $X$ is positive and finite.
From a (Poisson) point process $X$  on $\BS^d$ we obtain a hypersphere tessellation of $\BS^d$ by partitioning the unit sphere with the random hyperspheres $\{\BS^d\cap x^\perp:x\in X\}$, where $x^\perp$ denotes the orthogonal complement of $x$  (see again Fig.~\ref{Fig1} for an illustration).

We aim to show that, under appropriate assumptions, the shape of the Crofton cell $Z_0$, given a fixed positive lower bound for its
spherical volume, converges to the shape of a spherical cap as $\gamma_s\rightarrow\infty$. In fact, as has been explained in our previous work in Euclidean space, such results on limit shapes follow from more general deviation results, which hold for any fixed intensity and quantify the deviation of the shape of $Z_0$ from the shape of a spherical cap.
A functional $\vartheta:\overline{\cK}^d_s\rightarrow[0,\infty)$ is called a {\em deviation functional
for the class of spherical caps}, if it is continuous and if $\vartheta(K)=0$ holds, for some $K\in\overline{\cK}^d_s$ with
$\sigma_d(K)>0$, if and only if $K$ is a spherical cap (for precise definitions, see Section \ref{sec:genisop}). An example of such a deviation functional is the difference between spherical circumradius and the
spherical inradius of $K$. Another example, denoted by $\Delta_2^*$, measures the
deviation of the spherical radial function of $K$, with respect to a suitably chosen centre point, from its integral average in the $L^2$-sense.
A crucial tool in the probabilistic analysis of the asymptotic shape of the spherical Crofton cell of a   Poisson hypersphere
tessellation is the use of general  inequalities of isoperimetric type and stability improvements thereof, for {\em size and hitting functionals} in spherical space (explicit definitions are given in Section \ref{sec:genisop}). While isoperimetric results for a variety of geometric functionals in Euclidean space have been the subject of numerous investigations (see \cite{Schneider14} for a detailed and profound exposition of geometric inequalities and \cite{BoeHug1,BoeHug2,BoeHugFod,BoeHugFod2,BoeHugFod3} for some recent stability results, applications and further references), much less is known in the spherical setting.

In the following, we say that $\Phi$ is the hitting functional associated with (or determined by) an isotropic Poisson hypersphere
tessellation (or with an isotropic hypersphere process) if $\Phi=2U_1$ on $\overline{K}^d_s$, where $U_1$ is given by
$$
U_1(K)\defeq \frac{1}{2}\int_{\BS^{d}}\mathbf{1}\{ x^\bot\cap K\neq\emptyset\}\, \sdn(\dint x),\quad K\in\overline{K}^d_s.
$$
The number $\Phi(K)$ can be interpreted as the invariant measure of all hyperspheres hitting $K$ (see \eqref{eq:justhitting} for an explanation of this terminology). Thus it is a spherical
analogue of the Euclidean mean width functional. To resolve Kendall's problem in spherical space, with the spherical volume as the
size functional, we need a stability improvement of a geometric inequality of isoperimetric type involving spherical mean width $U_1$ and volume $\sigma_d$.
In \cite{SG}, the following inequality \eqref{basicsphere} is shown. It can be interpreted as a spherical version of the classical Urysohn inequality. The latter  provides a lower bound for the mean width functional (a multiple of the first intrinsic volume $V_1$)  in terms of the volume functional $V_d$ in Euclidean space $\R^d$, that is,
$V_1(K)\ge \left(V_d(K)/V_d(B)\right)^{1/d}\, V_1(B)$, where $B\subset\R^d$ is a Euclidean ball and $K\subset\R^d$ is a convex body
(see  \cite[Corollary 3.2]{HW2020} or \cite[p.~382]{Schneider14}). Equality holds if and only if $K$ is a ball. Stated in this form, the homogeneity of the involved functionals is crucial.
An equivalent form of the inequality states that $V_1(K)\ge  V_1(B)$, whenever $K$ is a convex body and $B$ is a Euclidean ball of equal
volume. In spherical
space, a corresponding result can be formulated as follows.

Let $ K\in\overline{K}^d_s$. If $C$ is a spherical cap with $\sigma_d(K)=\sigma_d(C)$, then
\begin{equation}\label{basicsphere}
U_1(K)\geq U_1(C),
\end{equation}
and equality holds if and only if $K$ is a spherical cap. Two proofs are provided in \cite{SG}. The second proof
exhibits an interesting connection to the Euclidean Blaschke--Santal\'o inequality \cite[Chap.~10.7]{Schneider14}. Furthermore, this proof can be
strengthened to yield the following
more general stability estimate.   For a spherical cap $C$, we denote its radius by $\alpha_C$. We write $\Delta_2$ for the restriction of $\Delta_2^*$ to the set $\cK^d_s$ of proper convex bodies of the unit sphere (explicit definitions are given in Section \ref{sec:geostab}).

\medskip

\textbf{Theorem A.} {\em Let $K\in {\cK}^d_s$, and let $C\subset \BS^d$ be a spherical cap
with $\sigma_d(K)=\sigma_d(C)>0$.
Then there is a constant $\beta_\circ >0$ such that
$$
U_1(K)\ge \left(1+ \beta_\circ \,\Delta_2(K)^2\right)U_1(C),
$$
where $ \beta_\circ $ depends on $d,\alpha_C$.}

\medskip

A version of this result is stated as Theorem \ref{stability}, where the dependence of $\beta_\circ$ on  $d$ and $\alpha_C$  is made explicit.

In Section \ref{sec:probabin}, we prove the following result, which is based on Theorem A. It provides not only the asymptotic shape (which is a spherical cap) of the spherical Crofton cell
given a lower bound for its volume, but also gives deviation inequalities for fixed intensities.

\medskip

\textbf{Theorem B.} {\em
Let $Z_0$ be the Crofton cell of the hypersphere  tessellation derived from an isotropic Poisson point process $X$ on $\BS^d$ with intensity $\gamma_s\in (0,\infty)$.
If $0<a<\omega_{d+1}/2$ and $ \varepsilon\in (0,1]$, then there are constants $c^*,\beta^* >0 $ such that
$$
\BP\left(\Delta_2^*(Z_0)\ge \varepsilon \mid \sigma_{d}(Z_0)\geq a\right)
\leq {c}^*\,\exp\left(-\beta^* \, \varepsilon^{2(d+1)}\,\gamma_s \right),
$$
where the constant ${c}^*$ depends  on $a,\varepsilon,d$ and $\beta^* $ depends on $a, d$.}

\medskip

For both, Theorem A and Theorem B, similar results with different stability exponents are obtained in
Section \ref{sec:geostab} and Section \ref{sec:probabin}, respectively, for
the inradius as the size functional and a suitably chosen deviation functional. In the case of a
general size functional and an associated deviation functional, an isoperimetric inequality and a corresponding
stability result are obtained in Section \ref{sec:genisop}. This finally leads to a general resolution of Kendall's problem, but
without an explicit bound on the stability order with respect to $\varepsilon$ as in Theorem B (see Theorem \ref{thm:genresult}).

As a consequence of
our approach, we also obtain the asymptotic distribution of the size functional of the Crofton  cell of an isotropic   Poisson hypersphere tessellation as the intensity goes to infinity.

\medskip

\textbf{Theorem C.} {\em
Let $Z_0$ be the Crofton cell of the hypersphere  tessellation derived from an isotropic Poisson point process $X$ on $\BS^d$ with intensity $\gamma_s\in (0,\infty)$.
Let $\Sigma$ be a general increasing and rotation invariant size functional. If
 $a>0$ is such that $\Sigma^{-1}([a,\infty))\neq\emptyset$, then
$$
\lim_{\gamma_s\to\infty}\gamma_s^{-1}\, \ln \BP\left(\Sigma(Z_0)\ge a\right) =-\tau(\Phi,\Sigma,a),
$$
where $\tau(\Phi,\Sigma,a)$ is the isoperimetric constant (introduced in Section \ref{sec:genisop}),  associated with the hitting functional
$\Phi=2 U_1$, the size functional $\Sigma$ and the threshold $a$.}

\medskip

In particular, this shows that the probability
$\BP(\Sigma(Z_0)\geq a)$ decays exponentially fast as $\gamma_s\to\infty$ (see Section \ref{sec:aymdistr}).

\medskip

Results similar to Theorem B and Theorem C (and to results stated in the following)
can also be obtained for tessellations derived from a binomial process
of size $N\ge d+1$, where the deterministic number $N$ (of points or subspaces)  replaces the intensity.
Since the arguments are similar, and preliminary versions of such results are contained in \cite{ARThesis},
we do not provide further details here.

\medskip

After investigating Crofton cells, it is a natural next step to look at typical cells.
Thus, in Section \ref{sec:typcells} we consider typical objects in spherical space.
Since $\BS^d$ is a homogeneous $\SO(d+1)$-space (see \cite[p. 584]{SW2008}), we could use the framework of
random measures on homogeneous spaces (see \cite{Last} and \cite{RZ}). Instead we will provide a more direct approach which yields
some additional insights. We start by recalling briefly the
Euclidean framework which has  been thoroughly studied.
A  process of compact convex particles in
Euclidean space $\R^d$ is a point process on the space $\cK^d$ of nonempty compact convex subsets of $\R^{d}$ (see \cite[Chap.~4.1]{SW2008}). If $\zeta$ is a (simple) stationary particle process (that is, its distribution is invariant under translations) with intensity $\gamma_{\zeta}$ and $c:\cK^{d}\rightarrow\R^{d}$ is a translation covariant centre function, then  a very intuitive representation for the distribution $\BQ$ of the {\em typical
particle} of $\zeta$ (see \cite[p. 106]{SW2008}) is
$$
\BQ(\cdot)=\frac{1}{\gamma_{\zeta}}\BE\sum_{K\in \zeta}\mathbf{1}\{K-c(K)\in\cdot\}\mathbf{1}\{c(K)\in[0,1]^d\},
$$
where  $\BQ$ is concentrated on sets having the Euclidean origin as their
centre. Here we implicitly use that $K-c(K)$ is the unique translate of $K$ whose centre is the origin.
In contrast, in spherical space there are infinitely many rotations $\varphi\in \SO(d+1)$ such that $\varphi \sfp=x$, for any fixed $x\in\BS^d$.
This is the reason why an additional randomization is naturally employed in the following definition of a typical particle of an isotropic  particle process
in spherical space. Let
$X'$ be a (simple) isotropic spherical particle process with intensity $\gamma_{X'}\in (0,\infty)$ (precise definitions are given in Section \ref{sec:typcells}).  Then the
distribution of the associated {\em typical particle} $Z$ can be defined by
$$
\BP(Z\in\cdot)=\frac{1}{\gamma_{X'}}\BE \sum_{K\in X'}\int_{\SO(d+1)}\mathbf{1}\{\varphi^{-1}K\in\cdot\}\
\kappa(c_s(K),\dint \varphi) ,
$$
where $c_s$ is a rotation covariant centre function, $\kappa$ is a probability kernel such that $\kappa(x,\cdot)$, for $x\in\BS^d$, is a probability
measure on $\SO(d+1)$, concentrated on the set $\{\varphi\in \SO(d+1): \varphi\sfp=x\}$. In addition, we set $c_s(\BS^d)=o$ (the zero vector of  $\R^{d+1}$) and choose $\kappa(o,\cdot)$ as the invariant probability measure on $\SO(d+1)$.
  Although $Z$ is not isotropic, its distribution
$\BP_Z$ is invariant under rotations fixing $\sfp$ and satisfies a disintegration result
for isotropic particle processes on $\BS^d$  (a Euclidean analogue can be found in \cite[Theorem 4.1.1]{SW2008}).
The disintegration result and the partial invariance together are characteristic for the distribution $\BP_Z$, as stated in Theorem \ref{thmdisintegration}.

In Section \ref{sec:typcellPois}, we interpret an isotropic tessellation $X'$ of $\BS^d$ as an isotropic particle process and use the
aforementioned disintegration result to obtain that $\BE\left[f(Z_0)\right]=\gamma_{X'}\BE\left[f(Z)\,\sdn(Z)\right]$,  for any measurable and rotation invariant function $f$ of the particles. Thus we relate the distributions of  the Crofton cell $Z_0$   and of the {\em typical cell} $Z$ of $X'$.
Using this relation, we transfer Theorem B and Theorem C to the typical cell of a Poisson hypersphere tessellation (the same can be done for the typical cell of the tessellation induced by a binomial hypersphere process of size $N\geq d+1$).

\medskip

Finally, we investigate the typical cell of an isotropic spherical Voronoi tessellation in Sections \ref{sec:sphPVT} and \ref{sec:typVor}.
For an isotropic spherical Poisson--Voronoi tessellations, we apply Mecke's characterization of Poisson processes
to show that the distribution of the typical cell of the Poisson--Voronoi
tessellation, induced by an isotropic Poisson process $X$ on $\BS^d$, is equal to the distribution of the Crofton cell
associated with a special
 Poisson hypersphere process $Y$. The hypersphere process $Y$ is the set of all hyperspheres having equal distance
to the spherical origin $\sfp$ and a point in $X$ and thus clearly $Y$ is not isotropic. This leads to a new functional $\widetilde{U}_{\sfp}$
on $\overline{\cK}_{\sfp}^d$, the set of all spherically convex bodies $K\in\overline{\cK}_{s}^d$ with $\sfp\in K\subset B_s(\sfp,\pi/2)$, defined by
$$
\widetilde{U}_{\sfp}(K)\defeq \int_{\BS^d}\mathbf{1}\{(x-\sfp)^\perp\cap K\neq\emptyset\}\, \sdn(\dint x),\quad K\in\overline{\cK}_{\sfp}^d.
$$
In this setting, we
measure the size with the centred spherical inradius defined by
$$
 r_{\sfp}(K)\defeq\max\{r\geq0: B_s(\sfp,r)\subset K\},\quad K\in\overline{\cK}_{\sfp}^d.
$$
Furthermore, let
$$
R_{\sfp}(K)\defeq\min\{r\geq0: K\subset B_s(\sfp,r)\},\quad K\in\overline{\cK}_{\sfp}^d,
$$
denote the centred spherical circumradius, and define a deviation functional $\vartheta_{\sfp}$
(for the class of spherical caps with centre $\sfp$) by
$$
\vartheta_{\sfp}(K)\defeq R_{\sfp}(K)-r_{\sfp}(K),\quad K\in\overline{\cK}_{\sfp}^d.
$$

Section \ref{subsec:ex3} is devoted to the following extremal and stability result for the size functional $\widetilde{U}_\sfp$.

\medskip

\textbf{Theorem D.} {\em If $a\in(0,\pi/2)$ and  $K\in\overline{\cK}_{\sfp}^d$ with
$r_{\sfp}(K)\geq a$,
%,and let $C\defeq B_s(\sfp,a)$.
then
\begin{equation*}
\widetilde{U}_{\sfp}(K)\geq\widetilde{U}_{\sfp}(B_s(\sfp,a))=\sdn(B_s(\sfp,2a)),
\end{equation*}
with equality if and only if $K=B_s(\sfp,a)$.
Furthermore, if $\vartheta_{\sfp}(K)\geq\varepsilon\in(0,1]$, then
\begin{equation*}
\widetilde{U}_{\sfp}(K)\geq \left(1+\tilde \beta_\circ\,\varepsilon^{\frac{d+1}{2}}\right)\widetilde{U}_{\sfp}(B_s(\sfp,a)),
\end{equation*}
where the constant $\tilde\beta_\circ$ depends on $a,d$.}

\medskip

Finally, in Section \ref{sec:sphPVT} we obtain an asymptotic result for the
typical cell $\ZV$ of the isotropic spherical Poisson--Voronoi tessellation induced by an isotropic Poisson process $X$ on $\BS^d$
with intensity $\gamma_s$, which is based on Theorem D.

\medskip

\textbf{Theorem E.} {\em  Let $X$ be an isotropic Poisson
process on $\BS^d$ with intensity $\gamma_s\in (0,\infty)$. If $a\in (0,\pi/2) $ and $\varepsilon\in(0,1]$, then
\begin{equation*}
\BP(\vartheta_{\sfp}(\ZV) \geq\varepsilon\mid r_{\sfp}(Z)\ge a)\leq
\tilde{c}\,\exp\left(- \tilde\beta \,\varepsilon^{\frac{d+1}{2}}\, \gamma_s\right),
\end{equation*}
where the constant $\tilde c>0$ depends  on $a,d,\varepsilon$ and  the constant $\tilde\beta >0$ depends  on $a,d$.}

\medskip

Asymptotic distributions of the spherical inradius or of general size functionals can be derived by
similar arguments as in the case of the Crofton cell.

\section{Hypersphere tessellations and the Crofton cell}\label{general1}
We work in the Euclidean space $\R^{d+1}$, $d\ge 2$, with scalar product $\langle\cdot\,,\cdot\rangle$ and norm $\|\cdot\|$.
The $d$-dimensional unit sphere is $\BS^d\defeq\{x\in\R^{d+1}:\|x\|=1\}$. The hyperspheres of $\BS^d$ are of the form $\BS_x\defeq\BS^d\cap x^\perp$, for some $x\in\R^{d+1}\setminus\{o\}$, where $x^\perp$ is the $d$-dimensional linear subspace totally orthogonal to $x$ and $o\in \R^{d+1}$ is the Euclidean origin (zero vector). The set of all hyperspheres of the unit sphere is $\BS^d_{d-1}\defeq\{\BS_x :x\in \BS^d\}$. For the interior, the closure and the boundary of a set $A\subset\BS^d$  (with respect to $\BS^d$) we write $\interior(A)$, $\cl(A)$ and $\partial A$.

The intrinsic (geodesic) distance of points $x,y\in\BS^d$ is denoted by $d_s(x,y)$ and given by $d_s(x,y)=2\arcsin(\|x-y\|/2)$. If $d_s(x,y)<\pi$, then the unique geodesic segment connecting $x$ and $y$ is denoted by $[x,y]$. For $\emptyset\neq A\subset\BS^d$ and $x\in\BS^d$, we write $d_s(A,x)\defeq \min\{d_s(a,x):a\in A\}$ for the distance of $x$ from $A$, and we set $A_\delta\defeq \{x\in \BS^d:d_s(A,x)\le \delta\}$ for $\delta\in [0,\pi]$. For compact sets $\emptyset\neq C,C'\subset\BS^d$ the spherical Hausdorff distance satisfies $\delta_s(C,C')=\min\{\varepsilon\in[0,\pi]: C\subset C'_\varepsilon,C'\subset C_\varepsilon\}$.

   The spherical convex hull of $A\subset\BS^d$ is denoted by $\conv(A)$ and obtained as the intersection of the positive hull of $A$ with $\BS^d$. In particular, for $x,y\in\BS^d$ with $d_s(x,y)<\pi$, we have $\conv(\{x,y\})=[x,y]$. Moreover, for $b\in\BS^d$ we write $\textrm{T}_b\BS^d$ for the tangent space of $\BS^d$ at $b$ (which can be identified with the linear subspace orthogonal to $b$) and  $\exp_b:\textrm{T}_b\BS^d\to\BS^d$ for the exponential map with base point $b$, in the sense of Riemannian geometry. The group of proper (orientation-preserving) rotations, which acts continuously and transitively on $\BS^d$, is denoted by $\SO(d+1)$ or $\Ih$ (for short).

Let $X$ be an isotropic Poisson process on $\BS^{d}$. We can view $X$ as a simple finite point process (a finite random counting measure)
on $\BS^d$ or as a finite random collection of points on $\BS^d$.
Since the spherical Lebesgue measure
$\sigma_{d}$ is (up to a constant) the only rotation invariant finite Borel measure on $\BS^{d}$, the intensity measure $\Theta=\BE X$ of $X$ satisfies $\Theta(\cdot)=\gamma_s
\,\sdn(\cdot)$, for some constant $\gamma_s\in[0,\infty)$.
Here we use that $\Theta$ is finite, since $X$ is finite and Poisson. In the following, we always
 assume that $\Theta\not\equiv 0$. Hence, $\BP(X(\BS^d)\ge 1)>0$, the {\em intensity} $\gamma_s$ of $X$ satisfies $\gamma_s\in (0,\infty)$ and
   $\BE X(\BS^{d})=\gamma_s$ is the expected number of points on the unit sphere.

Applying the measurable mapping
$h: \BS^{d} \rightarrow \BS^{d}_{d-1}$, $ x \mapsto \BS_x$,
to the points in $X$, we obtain the hypersphere process (the spherical hyperplane process) $\widetilde{X}\defeq h(X)$, where $h(X)$ denotes the image measure of $X$ under $h$. Clearly, $h(X)$ is again an isotropic and simple (a.s.) Poisson process (see \cite{LastPenrose,SW2008} for an introduction to Poisson processes in general spaces).
There is at least one hypersphere in (the support of) $\widetilde{X}$ if $\widetilde{X}\neq 0$, and in this case the hyperspheres of $\widetilde{X}$ partition $\BS^d$ into a finite collection of spherical polytopes with pairwise disjoint interiors. If $\widetilde{X}= 0$, then we obtain the empty tessellation, which is $\BS^d$ itself. Such a partition is referred to as a hypersphere (or spherical hyperplane) tessellation.

Recall that $\sfp$ denotes an arbitrary fixed point of the unit sphere (which we call the spherical origin or the pole). The {\em spherical Crofton cell} or {\em spherical zero cell} is the (a.s. uniquely determined) cell which contains
 $\sfp$ in its  interior. We will denote it by $Z_0$. Clearly, if $Z_0\neq\BS^d$, then $Z_0$ is almost surely  a spherical polytope, but $Z_0\in \mathcal{K}^d_s$ only if $X(\BS^d)\ge d+1$.

For $K\subset \BS^{d}$, we define
$\mathcal{H}_K\defeq\{L\in \BS^d_{d-1}: L\cap K\neq\emptyset\}$. Let $X$ be an isotropic Poisson point process with intensity $\gamma_s$ and $\widetilde{X}=h(X)$. If $K\in  \overline{\mathcal{K}}^d_s$, then $\mathcal{H}_K$ is Borel measurable, Campbell's theorem
(see \cite[Theorem 3.1.2]{SW2008}) yields that
\begin{equation}\label{eq:justhitting}
\BE \widetilde{X}(\mathcal{H}_K)
=\gamma_s\, 2\, U_1(K),\quad K\in  \overline{\mathcal{K}}^d_s,
\end{equation}
and the hitting functional $\Phi=2 U_1$ associated with $\widetilde{X}$ satisfies
$\Phi(C_0)=1$ for any hemisphere $C_0$ of $\BS^d$. To put the terminology into a broader context, let $X$ now
be a (not necessarily isotropic) point process on $\BS^d$ with finite and non-zero intensity measure $\BE X$. Then there are
a unique constant $\gamma_s> 0$ and a unique Borel probability measure $\kappa$ on $\BS^d$ such that
$\BE X=\gamma_s \, \kappa$. If $\widetilde X=h(X)$ is the induced hypersphere process, then
\begin{align*}
\BE\widetilde{X}(\mathcal{H}_K)=\gamma_s\int_{\BS^d}\mathbf{1}{\{\BS_x\cap K\neq\emptyset\}}\,\kappa(\dint x),
\end{align*}
in particular $\gamma_s=\BE \widetilde X(\BS^d_{d-1})$.
In this case, the hitting functional $\Phi$ associate with $\widetilde{X}$ is given by
$$
\Phi(K)=\int_{\BS^d}\mathbf{1}{\{\BS_x\cap K\neq\emptyset\}}\,\kappa(\dint x)=\frac{\BE\widetilde{X}(\mathcal{H}_K)}{\BE \widetilde X(\BS^d_{d-1})}.
$$
Hence,  $\Phi(K)$ is proportional to the mean number of hyperspheres in $\widetilde X$ hitting $K$.
If $X$ is isotropic, then $\kappa=\sdn$ and $\Phi=2U_1$. In the following, we write
\begin{equation*}
\mu(\cdot)\defeq \int_{\BS^{d}}\mathbf{1}{\{\BS_x\in\cdot\}}\, \sdn(\dint x)
\end{equation*}
for the rotation invariant  probability measure on  the Borel sets of $\BS^d_{d-1}$.

One of our principal aims in the following is to show that the Crofton cell $Z_0$, given a positive lower bound for its spherical volume, converges to a spherical cap as
$\gamma_s\rightarrow\infty$. This means that the conditional probability of $Z_0$ deviating from the shape of a spherical cap,
given $Z_0$ has spherical volume at least $a$, for some $a>0$, converges to $0$ as $\gamma_s\rightarrow\infty$. An explicit weak convergence result is stated and proved in Section \ref{sec:limshape}. More generally,
 we will quantify the deviation of $Z_0$ from a spherical cap.

Sections \ref{sec:genisop} - \ref{sec:approx} provide several geometric key results which are of interest in their
own right and are needed for the investigation of Kendall's problem in spherical space.
The latter is treated in Sections \ref{sec:probabin} - \ref{sec:sphPVT}. The results on typical cells of isotropic Poisson hypersphere processes, which are prepared by some general results on typical cells of isotropic particle processes, are covered in Section  \ref{sec:typcells}. Kendall's problem for the typical Voronoi cell of an isotropic Poisson process is treated in Section \ref{sec:sphPVT}.
Finally, some structural information on the typical Voronoi cell of
an isotropic point process on $\BS^d$ is contained in Section \ref{sec:typVor}.

\section{A general framework for isoperimetry}
\label{sec:genisop}

The probabilistic deviation results outlined in the introduction are based on geometric inequalities and related stability results.
Similarly as in Euclidean space (see, e.g., \cite[Chap.~12.2]{HS24}), we will describe a very general setting for stability results of isoperimetric type,
which are then applied to the solution of Kendall's problem in spherical space.

The main ingredients in our analysis are a hitting functional $\Phi$, a size functional $\Sigma$, and a deviation measure $\vartheta$.
In general, by a {\em hitting functional} we mean a map $\Phi:\overline{\cK}^{d}_s \to [0,\infty)$ that is continuous and such that $\Phi(K)=0$, for some $K\in
\overline{\cK}^{d}_s$, if and only if $K$ is a one-point  set. Here and in the following, the continuity on spherically convex bodies refers to the (spherical)  Hausdorff metric $\delta_s$. The first main example of a hitting functional is the hitting functional
$\Phi=2 U_1$ of an isotropic hypersphere process, as introduced in \eqref{eq:justhitting}, which is proportional to the spherical mean width functional. In this case, $\Phi$ is also
increasing and rotation invariant. Another example arises in the study of Voronoi tessellations as the hitting functional of a non-isotropic
hypersphere process and will be introduced later (on a restricted domain). Further examples are obtained by replacing $\sdn$ in the definition of $2U_1$
with a probability measure that is absolutely continuous with positive density with respect to $\sdn$.

A {\em size functional} is a continuous map $\Sigma:\overline{\cK}^{d}_s \to [0,\infty)$ which satisfies $\Sigma\not\equiv 0$ and $\Sigma(\{e\})=0$ for all
$e\in \mathbb{S}^d$. For the derivation of deviation inequalities and asymptotic distributions,
we will also assume that $\Sigma$ is increasing with respect to set inclusion. Examples of size functionals are volume, surface area, any of the functionals
$U_1,\ldots,U_{d-1}$ (see \cite[(6.62)]{SW2008}), inradius, diameter or width in a fixed direction. Note that the spherical intrinsic volumes $V_i$ (see \cite[Chap.~6.5]{SW2008} for an introduction) are not increasing
in general \cite[p.~262]{SW2008}.

Let $\Phi,\Sigma$ be fixed. Let $a>0$ be such that $\{K\in \overline{\cK}^{d}_s:\Sigma(K)\ge a\}=\Sigma^{-1}([a,\infty))\neq \emptyset$. We define the isoperimetric constant
$$
\tau(\Phi,\Sigma,a)\defeq\min\{\Phi(K):K\in \overline{\cK}^{d}_s,\Sigma(K)\ge a\}.
$$
Note that the minimum is attained, since $\Sigma^{-1}([a,\infty))$ is nonempty and compact (as $\overline{\cK}^{d}_s$ is compact and $\Sigma$ is continuous) and $\Phi$ is continuous. Moreover, we have $\tau(\Phi,\Sigma,a)>0$ since otherwise there would be some $K_0\in\overline{\cK}^{d}_s$ with $\Phi(K_0)=0$ and $\Sigma(K_0)\ge a>0$. But then $K_0=\{e_0\}$, for some $e_0\in\mathbb{S}^d$, and $0=\Sigma(\{e_0\})\ge a>0$, a contradiction.

If $a>0$ and $\Sigma^{-1}([a,\infty))\neq \emptyset$, then
$$\mathcal{E}(\Phi,\Sigma,a)\defeq\{K\in  \overline{\cK}^{d}_s:\Sigma(K)\ge a \text{ and }
\Phi(K)=\tau(\Phi,\Sigma,a)\}$$
is the nonempty set of {\em extremal bodies} associated with $\Phi$, $\Sigma$ and $a$. Clearly,
if $\Phi$ and $\Sigma$ are rotation invariant, then so is the class $\mathcal{E}(\Phi,\Sigma,a)$.
In the following sections, we simply write $\tau(a)$ and $\mathcal{E}(a)$ if $\Phi,\Sigma$ are clear from the context.

Finally, a continuous functional $\vartheta:\Sigma^{-1}([a,\infty))\to [0,\infty)$ such that
$\vartheta(K)=0$ if and only if $K\in\mathcal{E}(\Phi,\Sigma,a)$ is called a {\em deviation functional} for $\Phi,\Sigma,a$. A general,
canonical example is provided by
\begin{equation}\label{devfunc}
\vartheta(K)=\frac{\Phi(K)}{\tau(\Phi,\Sigma,a)}-1,\quad K\in \Sigma^{-1}([a,\infty)).
\end{equation}
For specific choices of $\Phi$ and $\Sigma$, other choices of deviation functionals will be more natural.
In particular, $\vartheta$ as given in \eqref{devfunc} is rotation invariant if this is true for $\Phi$ and $\Sigma$.

 \begin{proposition}\label{propneu1}
 Let $\Phi$ be a hitting functional, $\Sigma$  a size functional, and let $a>0$ be such that $\Sigma^{-1}([a,\infty))\neq \emptyset$. Let $\vartheta$ be a deviation functional for $\Phi,\Sigma,a$. Then there is an increasing (stability) function $f_a=f_{\Phi,\Sigma,a,\vartheta}:[0,\infty)\to [0,1]$ with $f_a(0)=0$, $f_a(t)>0$ for $t>0$ and such that
 $$
 \Phi(K)\ge \left(1+f_a(\varepsilon)\right)\tau(\Phi,\Sigma,a)
 $$
 for all $K\in \Sigma^{-1}([a,\infty))$ with $\vartheta(K)\ge \varepsilon\ge 0$.
 \end{proposition}

\begin{proof}
Under the assumptions of the proposition, we consider
$$
 \overline{\cK}^{d}_s(\Phi,\Sigma,a,\vartheta,\varepsilon)\defeq\{K\in \Sigma^{-1}([a,\infty)):\vartheta(K)\ge \varepsilon\}.
$$
Then $\overline{\cK}^{d}_s(\Phi,\Sigma,a,\vartheta,\varepsilon)\subset \overline{\cK}^{d}_s$ is compact. We can assume that this set is nonempty, since otherwise we simply define $f_a(\varepsilon)\defeq 1$ for any such $\varepsilon>0$. If $\overline{\cK}^{d}_s(\Phi,\Sigma,a,\vartheta,\varepsilon)\neq\emptyset$, we have
$$
\tau(\Phi,\Sigma,a,\vartheta,\varepsilon)\defeq\min\{\Phi(K):K\in \overline{\cK}^{d}_s(\Phi,\Sigma,a,\vartheta,\varepsilon)\}\ge \tau(\Phi,\Sigma,a).
$$
Assume that $\tau(\Phi,\Sigma,a,\vartheta,\varepsilon)=\tau(\Phi,\Sigma,a)$. Then there is some $K\in \overline{\cK}^{d}_s$ with $\Sigma (K)\ge a$, $\vartheta(K)\ge \varepsilon$ and $\Phi(K)=\tau(\Phi,\Sigma,a,\vartheta,\varepsilon)=\tau(\Phi,\Sigma,a)$, hence $K\in \mathcal{E}(\Phi,\Sigma,a)$.
By definition of $\vartheta$, we get $\vartheta(K)=0$, a contradiction. Thus we obtain $
\tau(\Phi,\Sigma,a,\vartheta,\varepsilon)> \tau(\Phi,\Sigma,a)
$
and define
$$
g_a(\varepsilon)\defeq\frac{\tau(\Phi,\Sigma,a,\vartheta,\varepsilon)}{\tau(\Phi,\Sigma,a)}-1>0
$$
and $f_a(\varepsilon)\defeq\min\{g_a(\varepsilon),1\}$. Then, for $K\in \Sigma^{-1}([a,\infty))$ with $\vartheta(K)\ge \varepsilon>0$, we have
$$
\Phi(K)\ge \tau(\Phi,\Sigma,a,\vartheta,\varepsilon)=(1+g_a(\varepsilon))\tau(\Phi,\Sigma,a)\ge (1+f_a(\varepsilon))\,\tau(\Phi,\Sigma,a),
$$
which proves the assertion.
\end{proof}

In the following section, we will provide specific versions of stability results with the
following choices of functionals.  Explicit definitions will be given in Section \ref{sec:geostab}.
In these specific situations, the domain of the functional has to be adjusted. Instead of using axiomatic
properties of functionals and deviation measures we will argue in a more direct way.

\begin{example}\label{ex1}
{\rm Let $\Sigma=\sigma_d$, $\Phi=2U_1$, $\vartheta=\Delta_2^*$ (an $L_2$-distance involving  spherical radial functions to be introduced in Section \ref{sec:geostab}).
Continuity of $\sigma_d$ on $\overline{\mathcal{K}}^d_s$ with respect to the Hausdorff metric is easy to see (\cite[p. 11]{Glas} or \cite[Theorem 6.5.2]{SW2008}), continuity of $U_1$ follows from \cite[(6.63) and Theorem 6.5.2]{SW2008}. A
direct argument as in the case of the functional $\widetilde{U}_\sfp$ (see Section \ref{sec:geostab}) can also be given. The definition
of the deviation measures $\Delta_2$ on $\mathcal{K}^d_s$ and $\Delta_2^*$ on $\overline{\cK}^d_s$ will be given in Section \ref{sec:geostab}. In particular, we will show that $\Delta_2$ and $\Delta_2^*$ are measurable (see Lemma \ref{eq:Delta2} and Remark \ref{rem2}), but continuity will not be required.}
\end{example}

\begin{example}\label{ex2}
{\rm Let $\Sigma=r$ (inradius), $\Phi=2U_1$, $\vartheta=R-r$, where $R$ is the circumradius (see \cite[Part 2, Proposition 2.7]{BH1999} or \cite{Dek} for basic properties of $R$). In contrast to the Euclidean case, the inball is uniquely determined if the inradius is positive (see \cite[Lemma 1]{White67}). Moreover, $r$ is a continuous functional on $\overline{\mathcal{K}}^d_s$. To verify this, let $K,K_i\in
    \overline{\mathcal{K}}^d_s$, $i\in\N$, with $K_i\to K$ as $i\to\infty$. It suffices to show that whenever
    $r(K_i)\to r_0\ge 0$, for some $r_0\ge 0$, as $i\to\infty$, then it follows that $r_0=r(K)$. For $i\in\N$ there is some
    $x_i\in K_i$ such that $B_s(x_i,r(K_i))\subset K_i$. There are a subsequence $(x_{i_j})_{j\in\N}$ and some $x_0\in \BS^d$
    such that $x_{i_j}\to x_0$ as $j\to\infty$. Thus it follows from \cite[Theorems 12.2.2 and 12.3.2]{SW2008} that $B_s(x_0,r_0)\subset K$. This yields that $r(K)\ge r_0$.  Assume that $r_0<r(K)\le \pi/2$. We can choose $\varepsilon>0$ such that $r_0+2\varepsilon \le r(K)$ for $i\ge i_0$. For $i\ge i_1\ge i_0$, we have $K\subset (K_i)_\epsilon$. Moreover, there is some $\bar{x}_0 \in\BS^d$ such that $B_s(\bar x_0,r_0+2\varepsilon)\subset K\subset (K_i)_\epsilon$ for $i\ge i_1$. Hence, by Lemma \ref{lemcancel}, $B_s(\bar x_0,r_0+\varepsilon)\subset K_i$ for $i\ge i_1$, which yields $r(K_i)\ge r_0+\varepsilon$ for $i\ge i_1$, and therefore $r_0\ge r_0+\varepsilon$, a contradiction. This finally shows that $r_0= r(K)$.

\begin{lemma}\label{lemcancel}
Let $A,B\in \overline{\mathcal{K}}^d_s$ and $\delta\in (0,\pi/2)$. If $A_\delta\subset B_\delta$, then $A\subset B$.
\end{lemma}

\begin{proof}
Suppose that $x\in A\setminus B$. Then $x\in A\subset A_\delta\subset B_\delta$, and hence $0<d_s(B,x)=:\rho\le \delta< \pi/2$.
 There are uniquely determined $b\in B$ and $u\in \textrm{T}_b\BS^d$ such that $d_s(B,x)=d_s(b,x)=\rho\in (0,\delta]$,  $\exp_b(\rho u)=x$ and $d_s(B,\exp_b(t u))=d_s(b,\exp_b(t u))=t$ for all $t\ge 0$ such that $\exp_b(t u)\in B_\delta$. For  $y\defeq\exp_b(\delta u)\in B_\delta$ we thus obtain $d_s(B,y)=\delta$, hence $y\in\partial B_\delta$. On the other hand, $\delta'\defeq d_s(A,y) \le d_s(x,y)=\delta-\rho<\delta$. Therefore $y\in A_{\delta'}$ and $B_s(y,\delta-\delta')\subset A_\delta\subset B_\delta$, which shows that $y\notin\partial B_\delta$, a contradiction.
\end{proof}

  The proof of the fact that $R$ is continuous on $\overline{\mathcal{K}}^d_s$ is straightforward.
  }
\end{example}

\begin{example}\label{ex3}
{\rm Let $\Sigma=r_{\sfp}$, $\Phi=\widetilde{U}_{\sfp}$ (or a multiple thereof), $\vartheta_{\sfp}=R_{\sfp}-r_{\sfp}$, where ${\sfp}\in\mathbb{S}^d$ is a fixed point. Only spherically convex
bodies are considered that contain ${\sfp}$ and are contained in the closed hemisphere centered at ${\sfp}$. As in the introduction, $r_{\sfp},R_{\sfp}$ denote the centred inradius and the centred circumradius, respectively. The functional $\widetilde{U}_\sfp$ will be discussed further  in Section \ref{sec:geostab}. The continuity of $r_{\sfp}$ follows as for $r$.}
\end{example}
%
%
%-------------------------------------------------------------------------------------------------------------------------------
%
%
%-------------------------------------------------------------------------------------------------------------------------------
%
\section{Geometric inequalities and stability results}
\label{sec:geostab}

In this section, we consider stability results which specify the general setting described in the preceding section.

\subsection{Framework of Example \ref{ex1}}
\label{subsec:ex1}

In the following, we use the notation and some of the results from
\cite{SG}, specifically
$$
D(x)\defeq\int_0^x\sin^{d-1}(t)\, \dint t,\quad x\in [0,\pi/2],
$$
and
$$
h(y)\defeq\tan^d(D^{-1}(y)),\quad y\in \text{im}(D).
$$
For $e\in \BS^d$, let $T_e\defeq e+e^\perp$ and define the open halfspace
$\He\defeq\{x\in \R^{d+1}:\langle x,e\rangle>0\}$, whose closure is   $\Hp\defeq \{x\in \R^{d+1}:\langle x,e\rangle\ge 0\}$. The  map $\proj_e:\BS^d\cap \He\to T_e$
with $\proj_e(u)\defeq\langle e,u\rangle^{-1}u$ is the radial projection to the tangent plane of $\BS^d$ at $e$.

For $K\in {\mathcal{K}}^d_s$ with $\sigma_d(K)>0$, the spherical polar
 $K^*\defeq\{u\in \BS^d:\langle u,x\rangle\le 0 \text{ for  $x\in K$}\}$ of $K$ is again in $ {\mathcal{K}}^d_s$ and $\interior(K^*)\neq \emptyset$.
If $e\in -\interior(K^*)$, then  $K\subset   \He$ (but it is not guaranteed that $e$ is contained in $K$).
The map
$$
F_K:-\interior(K^*)\to (0,\infty),\quad e\mapsto F_K(e)\defeq\int_K\langle e,u\rangle^{-(d+1)}\, \sigma_d(\dint u) ,
$$
assigns to $e\in -\interior(K^*)$ the volume of $\proj_e(K)$ in $T_e$ (see \cite[Sec.~3]{SG}). Since $F_K$ is continuous and $F_K(e)\to\infty$ as $e\to-\partial K^*$ (see \cite[Sec.~3]{SG}), the function $F_K$ on $-\interior(K^*)$
attains a minimum. The following lemma shows that $F_K$ is strictly geodesically convex, hence $F_K$ attains the minimum at a unique point $\fracc(K)\in -\text{int}(K^*)$. As usual, $F_K$ is called (strictly) geodesically convex, if for any unit speed geodesic $\gamma:I\to -\interior(K^*)$ (where $I\subset\R$ is an interval) the composition $F_K\circ\gamma:I\to (0,\infty)$ is (strictly) convex.

\begin{lemma}\label{Lem:sgc}
The map $F_K:-\interior(K^*)\to (0,\infty)$ is strictly geodesically convex.
\end{lemma}

\begin{proof}
Suppose that $e_1,e_2\in -\interior(K^*)$. If $u\in K$, then the geodesic segment $[e_1,e_2]$ is contained in $H^+(u)\cap \BS^d$. Let $u\in K$
be fixed, and let $\gamma$ be a unit speed geodesic from $e_1$ to $e_2$. Hence there is some $v\in \textrm{T}_{e_1}\BS^d$ with $\|v\|=1$ and some $t_0>0$ (independent of $u$) such that
$\gamma(t)\defeq \cos(t)e_1+\sin(t)v$, $t\in [0,t_0]$, satisfies $\gamma(0)=e_1$, $\gamma(t_0)=e_2$ and $\gamma(t)\in  H^+(u)\cap \BS^d$ for $t\in [0,t_0]$.
The map $\ell:[0,t_0]\to (0,\infty)$, $\ell(t)\defeq \langle \gamma(t),u\rangle$, is strictly concave on $[0,t_0]$, since $\ell''(t)=-\ell(t)<0$ for $t\in [0,t_0]$.
Since $(0,\infty)\to(0,\infty)$, $s\mapsto s^{-(d+1)}$, is strictly decreasing and convex, the composition $t\mapsto\ell(t)^{-(d+1)}$ is strictly convex on $[0,t_0]$.

This shows that, for each $u\in K$, the map $t\mapsto \langle\gamma(t),u\rangle^{-(d+1)} >0$ is strictly convex on $[0,t_0]$. Since $\sigma_d(K)>0$, it follows that $F\circ\gamma$ is strictly convex, which (by definition) is the asserted strict geodesic convexity of $F$.
\end{proof}

 The point $\fracc(K)$, which  in \cite{Besau24} is called the GHS-centre of $K$,  is the unique $e\in -\interior(K^*)$ such that
$$
\int_K\langle e,u\rangle^{-(d+2)}\langle v,u\rangle \, \sigma_d(\dint u)=0\quad \text{for }v\in\BS_e.
$$
Equivalently, it is the unique $e\in -\interior(K^*)$ such that
$$
\frac{1}{\mathcal{H}^d(\proj_e(K))}\int_{\proj_e(K)}x\,\mathcal{H}^d(\dint x)=e,
$$
which means that $e$ is the centre of mass of $\Pi_e(K)$ in $T_e$, where $\mathcal{H}^d$ denotes the $d$-dimensional Hausdorff measure. It follows that $\fracc(K)\in \interior(K)$, and thus $\fracc(K)\in \dom(K)\defeq\interior(K)\cap (-\interior(K^*))$. The preceding statements are covered by \cite[Sec.~3]{SG} or follow as straightforward consequences.
Moreover, if $K\in \overline{\mathcal{K}}^d_s$, then a separation argument shows that $\relint(K)\cap (-\relint(K^*))\neq \emptyset$, and hence $K\cap (-K^*)\neq\emptyset$.

For $e\in \dom(K)$, the positive and continuous function $\alpha_{K,e}:\mathbb{S}_e\to (0,\pi/2)$, defined by
\begin{equation}\label{radialrep}
\partial (\proj_e(K))=\left\{e+\tan(\alpha_{K,e}(u))u:u\in \mathbb{S}_e\right\},
\end{equation}
is the spherical radial function of $K$, and $\tan\circ\,\alpha_{K,e}$ is the radial function $\varrho(\proj_e(K),\cdot)$ of $\proj_e(K)$ with
respect to the origin $e$ in $T_e$ (as a function on $\BS_e$). We consider the set
$$
\mathcal{D}\defeq \left\{(K,e,u)\in \mathcal{K}^d_s\times\BS^d\times\BS^d:e\in\dom(K),u\in \BS_e\right\},
$$
which is (rotation) invariant with respect to an application of the same rotation $\sigma\in\SO(d+1)$ to all three arguments. For each $n\in\N$,
the set
$$
\mathcal{D}_n\defeq \left\{(K,e,u)\in \mathcal{K}^d_s\times\BS^d\times\BS^d:B_s(e,1/n)\subset K,B_s(-e,1/n)\subset K^*,u\in \BS_e\right\}
$$
is closed. To see this, one can use that $B_s(-e,1/n)\subset K^*$ if and only if $K\subset B_s(e,\pi/2-1/n)$ and the fact that inclusions are preserved under convergence of closed sets (see \cite[Theorem 12.2.2]{SW2008} and \cite[Theorem 12.3.2]{SW2008}, where the latter remains true in spherical space with the same proof). Therefore, $\mathcal{D}=\bigcup_{n\ge 1}\mathcal{D}_n$ is measurable as a countable union of closed sets.
Since $\alpha_{K,e}(u)=\alpha_{\sigma K,\sigma e}(\sigma u)$ for $\sigma\in \SO(d+1)$ and
$(K,e,u)\in \mathcal{D}$ and $\tan^{-1}(\varrho(\proj_{e}L,v))$ depends continuously on the argument $(L,v)\in\{(K,w)\in\mathcal{K}^d_s\times\BS_{e}:  e\in\dom(K)\}$, the map
$$\alpha:\mathcal{D}\to (0,\pi/2),\quad  (K,e,u)\mapsto \alpha_{K,e}(u),
$$
is continuous.

Let $\kappa_d$ denote the volume of the $d$-dimensional Euclidean unit ball, hence $d\kappa_d=\omega_d$. Using \cite[Lemma 6.5.1]{SW2008} and writing $\sden\defeq\omega_d^{-1}\sigma_{d-1}$ for the normalized spherical Lebesgue measure on any hypersphere of $\BS^d$, for $e\in \dom(K)$ we obtain
$$
\sigma_d(K)=\int_{\mathbb{S}_e}\int_0^{\alpha_{K,e}(u)}\sin^{d-1}(t)\, \dint t\, \sigma_{d-1}(\dint u),
$$
and thus
\begin{equation}\label{eq:meanvalue}
\frac{\sigma_d(K)}{d\kappa_d}=\int_{\mathbb{S}_e}D(\alpha_{K,e}(u))\, \sden(\dint u).
\end{equation}
If $C\subset\BS^d$ is a non-degenerate spherical cap contained in an open hemisphere, then there is a constant $\alpha_C\in (0,\pi/2)$ such that
$$
\frac{\sigma_d(C)}{d\kappa_d}=\int_0^{\alpha_C}\sin^{d-1}(t)\, \dint t=D(\alpha_C),\qquad
h\left(\frac{\sigma_d(C)}{d\kappa_d}\right)=\tan^d(\alpha_C).
$$
If $C^*\subset \mathbb{S}^d$ is the polar of $C$, then $\alpha_{C^*}+\alpha_C=\pi/2$.

For  $K\in {\mathcal{K}}^d_s$ and $e\in \dom(K)$, we define $\overline{D\circ \alpha_{K,e}}$ as the integral mean of $D\circ \alpha_{K,e}$ with respect to $\sden$ over
$\mathbb{S}_e$, as given in \eqref{eq:meanvalue}, and
$$
\Delta_2(K)\defeq \left\| D\circ \alpha_{K,\fracc(K)}-\overline{D\circ \alpha_{K,\fracc(K)}}\right\|_{L^2(\mathbb{S}_{\fracc(K)},\sden)}.
$$
Thus $\Delta_2(K)$ measures the deviation of the shape of $K$ from the shape of a spherical cap in the $L^2$ sense. Clearly, $\Delta_2(K)=0$ if and only if $K$ is a spherical cap.

For $K\in \cK^d_s$ and $e\in \dom(K)$, we define
$$
\underline{\alpha}_e(K)\defeq\min\{\alpha_{K,e}(u):u\in \mathbb{S}_e\},\qquad \overline{\alpha}_e(K)\defeq\max\{\alpha_{K,e}(u):u\in \mathbb{S}_e\}.
$$
If $e=\fracc(K)$, we simply write $\underline{\alpha}(K),\overline{\alpha}(K)$ for $\underline{\alpha}_{\fracc(K)}(K),\overline{\alpha}_{\fracc(K)}(K)$, respectively, and define
$$
\Delta_0(K)\defeq \overline{\alpha}(K)-\underline{\alpha}(K),
$$
which also measures the deviation of the shape of $K$ from the shape of a spherical cap.

\begin{lemma}\label{eq:Delta2}
The map $\Delta_2:\cK^d_s\to [0,\pi/2]$ is measurable and $\Delta_2(K)\le \Delta_0(K)$ for $K\in\cK^d_s$.
\end{lemma}

\begin{proof}
For $e\in \dom(K)$ we have
$$
\alpha_{K,e}(u)\in [\underline{\alpha}_e(K),\overline{\alpha}_e(K)],\quad u\in \mathbb{S}_e.
$$
Since $D,D'$ are strictly increasing, we get
\begin{align*}
\left|D(\alpha_{K,e}(u))-\frac{\sigma_d(K)}{d\kappa_d}\right|&\le D(\overline{\alpha}_e(K))-D(\underline{\alpha}_e(K))
\le D'(\overline{\alpha}_e(K))(\overline{\alpha}_e(K)-\underline{\alpha}_e(K))\\
&=\sin^{d-1}(\overline{\alpha}_e(K))(\overline{\alpha}_e(K)-\underline{\alpha}_e(K))
\le \overline{\alpha}_e(K)-\underline{\alpha}_e(K),
\end{align*}
and thus
\begin{equation}\label{delta-absch}
\Delta_2(K)\le \Delta_0(K)\le \pi/2.
\end{equation}
Next we show that $\Delta_2$ is measurable. Since $\alpha:\mathcal{D}\to [0,\pi/2]$ is continuous, the map
$$
\{(K,e)\in\cK^d_s\times\BS^d:e\in\dom(K)\}\to [0,\infty),\quad (K,e)\mapsto \int_{\BS_e}(D\circ\alpha_{K,e})^2(u)\, \sden(\dint u),
$$
is continuous, and therefore this is also true for the map
$$
\{(K,e)\in\cK^d_s\times\BS^d:e\in\dom(K)\}\to [0,\pi/2],\quad (K,e)\mapsto \left\| D\circ \alpha_{K,e}-\overline{D\circ \alpha_{K,e}}\right\|_{L^2(\mathbb{S}_e,\sden)}.
$$
The remaining assertion follows, once we have verified that the map $\fracc:\cK^d_s\to \BS^d$, $K\mapsto \fracc(K)$, is measurable. For this, we consider
$$
\mathcal{D}_n^* \defeq\left\{K\in\cK^d_s:B_s(\fracc(K),1/n)\subset K,B_s(-\fracc(K),1/n)\subset K^*\right\}
$$
for $n\in\N$. We claim that $\mathcal{D}_n^*$ is closed in $\overline{\mathcal{K}}^d_s$. For this, let $K_i\in \mathcal{D}_n^*$, $i\in\N$, and $K_i\to K_0\in \overline{\mathcal{K}}^d_s$. By assumption, $B_s(\fracc(K_i),1/n)\subset K_i$ and $K_i\subset B_s(\fracc(K_i),\pi/2-1/n)$ for $i\in\N$. Since $\BS^d$ is compact, $e_i\defeq \fracc(K_i)\to e_0$, for some $e_0\in\BS^d$, for $i\to\infty$ and $i\in I$, where $I\subset\N$ is an infinite subset. We conclude that $B_s(e_0,1/n)\subset K_0\subset B_s(e_0,\pi/2-1/n)$, in particular $K_0\in\cK^d_s$ and $\langle e_i,x\rangle \ge \sin(1/n)$ for $x\in K_i$ and $i\in\N_0$. Hence it follows from
$$
\int_{K_i}\langle \fracc(K_i),u\rangle^{-(d+2)}\langle v,u\rangle \, \sigma_d(\dint u)=0\quad \text{for }v\in\BS_{\fracc(K_i)},i\in I,
$$
and the dominated convergence theorem that
$$
\int_{K_0}\langle e_0,u\rangle^{-(d+2)}\langle v,u\rangle \, \sigma_d(\dint u)=0\quad \text{for }v\in\BS_{e_0}.
$$
Therefore we conclude that $e_0=\fracc(K_0)$ (recall the discussion after \eqref{Lem:sgc}). The argument shows that $\mathcal{D}_n^*$ is closed and the restriction of $\fracc$ to $\mathcal{D}_n^*$ is continuous and hence measurable. Since $\cK^d_s=\bigcup_{n\in\N}\mathcal{D}_n^*$, the assertion follows.
\end{proof}

\begin{theorem}\label{stability}
If $K\in{\cK}_s^d$ and  $C\subset\BS^d$ is a spherical cap with $\sigma_d(K)=\sigma_d(C)>0$, then
$$
U_1(K)\geq \left(1+{\beta}\,\Delta_2(K)^2\right)U_1(C),
$$
where
$$
{\beta}\defeq \min\left\{\frac{\frac{2}{d}\binom{d+1}{2}\sin^{d+1}(\alpha_C)\tan^{-2d}(\alpha_C)
}{1+\binom{d+1}{2}\left(\frac{\pi}{2}\right)^2\tan^{-d}
(\alpha_C)}, \frac{8}{\pi^2} \, D\left(\frac{\pi}{2}-\alpha_C\right)\right\}.
$$
\end{theorem}

\begin{proof}
Throughout the proof, let $K\in {\cK}_s^d$ be fixed. We simply write $\fracc$ for $\fracc(K)\in \dom(K)$, put $\alpha\defeq\alpha_{K,\fracc(K)}$ and restrict the domain of $D$ to $(0,\pi/2)$ so that $\text{im}(D)=D((0,\pi/2))=(0,D(\pi/2))$, since $D$ is strictly increasing and continuous. Since $U_1$ is rotation invariant,
we can assume $C$ to be centred at $\fracc$. Note that due to the assumptions we have $\alpha_C\in (0,\pi/2)$.

We continue to use the preceding notation and otherwise refer to \cite[Sec.~3]{SG}. Then
$x_0\defeq\overline{D\circ\alpha}=\sigma_d(K)/(d\kappa_d)= \sigma_d(C)/(d\kappa_d)\in \text{im}(D)$,
since $\alpha(u)\in(0,\pi/2)$, for $u\in \mathbb{S}_{\fracc}$, $\sden$ is a probability measure, and $D$ is strictly increasing and continuous.
For any $z\in \text{im}(D)$, we have
$$
h(z)-h(x_0)=h'(x_0)(z-x_0)+\frac{1}{2}h''(x_0+\theta(z-x_0))(z-x_0)^2,
$$
for some $\theta=\theta(x_0,z)\in (0,1)$. Since
$$
h'(y)=\frac{d}{\cos^{d+1}(D^{-1}(y))}\geq d
$$
and
$$
h''(y)=\frac{d(d+1)}{\cos^{d+2}(D^{-1}(y))}\frac{1}{\sin^{d-2}(D^{-1}(y))}\ge d(d+1),
$$
for $y\in \text{im}(D)$, we deduce that
\begin{equation}\label{eqstab1}
h(z)-h(x_0)\ge h'(x_0)(z-x_0)+\binom{d+1}{2} (z-x_0)^2.
\end{equation}
Moreover, it follows that the functions $h$ and $h'$ are strictly increasing.
Substituting $z=D(\alpha(u))$, $u\in \mathbb{S}_{\fracc}$, in \eqref{eqstab1}, and then integrating \eqref{eqstab1} with respect to $\sden$
over $\mathbb{S}_{\fracc}$, we obtain
$$
\int_{\mathbb{S}_{\fracc}}h(D(\alpha(u)))\, \sden(\dint u)-h\left(\frac{\sigma_d(K)}{d\kappa_d}\right)
\ge 0+\binom{d+1}{2}\Delta_2(K)^2.
$$
Using that
$$
h\left(\frac{\sigma_d(K)}{d\kappa_d}\right)=h\left(\frac{\sigma_d(C)}{d\kappa_d}\right)=\tan^d(\alpha_C)
,
$$
the radial representation \eqref{radialrep}, and hence
$$
\mathcal{H}^d(\proj_{\fracc}(K))=\frac{1}{d}\int_{\mathbb{S}_{\fracc}}\tan^d(\alpha(u))\, \sigma_{d-1}(\dint u)
=\frac{1}{d}\int_{\mathbb{S}_{\fracc}}h(D(\alpha(u)))\, \sigma_{d-1}(\dint u),
$$
we conclude
\begin{equation}\label{eqa}
\frac{\mathcal{H}^d(\proj_{\fracc}(K))}{\kappa_d}=\int_{\mathbb{S}_{\fracc}}h(D(\alpha(u)))\, \sden(\dint u)
\ge \left(1+\beta_1\Delta_2(K)^2\right)h\left(\frac{\sigma_d(K)}{d\kappa_d}\right),
\end{equation}
where $\beta_1\defeq\binom{d+1}{2}\tan^{-d}(\alpha_C)$. Next, we recall some relations from \cite{SG}. The assumption $\sigma_d(K)= \sigma_d(C)>0$ and the equality
case of \cite[(27)]{SG} imply that
\begin{equation}\label{eqn1}
h\left(\frac{\sigma_d(K)}{d\kappa_d}\right)=h\left(\frac{\sigma_d(C)}{d\kappa_d}\right)=\frac{\mathcal{H}^d(\proj_{\fracc}(C))}{\kappa_d},
\end{equation}
and the equality cases of \cite[(26)]{SG} and \cite[(30)]{SG} yield
\begin{equation}\label{eqn2}
h\left(\frac{\sigma_d(C^*)}{d\kappa_d}\right)=\frac{\kappa_d}{\mathcal{H}^d(\proj_{\fracc}(C))}.
\end{equation}
Hence, combination of \eqref{eqa}, \eqref{eqn1} and \eqref{eqn2} gives
\begin{equation}\label{eqn12neu}
\frac{\kappa_d}{\mathcal{H}^d(\proj_{\fracc}(K))}
\le  \frac{1}{1+\beta_1\Delta_2(K)^2}h\left(\frac{\sigma_d(C^*)}{d\kappa_d}\right).
\end{equation}
Now we use $(26)$ and $(30)$ from \cite{SG} for the first inequality and \eqref{eqn12neu} for the second inequality to obtain
\begin{align}
h\left(\frac{\sigma_d(K^*)}{d\kappa_d}\right)&\le \frac{\kappa_d}{\mathcal{H}^d(\proj_{\fracc}(K))} \le
%\le \frac{1}{1+\beta_1\Delta_2(K)^2}\frac{1}{h\left(\frac{\sigma_d(K)}{d\kappa_d}\right)}=
\frac{1}{1+\beta_1\Delta_2(K)^2}h\left(\frac{\sigma_d(C^*)}{d\kappa_d}\right)
%&= \frac{1}{1+\beta_1\Delta_2(K)^2}\frac{\kappa_d}{\mathcal{H}^d(\proj_{\fracc}(C))}\nonumber
%= \frac{1}{1+\beta_1\Delta_2(K)^2}h\left(\frac{\sigma_d(C^*)}{d\kappa_d}\right)\nonumber\\
%&=\left(1-\frac{\beta_1}{1+\beta_1\Delta_2(K)^2}\Delta_2(K)^2\right)h\left(\frac{\sigma_d(C^*)}{d\kappa_d}\right)
\nonumber\\
&\le \left(1-\beta_2\Delta_2(K)^2\right)h\left(\frac{\sigma_d(C^*)}{d\kappa_d}\right),\label{eqstab2}
\end{align}
where for the last step we used that
$$
\beta_2\defeq\frac{\beta_1}{1+\beta_1(\pi/2)^2}\le \frac{\beta_1}{1+\beta_1\Delta_2(K)^2},
$$
by Lemma \ref{eq:Delta2}. Since   $\alpha_{C^*}=\pi/2-\alpha_C$, we get
\begin{align}\label{eq:beta3}
\beta_3:&=\min\left\{\beta_2\frac{\tan^{-d}(\alpha_C)}{D(\frac{\pi}{2}-\alpha_C)}
\frac{\sin^{d+1}(\alpha_C)}{d},\left(\frac{2}
{\pi}\right)^2\right\}\nonumber\\
&\le \beta_2\frac{\tan^d(\alpha_{C^*})}{D(\alpha_{C^*})}\frac{\cos^{d+1}(\alpha_{C^*})}{d}
= \beta_2\frac{h\left(\frac{\sigma_d(C^*)}{d\kappa_d}\right)}{\frac{\sigma_d(C^*)}{d\kappa_d}}h'\left(\frac{\sigma_d(C^*)}
{d\kappa_d}\right)^{-1}.
\end{align}
The minimum in the definition of $\beta_3$ is taken  to ensure that $1-\beta_3\Delta_2(K)^2\ge 0$.
From \eqref{eq:beta3}, the mean value theorem and the fact that $h$ and $h'$ are increasing, we deduce that
\begin{align}\label{eqstab3}
&h\left(\frac{\sigma_d(C^*)}{d\kappa_d}\right)-h\left((1-\beta_3\Delta_2(K)^2)\frac{\sigma_d(C^*)}{d\kappa_d}\right)\nonumber\\
&\quad \le h'\left(\frac{\sigma_d(C^*)}{d\kappa_d}\right)\beta_3\Delta_2(K)^2\frac{\sigma_d(C^*)}{d\kappa_d}
 \le \beta_2\Delta_2(K)^2h\left(\frac{\sigma_d(C^*)}{d\kappa_d}\right) .
\end{align}
Combining \eqref{eqstab2} and \eqref{eqstab3}, we get
$$
h\left(\frac{\sigma_d(K^*)}{d\kappa_d}\right)\le h\left((1-\beta_3\Delta_2(K)^2)\frac{\sigma_d(C^*)}{d\kappa_d}\right),
$$
and hence
$$
 \sdn(K^*) \le (1-\beta_3\Delta_2(K)^2)\sdn(C^*) .
$$
Since
$$
\frac{1}{2}-U_1(K)=\sdn(K^*) ,\quad \frac{1}{2}-U_1(C)=\sdn(C^*)
$$
by \cite[(20)]{SG},
we deduce that
$$
U_1(K)\ge \frac{1}{2}\beta_3\Delta_2(K)^2+(1-\beta_3\Delta_2(K)^2)U_1(C).
$$
Finally, we use
$$
\frac{1}{2}-U_1(C)=\sdn(C^*) =D(\alpha_{C^*})\ge  2D(\alpha_{C^*})U_1(C),
$$
and therefore
$$
\frac{1}{2}\ge \left(1+2D\left(\frac{\pi}{2}-\alpha_C\right)\right)U_1(C),
$$
to get
$$
U_1(K)\ge\left[\beta_3\Delta_2(K)^2+2D\left(\frac{\pi}{2}-\alpha_C\right)\beta_3
\Delta_2(K)^2+1-\beta_3\Delta_2(K)^2\right]U_1(C),
$$
and thus
$$
U_1(K)\ge \left(1+2D\left(\frac{\pi}{2}-\alpha_C\right)\beta_3\Delta_2(K)^2\right)U_1(C),
$$
which yields the assertion with $\beta=2D ( {\pi}/{2}-\alpha_C )\beta_3$.
\end{proof}

\begin{remark}\label{rem1}
{\rm
It is easy to see that
$$
\beta\ge \min\left\{\frac{\sin^{d+1}(\alpha_C)}{\tan^{2d}(\alpha_C)+2d\tan^d(\alpha_C)},
0.4^d\left(\frac{\pi}{2}-\alpha_C\right)^d\right\}.
$$
}
\end{remark}

\begin{remark}\label{rem2}
{\rm For  $K\in \overline{\mathcal{K}}^d_s$ with $\sigma_d(K)>0$, we define $\Delta_2^*(K)\defeq 0$, if $K$ is a hemisphere, and otherwise we set $\Delta_2^*(K)\defeq\Delta_2((K\cap (-K))^\perp\cap K)$, where $\Delta_2$ is considered in $\BS^d\cap (K\cap (-K))^\perp$. Again $\Delta_2^*$ is measurable (since all involved operations such as taking intersections, orthogonal complements, dimension are measurable) and an extension of Theorem \ref{stability} is obtained as a corollary, in which $C$ has to be replaced by a spherical wedge of the appropriate dimension.
 }
\end{remark}

\subsection{Framework of Example \ref{ex2}}
\label{subsec:ex2}

For $K\in\overline{\cK}_s^d$ and $e\in K$, let $r_e(K)$ denote the (spherical) inradius and $R_e(K)$ the circumradius of $K$, with respect to $e$ as the centre of the insphere and the circumsphere, respectively. We consider the size
functional
$$
\Sigma_r(K)\defeq\max\{r_e(K):e\in K\}=r(K),
$$
which is the inradius $r(K)$ of $K$.
The following lemma shows that $R_e(K)\le \pi/2$ if $e\in K$ is (the unique point in $K$ provided that $\interior(K)\neq\emptyset$) such that $r_e(K)=\Sigma_r(K)$.

\begin{lemma}\label{maxinr}
Let $K\in\overline{\cK}_s^d$. If $e\in K$ is such that $r_e(K)=\Sigma_r(K)$, then
$K\subset B_s(e,\pi/2)$.
\end{lemma}

\begin{proof}
For the proof we can assume that $\interior(K)\neq\emptyset$ (otherwise we argue in a lower-dimensional great subsphere). If $r_e(K)=\pi/2$, then $K=B_s(e,r_e(K))$ and the assertion is clear. From now on we can assume that $0<r_e(K)<\pi/2$.
By assumption and a separation argument, applied in the hyperplane spanned by $\partial B_s(e,r_e(K))$, it follows that $\cos(r_e(K))e\in \text{conv}(\partial K\cap B_s(e,r_e(K)))$ (here the convex hull is taken in $\R^{d+1}$).
Hence there exist $k\in\N$,  $\lambda_1,\ldots,\lambda_k\in[0,1]$ with  $\sum_{i=1}^k\lambda_i=1$, and $a_1,\ldots,a_k\in\text{conv}(\partial K\cap B_s(e,r_e(K)))\subset\mathbb{S}^d$ such that
\begin{equation}\label{eq1.1}
\sum_{i=1}^k\lambda_i a_i=\cos(r_e(K))e.
\end{equation}
Furthermore, since $B_s(e,r_e(K))$ supports $K$ from inside (and hence the same is true for the convex cones generated by these sets), we have
\begin{equation}\label{eq1.2}
\langle a_i-(\cos(r_e(K)))^{-1} e,x\rangle \le 0, \quad x\in K,\,i=1,\ldots,k.
\end{equation}
From \eqref{eq1.2}, we get
$$
\left\langle \sum_{i=1}^k  \lambda_i a_i- (\cos(r_e(K)))^{-1} e,x\right\rangle \le 0, \quad x\in K.
$$
Hence,  \eqref{eq1.1} implies that $\left(\cos(r_e(K))-(\cos(r_e(K)))^{-1}\right)\langle e , x\rangle\le 0$, that is, $\langle e,x\rangle \ge 0$, for
all $x\in K$.
\end{proof}

\bigskip

In order to measure the deviation of the shape of $K\in\overline{\mathcal{K}}^d_s$ from a spherical shape, we  now use the deviation measure
$$
\vartheta_r(K)\defeq\min\{R_e(K)-r_e(K):e\in K,r_e(K)=\Sigma_r(K)\}.
$$
Note that if  $\interior(K)\neq\emptyset$, then the insphere of $K$ is uniquely determined (see Example \ref{ex2}) and  $\vartheta_r(K)=R_e(K)-r_e(K)$, where $e$ is the center of the insphere of $K$.
The obvious geometric inequality $U_1(K)\ge  U_1(B_s(e,a))$, for $K\in\overline{\mathcal{K}}^d_s$ with $\Sigma_r(K)\ge a$, will be improved by a stability result in Theorem \ref{stab2} . The following simple lemma will be useful.

\begin{lemma}\label{lem:UsBa}
If $e\in\BS^d$ and $a\in [0,\pi/2]$, then
\begin{equation}\label{eq:UsBa}
U_1(B_s(e,a))=\frac{\omega_d}{\omega_{d+1}}\int_0^a\cos^{d-1}(s)\, \dint s\in [0,1/2].
\end{equation}
Moreover,
\begin{equation}\label{eq:UsBa2}
\frac{\omega_d}{\omega_{d+1}}\left(\frac{2}{\pi}\right)^{d-1}a\cdot \max\{1,(\pi/2-a)^{d-1}\}\le U_1(B_s(e,a))\le \frac{\omega_d}{\omega_{d+1}}a.
\end{equation}
\end{lemma}

\begin{proof}
By definition, we have
\begin{align*}
U_1(B_s(e,a))&=\frac{1}{2\omega_{d+1}}\int_{\BS^d}\indi\{\BS_x\cap B_s(e,a)\neq \emptyset\}\, \sigma_d(\dint x)\\
&=\frac{\omega_d}{2\omega_{d+1}}\int_0^\pi\sin^{d-1}(t)\indi\left\{t-\frac{\pi}{2}\in [-a,a]\right\}\, \dint t\\
&=\frac{\omega_d}{2\omega_{d+1}}\int_{-a}^a\cos^{d-1}(s)\, \dint s,
%&=\frac{\omega_d}{\omega_{d+1}}\int_0^a\cos^{d-1}(s)\, \dint s.
\end{align*}
which yields \eqref{eq:UsBa}.
(For $a=\pi/2$, the last integral equals $\frac{\omega_{d+1}}{\omega_d}$.) The upper bound for $U_1(B_s(e,a))$ is now obvious. For the lower bound, we use Jensen's inequality to get (for $a\in (0,\pi/2]$)
\begin{equation}\label{eq:UsBa3}
\frac{1}{a}\int_0^a \cos^{d-1} (s)\, \dint s\ge \left(\frac{1}{a}\int_0^a\cos(s)\, \dint s\right)^{d-1}=\left(\frac{\sin(a)}{a}\right)^{d-1}\ge\left(\frac{2}{\pi}\right)^{d-1}.
\end{equation}
On the other hand, using $\cos(t)\ge \frac{2}{\pi}(\frac{\pi}{2}-t)$ for $t\in [0,\pi/2]$ and the mean value theorem, we get
\begin{align}\label{eq:UsBa4}
\int_0^a\cos^{d-1} (s)\, \dint s
&\ge \left(\frac{2}{\pi}\right)^{d-1}\int_0^a\left(\frac{\pi}{2}-s\right)^{d-1}\, \dint s
= \left(\frac{2}{\pi}\right)^{d-1}\frac{1}{d}\left[\left(\frac{\pi}{2}\right)^d-\left(\frac{\pi}{2}-a\right)^d\right]\nonumber\\
&\ge \left(\frac{2}{\pi}\right)^{d-1}\left(\frac{\pi}{2}-a\right)^{d-1}a.
\end{align}
Combination of \eqref{eq:UsBa3} and \eqref{eq:UsBa4}  yields the lower bound in \eqref{eq:UsBa2}.
\end{proof}

\begin{theorem}\label{stab2}
Let $K\in\overline{\mathcal{K}}^d_s$, $a>0$, and $\varepsilon\in(0,1]$.
If $\Sigma_r(K)\ge a$ and $\vartheta_r(K)\ge \varepsilon$, then
$$
U_1(K)\ge \left(1+ c_\circ(a,d) \, \varepsilon^{\frac{d+1}{2}} \right)  U_1(B_s(e,a)),
$$
where $c_\circ(a,d)\ge  (8d)^{-1}(3\pi^{-4})^{d-1}a^{d-2}\left(\tfrac{\pi}{2}-a\right)^{d-1}\ge 2\cdot 0.02^da^{d-2}\left(\tfrac{\pi}{2}-a\right)^{d-1}$.
\end{theorem}

\begin{proof}
Let $K\in\overline{\mathcal{K}}^d_s$ be such that $\Sigma_r(K)\ge a$ and $\vartheta_r(K)\ge \varepsilon$. Let $e\in K$ be such that $r_e(K)=\Sigma_r(K)$ so that $a+\varepsilon\le r_e(K)+\varepsilon\le R_e(K)\le \pi/2$, where Lemma \ref{maxinr} was used.
Then there is some $z_0\in K$ such that $\conv(B_s(e,a)\cup\{z_0\})\subset K$ and $d_s(e,z_0)=a+\varepsilon>0$. Hence,
writing
$$
A\defeq\{x\in \BS^d:\BS_x\cap B_s(e,a)=\emptyset, \BS_x\cap \conv(B_s(e,a)\cup\{z_0\})\neq\emptyset\},
$$
we have
\begin{align*}
U_1(K)&\ge U_1(\conv(B_s(e,a)\cup\{z_0\}))=U_1(B_s(e,a))+\frac{1}{2\omega_{d+1}}
\int_{\mathbb{S}^{d}}\mathbf{1}\{x\in A \}\, \sigma_d(\dint x).
\end{align*}
For given $e\in\mathbb{S}^d$, $u\in \BS_e$, $a>0$, and $\varepsilon\in (0,1]$, we define
$$
\delta(e,u,a,\varepsilon)\defeq\sigma_1([\conv(B_s(e,a)\cup\{z_0\})\setminus B_s(e,a)]\cap
\{\lambda \, e+\mu\, u:\lambda\in\R,\mu\ge 0\})
$$
and
$$
C(e,a,\varepsilon)\defeq\{u\in \mathbb{S}_e:\delta(e,u,a,\varepsilon)\ge \varepsilon/2\},
$$
where $\sigma_1$ is the $1$-dimensional spherical Lebesgue measure on geodesic arcs.
By symmetry, $C(e,a,\varepsilon)\subset\mathbb{S}_e$  is a spherical cap
with centre $(z_0-\langle z_0,e\rangle e)/\sqrt{1-\langle z_0,e\rangle^2}$. We will show below that the
spherical radius $\omega$ of this cap satisfies
\begin{equation}\label{eq:auxbound}
\tan(\omega)\ge \sqrt{2}\pi^{-2}\, a\sqrt{\varepsilon}.
\end{equation}
Once this is shown, it follows that
\begin{align}\label{eq:lowvolbd}
\sigma_{d-1}(C(e,a,\varepsilon))&\ge \kappa_{d-1}\sin^{d-1}(\omega)%=\frac{\omega_{d-1}}{d-1}\left(\frac{\tan(\omega)}{\sqrt{1+\tan^2(\omega)}}\right)^{d-1}\nonumber\\
\ge \frac{\omega_{d-1}}{d-1}\left(\pi^{-2}a\sqrt{\varepsilon}\right)^{d-1}\ge \frac{\omega_{d}}{2d}\left(\pi^{-2}a\sqrt{\varepsilon}\right)^{d-1},
\end{align}
where for the second inequality we used that
$$
\sin(\omega)=\frac{\tan(\omega)}{\sqrt{1+\tan^2(\omega)}}\ge \frac{1}{\sqrt{2}}\tan(\omega)\ge \pi^{-2}a\sqrt{\varepsilon},\quad\text{if }0\le \omega\le \pi/4,
$$
by \eqref{eq:auxbound},
$$
\sin(\omega)\ge \frac{1}{2}\sqrt{2}\ge  \pi^{-2}a\sqrt{\varepsilon},\quad\text{if }\pi/4\le \omega\le \pi/2.
$$
and
$$
\frac{\omega_{d-1}}{(d-1)\omega_d}\ge \frac{1}{2d},\quad d\ge 2.
$$
We have
\begin{align}\label{eq:lb2}
\int_{\mathbb{S}^d}\mathbf{1}\{x\in A\}\, \sigma_d(\dint x)&\ge \int_{C(e,a,\varepsilon)}\int_{\frac{\pi}{2}+a}^{\frac{\pi}{2}+a+\delta(e,u,a,\varepsilon)}
\sin^{d-1} (t)\, \dint t\, \sigma_{d-1}(\dint u)\nonumber\\
&\ge \int_{C(e,a,\varepsilon)}\int_{\frac{\pi}{2}+a}^{\frac{\pi}{2}+a+\frac{\varepsilon}{2}}
\sin^{d-1} (t)\, \dint t\, \sigma_{d-1}(\dint u)\nonumber\\
&=\int_{C(e,a,\varepsilon)}\int_{a}^{a+\frac{\varepsilon}{2}}
\cos^{d-1} (t)\, \dint t\, \sigma_{d-1}(\dint u)\nonumber\\
&=\sigma_{d-1}(C(e,a,\varepsilon))\int_{a}^{a+\frac{\varepsilon}{2}}
\cos^{d-1} (t)\, \dint t.
\end{align}
Using Jensen's inequality, basic trigonometric identities and $\varepsilon\le \pi/2-a$, we get
\begin{align}\label{eq:lb3}
\int_{a}^{a+\frac{\varepsilon}{2}}
\cos^{d-1} (t)\, \dint t&\ge \frac{\varepsilon}{2}\left(\frac{2}{\varepsilon}\int_{a}^{a+\frac{\varepsilon}{2}}
\cos t\, \dint t\right)^{d-1}
=\frac{\varepsilon}{2}\left(\frac{2}{\varepsilon}\,2\cos\left(a+\frac{\varepsilon}{4}\right)
\sin\left(\frac{\varepsilon}{4}\right)\right)^{d-1}\nonumber\\
&\ge\frac{\varepsilon}{2}\left(\frac{2}{\varepsilon}\,\frac{2}{\pi}\,\frac{\varepsilon}{4}\,2\,\cos\left(a+\frac{\pi}{8}-\frac{a}{2}\right)\right)^{d-1}
\nonumber\\
&=\frac{\varepsilon}{2}\left(\frac{2}{\pi}\right)^{d-1}\sin^{d-1}
\left(\frac{3}{4}\left(\frac{\pi}{2}-a\right)\right)\nonumber\\
&\ge \frac{3^{d-1}}{2\pi^{2d-2}}\varepsilon \left(\frac{\pi}{2}-a\right)^{d-1}.
\end{align}
Hence, we finally deduce from \eqref{eq:lb2}, \eqref{eq:lowvolbd} and \eqref{eq:lb3} that
\begin{align*}
\int_{\mathbb{S}^d}\mathbf{1}\{x\in A\}\, \sigma_d(\dint x)&
\ge\frac{1}{4d}{3^{d-1}}{\pi^{4(1-d)}}\omega_d
\varepsilon^{\frac{d+1}{2}}a^{d-1} \left(\frac{\pi}{2}-a\right)^{d-1}.
\end{align*}
From the upper bound in \eqref{eq:UsBa2} we finally get
$$
U_1(K)\ge \left(1+ c(d) a^{d-2}\left(\frac{\pi}{2}-a\right)^{d-1}\, \varepsilon^{\frac{d+1}{2}} \right)U_1(B_s(e,a)),
$$
where $c(d)\ge (8d)^{-1}(3\pi^{-4})^{d-1}\ge 2\cdot  0.02^{d}$.

To complete the argument, we have to verify \eqref{eq:auxbound}, that is, the asserted lower bound for the spherical radius of the cap $C(e,a,\varepsilon)$. By symmetry, it
is sufficient to consider the case $d=2$. Then the boundary of $\conv(B_s(e,a)\cup\{z_0\})$ is the union
of  two geodesic segments, denoted by $[z_0,p]$, $[z_0,\bar p]$, and the arc of $\partial B_s(e,a)$ connecting $p$ and
$\bar p$ which does not meet the geodesic segment $[e,z_0]$. Let $y\in [e,z_0]$ be such that $d_s(y,e)=a+\varepsilon/2$.
Further, let $q$ and $\bar q$ be the intersections of the geodesic through $y$ orthogonal to $[e,z_0]$ with $[p,z_0]$
and $[\bar p,z_0]$, respectively. Finally, let $\alpha\defeq\angle(e,z_0,p)=\angle (y,z_0,q)$ and $\omega\defeq\angle(q,e,y)$. Applying the sine rule in the right spherical triangle $\triangle(p,e,z_0)$ and Napier's rules for the right spherical triangles $\triangle(z_0,y,q)$ and $\triangle(e,q,y)$ (in this order), we obtain
$$
\sin(\alpha)=\frac{\sin(a)}{\sin(a+\varepsilon)},\quad \sin\left(\frac{\varepsilon}{2}\right)=\tan(d_s(y,q))\cot(\alpha),
\quad \sin\left(a+\frac{\varepsilon}{2}\right)=\tan(d_s(y,q)) \cot(\omega).
$$
Combining these relations, we deduce
\begin{align*}
\tan(\omega)&=\sin\left(\frac{\varepsilon}{2}\right)\frac{\sin(a)}{\sin\left(a+\frac{\varepsilon}{2}\right)}
\frac{1}{\sqrt{\sin(a+\varepsilon)+\sin(a)}}\frac{1}{\sqrt{\sin(a+\varepsilon)-\sin(a)}}\\
&\ge \frac{2}{\pi}\,\frac{\varepsilon}{2}\,\frac{2}{\pi}\, a\,\frac{1}{\sqrt{2}}\,  \frac{1}{\sqrt{\varepsilon}}
= \sqrt{2}\pi^{-2} a\sqrt{\varepsilon},
\end{align*}
which is the required bound.
\end{proof}

\subsection{Framework of  Example \ref{ex3}}
\label{subsec:ex3}

For $e\in \mathbb{S}^d$, we consider
$$\overline{\mathcal{K}}^d_{e}\defeq\{K\in \overline{\mathcal{K}}^d_s:
 e\in K\subset B_s(e,\pi/2\}
$$
and define
$$
\widetilde{U}_e(K)\defeq\int_{\mathbb{S}^d}\mathbf{1}\{(x-e)^\perp\cap K\neq\emptyset\}\, \sdn(\dint x),\quad K\in \overline{\mathcal{K}}^d_e.
$$
In analogy to the measure $\mu$ on $\BS^d_{d-1}$, we also introduce
$$
\tilde{\mu}_e(\cdot)\defeq \int_{\mathbb{S}^d}\mathbf{1}\{\BS_{x-e}\in\cdot\}\, \sdn(\dint x),
$$
so that $\widetilde{U}_e(K)=\tilde{\mu}_e(\mathcal{H}_K)$ for $K\in\overline{\mathcal{K}}^d_s$.

For $a\in (0,\pi/2)$, we have $(x-e)^\perp\cap B_s(e,a)\neq\emptyset$ if and only if $x\in B_s(e,2a)$, and hence it follows that
$\widetilde U_e (B_s(e,a))=\sigma_d(B_s(e,2a))$. H\"older's inequality and the relation $1-\cos(2a)=2\sin^2(a)$ imply that
$$
\sigma_d(B_s(e,2a))=\omega_d\int_0^{2a}\sin^{d-1}(t)\, \dint t\ge \omega_d\, 2a\left(\frac{1-\cos(2a)}{2a}\right)^{d-1}\ge 2\omega_d\left(\frac{2}{\pi}\right)^{2(d-1)}a^d.
$$
For an upper bound, we use $\sin(t)\le t$, $t\in [0,\pi]$, to get $\sigma_d(B_s(e,2a))\le d^{-1}\omega_d2^da^d$.
 Thus we have
\begin{equation}\label{lubounds}
2\frac{\omega_d}{\omega_{d+1}}\left(\frac{2}{\pi}\right)^{2(d-1)}a^d\le \widetilde U_e (B_s(e,a))=\sigma_d(B_s(e,2a))\le  d^{-1}\frac{\omega_d}{\omega_{d+1}}2^da^d.
\end{equation}

\begin{lemma}\label{Utildecont}
Let $e\in\BS^d$.
The functional $\widetilde U_e$ is continuous on $\overline{\mathcal{K}}^d_e$, $\tilde{\mu}_e\le 2^d\mu$, and
 $\widetilde U_e\le 2^{d+1} U_1$.
\end{lemma}

\begin{proof}
We first prove the second assertion in a slightly stronger form. A point $x\in \BS^d$ can be parameterized in the form
$x=\cos(\varphi)\, e+\sin(\varphi)\, u$ with $\varphi\in[0,\pi]$ and $u\in\BS_e$. Then we obtain
\begin{align*}
&\int_{\BS^d}\mathbf{1}\left\{(x-e)^\perp\cap \BS^d\in \cdot\right\}\, \sigma_d(\dint x)\\
&\qquad = \int_{\BS_e}\int_0^\pi \mathbf{1}\left\{\text{lin}\left\{\cos\left(\frac{\varphi}{2}\right)\,
e+\sin\left(\frac{\varphi}{2}\right)\, u,e^\perp\cap u^\perp\right\}\cap
\BS^d\in\cdot\right\}\sin^{d-1}(\varphi)\, \dint \varphi\, \sigma_{d-1}(\dint u)\\
&\qquad = 2\int_{\BS_e}\int_0^{\pi/2} \mathbf{1}\left\{\text{lin}\{\cos (s)\, e+\sin (s)\, u,e^\perp\cap u^\perp\}\cap
\BS^d\in\cdot\right\}\sin^{d-1}(2s)\,
\dint  s\, \sigma_{d-1}(\dint u)\\
&\qquad = 2\int_{\BS_e}\int_{\pi/2}^{\pi} \mathbf{1}\left\{\text{lin}\{\sin (t)\, e-\cos (t)\, u,e^\perp\cap u^\perp\}\cap
\BS^d\in\cdot\right\}\sin^{d-1}(2t-\pi)\, \dint  t\, \sigma_{d-1}(\dint u).
\end{align*}
Since $\sin(2t-\pi)\le 2\sin(t)$ for $t\in [\pi/2,\pi]$, we obtain
\begin{align}
&\int_{\BS^d}\mathbf{1}\left\{(x-e)^\perp\cap \BS^d\in \cdot\right\}\, \sigma_d(\dint x)\nonumber\\
&\qquad \le  2^d\int_{\BS_e}\int_{\pi/2}^{\pi} \mathbf{1}\left\{\text{lin}\{\sin (t)\, e-\cos (t)\, u,e^\perp\cap u^\perp\}
\cap \BS^d\in\cdot\right\}\sin^{d-1}(t)\, \dint  t\, \sigma_{d-1}(\dint u)\nonumber\\
& \qquad = 2^d\int_{\BS_e}\int_{\pi/2}^{\pi} \mathbf{1}\left\{(\cos (t)\, e+\sin (t)\, u)^\perp \cap
\BS^d\in\cdot\right\}\sin^{d-1}(t)\, \dint  t\, \sigma_{d-1}(\dint u)\nonumber\\
&\qquad \le 2^d\int_{\BS^d}\mathbf{1}\left\{x^\perp\cap \BS^d\in \cdot\right\}\,
\sigma_d(\dint x),\label{starast}
\end{align}
which yields the second and the third assertion.

Now we prove the continuity assertion. For $K\in \overline{\mathcal{K}}^d_e$, we put
$A_K\defeq\{z\in\BS^d:(z-e)^\perp\cap K\neq \emptyset\}$.

\noindent
Let $K,K_i\in \overline{\mathcal{K}}^d_e$, $i\in\N$, with  $K_i\to K$ as $i\to\infty$.
For $x\in \BS^d$, we distinguish the following cases.

If $(x-e)^\perp\cap \text{relint}(K)\neq\emptyset$ and $x\neq e$, then
$(x-e)^\perp\cap K_i\neq\emptyset$ if $i$ is sufficiently large,
and hence $\mathbf{1}_{A(K_i)}(x)\to \mathbf{1}_{A(K)}(x)$ as $i\to\infty$.

If $(x-e)^\perp\cap  K=\emptyset$, then $(x-e)^\perp\cap K_i=\emptyset$ if $i$ is sufficiently large,
and hence $\indi_{A(K_i)}(x)\to \indi_{A(K)}(x)$ as $i\to\infty$.

In addition, \eqref{starast} yields
\begin{align*}
&\sigma_d\left(\left\{z\in\BS^d:(z-e)^\perp\cap\text{relint}(K)=\emptyset,(z-e)^\perp\cap K \neq\emptyset\right\}
\right)\\
&\qquad \le  2^d\sigma_d\left(\left\{x\in \BS^d:x^\perp\cap\text{relint}(K)=\emptyset,x^\perp\cap K\neq\emptyset\right\}\right)\\
%&\qquad =2^d\nu_{d+1}\left(\left\{\rho\in SO_{d+1}:(\rho %e^\perp)\cap\text{relint}(K)=\emptyset,(\rho e^\perp)\cap K\neq\emptyset\right\}\right)\\
&\qquad\le 2^d\nu_{d+1}\left(\left\{\rho\in \SO_{d+1}:K\text{ and }\rho e_o^\perp\text{ touch each other}\right\}\right)=0,
\end{align*}
where $e_o\in\BS^d$ is arbitrary and fixed and we used \cite[Korollar 5.2.1]{Glas} (or \cite[Corollary 2.3.11]{Schneider14}) for the final equality.

Thus we have $\indi_{A(K_i)}(x)\to \indi_{A(K)}(x)$ as $i\to\infty$ for $\sigma_d$-almost all $x\in\BS^d$, so
that the assertion follows from the dominated convergence theorem.
\end{proof}

For a given $e\in \mathbb{S}^d$ and  $K\in\overline{\mathcal{K}}^d_e$, we define the size functional $\Sigma_e$ by $\Sigma_e(K)\defeq r_e(K)$ and the deviation functional $\vartheta_e$  by $\vartheta_e(K)\defeq R_e(K)-r_e(K)$, that is, both functionals are considered on a restricted domain in comparison with general size and hitting functionals, respectively.

\begin{theorem}\label{stab3}
Let $e\in \mathbb{S}^d$. Let $K\in\overline{\mathcal{K}}^d_e$  and $a\in (0,\pi/2)$. If $\Sigma_e(K)\ge a$, then
$$\widetilde U_e(K)\ge \widetilde U_e(B_s(e,a))=\sdn(B_s(e,2a)).$$
If also  $\vartheta_e(K)\ge \varepsilon\in (0,1]$,
then
$$
\widetilde U_e(K)\ge \left(1+c(a,d)\, \varepsilon^{\frac{d+1}{2}} \right)\,\widetilde U_e(B_s(e,a)),
$$
where $c(a,d)\ge (2\pi)^{-2d}a^{-1}
(\min\left\{a,\tfrac{\pi}{2}-a\right\})^{d-1}$.
\end{theorem}

\begin{proof}
If $\Sigma_e(K)\ge a\in (0,\pi/2)$ and $\vartheta_e(K)\ge \varepsilon\in(0,1]$, then there is some $z_0\in K\subset B_s(e,\pi/2)$ with $d_s(e,z_0)=a+\varepsilon\le \pi/2$ and
$\conv(B_s(e,a)\cup\{z_0\})\subset K$. We put
$$
\widetilde{A}\defeq \{x\in \mathbb{S}^d:(x-e)^\perp\cap B_s(e,a)=\emptyset, (x-e)^\perp\cap \conv(B_s(e,a)\cup\{z_0\})\neq\emptyset\}.
$$
Arguing as before, we get
$$
\widetilde U_e(K)\ge \widetilde U_e(B_s(e,a))+\int_{\mathbb{S}^{d}}\mathbf{1}\{x\in \widetilde{A}\}\, \sdn(\dint x).
$$
Moreover
\begin{align*}
\int_{\mathbb{S}^{d}}\mathbf{1}\{x\in \widetilde{A} \}\, \sigma_d(\dint x)&\ge
\int_{\mathbb{S}^d_e}\int_0^\pi \mathbf{1}\left\{\cos\left(\frac{\varphi}{2}\right)\, e+
\sin\left(\frac{\varphi}{2}\right)\, u\in \conv(B_s(e,a)\cup\{z_0\})\setminus B_s(e,a)\right\}\\
&\qquad\qquad\qquad\qquad \times \sin^{d-1}(\varphi)\, \dint\varphi\,\sigma_{d-1}(\dint u)\\
&\ge
\int_{C(e,a,\varepsilon)}\int_0^\pi \mathbf{1}\{\varphi/2\in (a,a+\varepsilon/2)\}\sin^{d-1}(\varphi)\, \dint \varphi\,\sigma_{d-1}(\dint u).
\end{align*}
For $\varphi \in (2a,2a+\varepsilon)$, we have
\begin{align*}
\sin\varphi&\ge \min\{\sin(2a),\sin(2a+\varepsilon)\}\ge \min\{\sin(2a),\sin(a+\pi/2)\}
\ge \frac{1}{2}\min\{a,\pi/2-a\}.
\end{align*}
The second inequality is clear if the minimum is attained by $\sin(2a)$. If $\sin(2a+\varepsilon)< \sin(2a)$, then $\pi/2<2a+\varepsilon=a+(a+\varepsilon)\le a+\pi/2$ and therefore $\sin(a+\pi/2)\le \sin(2a+\varepsilon)$. For the third inequality, we observe that if $a\ge \pi/4$, then $\pi/2\le 2a\le a+\pi/2$. Hence
\begin{align*}
\int_{\mathbb{S}^{d}}\mathbf{1}\{x\in \widetilde{A} \}\, \sdn(dx)&\ge
2^{1-d}(\min\{a,\pi/2-a\})^{d-1}\frac{\sigma_{d-1}(C(e,a,\varepsilon))}{\omega_{d+1}} \,\varepsilon\\
 &\ge \widetilde U_e(B_s(e,a))(2\pi)^{-2d}a^{-1}(\min\{a,\pi/2-a\})^{d-1} \,\varepsilon^{\frac{d+1}{2}},
 \end{align*}
 where we used \eqref{eq:lowvolbd} and \eqref{lubounds}.
\end{proof}

\section{Approximation results}
\label{sec:approx}

In this section, we provide two auxiliary results. First, we
establish an approximation result by means of spherical polytopes, from which we subsequently derive a crucial consequence for hitting functionals. Recall that $\delta_s$
 denotes the spherical Hausdorff distance.

\begin{lemma}\label{approx}
Let $K\in\overline{\cK}_s^{d}$. Then there are constants $k_1$ and $b_1$, depending only on $d$, such that for all $k\geq k_1$ there
is a spherical polytope $Q\in\cK^d_s$ with $k$ vertices (which can be chosen on the boundary of $K$ if $K\in{\cK}_s^{d}$) satisfying
$$\delta_s(K,Q)\leq b_1k^{-2/(d-1)}.$$
\end{lemma}

\begin{proof}
Since every convex body in $\overline{\cK}_s^{d}$ can be approximated by proper convex bodies, it is sufficient to consider the case where $K\in\cK^d_s$.

Let $\varepsilon\in (0,\pi/2)$, and first assume that the circumradius of $K$ satisfies $R(K)\leq\pi/2-\varepsilon$.  Below $\varepsilon$ will be chosen as a function of $d$. Let
$z(K)\in\BS^d$ denote the centre of the spherical
circumball of $K$. As in Section \ref{subsec:ex1} we consider the radial map
\begin{equation*}
\proj_{z(K)}:  \BS^d\cap H^+(z(K))\rightarrow  T_{z(K)}=z(K)+z(K)^\perp,\quad  x \mapsto \langle x,z(K)\rangle^{-1}\, x.
\end{equation*}
The image $\proj_{z(K)}(K)$ is contained in a $d$-dimensional Euclidean ball in $z(K)+z(K)^\perp$ centered at $z(K)$ with radius bounded from above by
$R_\circ \leq \cot(\varepsilon)$.
Applying the main result from \cite{Poly}, we  get constants $k_0=k_0(d)$ and $b_0=b_0(d)$ such that the following is
true. For $k\in\N$, $k\geq k_0$, there is a polytope $Q_0\subset  z(K)+z(K)^\perp $ with at most $k$ vertices, located on the boundary of
$\proj_{z(K)}(K)$, satisfying
\begin{equation*}
\delta(R_\circ^{-1}\proj_{z(K)}(K),R_\circ^{-1}Q_0)\leq b_0 k^{-2/(d-1)}.
\end{equation*}
Here $\delta$ denotes the Hausdorff distance in $\R^{d+1}$. The polytopes $R_\circ^{-1}\proj_{z(K)}(K)$ and $R_\circ^{-1}Q_0$ lie in a
 common affine
subspace parallel to $z(K)^\bot$.
Therefore $Q_0\subset \proj_{z(K)}(K)\subset Q_0 +R_\circ\, b_0\, k^{-2/(d-1)}B_E$, where $B_E$ is the unit ball in $z(K)^\bot$, and thus
$$\delta(\proj_{z(K)}(K),Q_0)\leq \cot(\varepsilon)\, b_0\, k^{-2/(d-1)}.$$
The mapping
 $\Pi_{z(K)}^{-1}: z(K)+z(K)^\perp\rightarrow  \BS^d\cap H^+(z(K))$, $x\mapsto \|x\|^{-1}{x}$,
is Lipschitz continuous with Lipschitz constant at most $2$. Since also
$({\pi}/{2})\,\|x-y\|\geq2\,\arcsin\left( {\|x-y\|}/{2}\right)=d_s(x,y)$, for $x,y\in\BS^d$,
 the spherical polytope $Q\defeq \Pi_{z(K)}^{-1}(Q_0)$ satisfies $Q\subset \Pi_{z(K)}^{-1}(\proj_{z(K)}(K))=K$ and
$$
\delta_s(Q,K)\leq \cot(\varepsilon)\, b_0\,\pi\, k^{-2/(d-1)},\quad k\geq k_0.
$$
For arbitrary $K\in\ {\cK}_s^{d}$, we divide $K$ into $2^d$ pieces $K\cap \cap_{i=1}^dH^*(\epsilon_i e_i)$, $\epsilon_i\in\{-1, 1\}$, where $e_0,e_1,\ldots,e_d$ is an orthonormal basis of $\R^{d+1}$ with $K\subset H^+(e_0)$.
Then each piece of $K$ is contained in a regular spherical $d$-simplex of edge-length $\pi/2$, which is the spherical convex
hull of $d+1$ unit vectors $e_0,\epsilon_1e_1,\ldots,\epsilon_de_d$. Its circumradius is $\arccos(1/\sqrt{d+1})\in(\pi/4,\pi/2)$; see \cite[Theorem 2]{DW}. Defining
$$\varepsilon\defeq\pi/2-\arccos\left( {\sqrt{d+1}^{-1}}\right),$$
the individual pieces satisfy $R(K_i)\leq\frac{\pi}{2}-\varepsilon$,
$ i=1,\ldots,2^{d}$. Applying the reasoning above to every piece, we obtain spherical polytopes $Q_i$, $ i=1,\ldots,2^{d}$, each having at most $k$ vertices,  such that
$$\delta_s(Q_i,K_i)\leq \cot(\varepsilon)\, b_0\,\pi\, k^{-2/(d-1)},\quad i=1,\ldots,2^{d},\ k\geq k_0.$$
Defining $Q\defeq \conv(\bigcup_{i=1}^{2^{d}}Q_i)$ (here each spherical polytope $Q_i$ can be replaced by the set of its vertices), we obtain a proper spherical polytope with at most $2^{d}\, k$ vertices and
$$
\delta_s(Q,K)\leq \cot(\varepsilon)\, b_0\,\pi\, k^{-2/(d-1)}.
$$
Since $\cot(\varepsilon)={\cot\left(\pi/2-\arccos(1/\sqrt{d+1})\right)}=\sqrt{d}$,
the assertion follows with $k_1=2^{d}\, k_0$ and $b_1=b_0\,\pi\,\sqrt{d}\, 4^{d/(d-1)}$. Possibly enlarging $Q$, we can assume that all vertices of $Q$ are in the boundary of $K$.
\end{proof}

\medskip

The following lemma states that the values of $\Phi$ on spherical polytopes can be approximated by the values
of $\Phi$ on proper spherical polytopes with a controlled number of vertices (extreme points). We write $\text{ext}(P)$ for the set
of extreme points of a proper spherical polytope $P$ and $f_0(P)$ for the number of its extreme points.

\begin{lemma}\label{Ecken}
Let $\Phi$ be a hitting functional, $\Sigma$ a size functional, and let $a,\alpha>0$. Then there is an integer $\nu\in\N$, depending only on $\Phi, \Sigma,d,a,\alpha$, such that
for every spherical polytope $P\in \cK^d_s$ with $\Sigma(P)\geq a$ there is a proper spherical polytope $Q=Q(P)$ satisfying
$\ext(Q)\subset\ext(P)$, $f_0(Q)\leq\nu$ and
$$
\Phi(Q)\geq(1-\alpha)\Phi(P).
$$
Furthermore, the mapping $P\mapsto Q(P)$ can be chosen to be measurable.
\end{lemma}

\begin{proof}
Since the functional $\Phi:\overline{\cK}^d_s\to [0,\infty)$ is continuous
on the compact space $\overline{\cK}_s^d$  (with respect to the spherical Hausdorff metric $\delta_s$), $\Phi$ is  uniformly continuous.
Let $\Phi,\Sigma,\alpha,a$   be as in the statement of the lemma and define $\varepsilon\defeq\alpha\,\tau(\Phi,\Sigma,a)$. From the uniform continuity of $\Phi $ it follows that
there is some $\delta=\delta(\varepsilon)>0$
  such that $\left|\Phi(K)-\Phi(K')\right|\leq\varepsilon=\alpha\,\tau(\Phi,\Sigma,a)$
for all $K,K'\in\overline{\cK}_s^d$ with $\delta_s(K,K')\leq\delta(\varepsilon)$. Let $P$ be a spherical polytope with $\Sigma(P)\ge a$.
From Lemma \ref{approx} we obtain a spherical
polytope $Q=Q(P)$ and a number $\nu=\nu(\Phi,\Sigma,d,a,\alpha)$ such that
$\text{ext}(Q)\subset\text{ext}(P) $, $f_0(Q)\le \nu$ and  $ \delta_s(Q,P)\leq\delta(\varepsilon)$.
Since $\Sigma(P)\geq a$, we
conclude that $\Phi(P)\ge \tau(\Phi,\Sigma,a)$, and therefore
$$
\Phi(P)-\Phi(Q)\le |\Phi(P)-\Phi(Q)|\leq\varepsilon=\alpha\,\tau(\Phi,\Sigma,a)\leq\alpha\,\Phi(P),
$$
which yields the first assertion.

Identifying each spherical polytope with a Euclidean polytope which is the convex hull of the Euclidean origin and the vertices of the spherical polytope, we conclude that the second assertion follows as in \cite[Lemma 4.2]{HRS}.
\end{proof}

\section{Probabilistic inequalities}
\label{sec:probabin}
After the geometric preparations of the preceding sections, we can proceed with estimating the conditional probabilities involved in
resolving Kendall's problem in spherical space.  In each case, the conditional probability  is the ratio of two probabilities. The probability in the denominator
is easy to treat. In the following lemma, we provide an upper bound for the numerator, which is the main step in the probabilistic estimate.

Throughout this section, we consider a hypersphere  tessellation, generated by the hypersphere process $\widetilde{X}=h(X)$, where $X$ is an isotropic Poisson point process on $\BS^d$ with intensity $\gamma_s\in (0,\infty)$. Let $\Phi=2U_1$, which we called the hitting functional associated with $\widetilde{X}=h(X)$ (which is rotation invariant). We write again $\mu$ for the image measure of $\sdn$ under $h$.
Let $\Sigma$ be a general increasing size functional, and let $a>0$ be such that $\Sigma^{-1}([a,\infty))\neq\emptyset$. Let $\vartheta$ be a deviation functional for $\Phi,\Sigma,a$. We write $\tau(a)=\tau(\Phi,\Sigma,a)$ for the isoperimetric constant based on
the parameters $\Phi,\Sigma$ and $a>0$ and $f_a$ for the stability function as in Proposition \ref{propneu1}. Note that if $\Sigma$ is rotation invariant, then the class of extremal bodies is also rotation invariant, and the same is true for the canonical deviation functional, defined in \eqref{devfunc}. The Crofton cell $Z_0$ is the almost surely unique cell of the tessellation containing  $\sfp$ (in its interior).

\begin{lemma}\label{zaehler}
Let $a>0$, $\varepsilon \in (0,1]$ and
$$
\cK_{a,\varepsilon}\defeq\{K\in\overline{\cK}_s^{d}: \Sigma(K)\ge a, \vartheta(K)\geq\varepsilon\}.
$$
If $\bar\kappa\in(0,1)$, then
\begin{equation}\label{phszaehler}
\BP(Z_0\in\cK_{a,\varepsilon})\leq c_1\, \max\{1,\gamma_s\}^{d \nu}\,\exp\left(-\gamma_s\left(1+
(1-\bar\kappa)f_a(\varepsilon)\right)\tau(a)\right),
\end{equation}
where the constants $c_1$ and $\nu$ depend on $a,d,\varepsilon,\bar\kappa,\Sigma,\vartheta$.
\end{lemma}

\begin{proof}
By assumption, we have  $\Sigma^{-1}([a,\infty))\neq\emptyset$.
Let $N\in\mathbb{N}$. For $H_1,\ldots,H_N\in \mathcal{H}_{\mathbb{S}^d}$ such that $\sfp\notin H_i$, for $i=1,\ldots,N$, we define
$H_{(N)}\defeq(H_1,\ldots,H_N)$ and let $P(H_{(N)})$
denote the spherical Crofton cell of the tessellation induced by $H_1,\ldots,H_N$.
In what follows, we consider $H_1,\ldots,H_N\in \mathcal{H}_{\mathbb{S}^d}$
such that $P(H_{(N)})\in\cK_{a,\varepsilon}\cap\cK_s^d$.
This requires $N\geq d+1$. If $N\geq d+1$ and $H_1,\ldots,H_N$ are i.i.d. with a distribution which has a density with
respect to the invariant measure on $\BS^d_{d-1}$, then $P(H_{(N)})\in\cK_s^d$ is satisfied almost surely.

Define $\alpha\defeq\bar\kappa f_a(\varepsilon)/(1+f_a(\varepsilon))$, hence $(1-\alpha)(1+f_a(\varepsilon))=1+\bar\alpha$, where $\bar\alpha\defeq(1-\bar\kappa)f_a(\varepsilon)$.
By Lemma \ref{Ecken} and Proposition \ref{propneu1}, there are $\nu=\nu(\Sigma,\vartheta,d,a,\varepsilon,\bar\kappa)$ vertices of $P(H_{(N)})
\in\cK_{a,\varepsilon}\cap\cK_s^d$ such that
the spherical convex hull $Q(H_{(N)})$ of these vertices satisfies
\begin{align*}
1\geq\Phi(Q(H_{(N)}))&\geq(1-\alpha)\Phi(P(H_{(N)}))\\
&\geq(1-\alpha)(1+f_a(\varepsilon))\tau(a)\\
&=(1+\bar\alpha)\tau(a),
\end{align*}
where we used $\Phi(\cdot)\leq1$. By Lemma \ref{Ecken}, we can assume that the map
$(H_1,\ldots,H_N)\mapsto Q(H_{(N)})$
is measurable. Since $\mu=h(\sdn)$,
every vertex of $Q(H_{(N)})$ lies $\mu^N$-almost surely in exactly $d$ of these hyperspheres. The remaining
subspheres do not hit
$Q(H_{(N)})$. Hence, the number of hyperspheres hitting $Q(H_{(N)})$ is $j\in\{d+1,\ldots,d\,\nu\}$.
Without loss of generality we assume $H_1\cap Q(H_{(N)})\neq
\emptyset, \ldots,H_j\cap Q(H_{(N)})\neq\emptyset$.
Then there are subsets $J_1,\ldots, J_{f_0(Q(H_{(N)}))}$ of $\{1,\ldots,j\}$,
each of cardinality $d$, such that the intersections
$$\bigcap_{l\in J_i}H_l,\quad i=1,\ldots,f_0(Q(H_{(N)}))\leq\nu,$$
give the vertices of $Q(H_{(N)})$. In the following, we denote by $\sum_{(J_1,\ldots,J_\nu)}$  the sum over all $\nu$-tuples of
subsets of $\{1,\ldots,j\}$ with $d$ elements. Note that for $K\in\overline{\cK}_s^{d}$ we have
$$\int_{\cH_{\BS^d}}\mathbf{1}\{H\cap K=\emptyset\}\, \mu(\dint H)=1-\Phi(K).$$
Assuming $N\geq d+1$ and using that $\Phi(\BS^d)=1$, we obtain
\begin{align}\label{grossd}
&\BP(Z_0\in\cK_{a,\varepsilon}  \mid \widetilde{X}(\cH_{\BS^d})=N)
=\int_{\cH^N_{\BS^d}}\mathbf{1}\{P(H_{(N)})\in\cK_{a,\varepsilon}\cap\cK_s^d\}\, \mu^N(\dint (H_1,\ldots,H_N))\nonumber\\
&\leq\sum_{j=d+1}^{d \nu}\binom{N}{j}\int_{\cH^N_{\BS^d}}\mathbf{1}\{P(H_{(N)})\in\cK_{a,\varepsilon}\cap\cK_s^d\}
\mathbf{1}\{H_i\cap Q(H_{(N)})\neq\emptyset,\ i=1,\ldots,j\}\nonumber\\
&\qquad\qquad\qquad\qquad\times\mathbf{1}\{H_i\cap Q(H_{(N)})=\emptyset,\ i=j+1,\ldots,N\}\,  \mu^N(\dint (H_1,\ldots,H_N))\nonumber\\
&\leq\sum_{j=d+1}^{d \nu}\binom{N}{j}\sum_{(J_1,\ldots,J_\nu)}\int_{\cH_{\BS^d}^j}\int_{\cH_{\BS^d}^{N-j}}
\mathbf{1}\{\Phi(\conv\bigcup_{r=1}^\nu\bigcap_{i\in J_r}H_i)\geq (1+\bar\alpha)\tau(a)\}\nonumber\\
&\qquad\qquad\qquad\qquad \times\mathbf{1}\{H_l\cap\conv\bigcup_{r=1}^\nu\bigcap_{i\in J_r}H_i=\emptyset,\ l=j+1,\ldots,N\}\nonumber\\
&\qquad\qquad\qquad\qquad\times\mu^{N-j}(\dint (H_{j+1},\ldots,H_N))\, \mu^j(\dint (H_1,\ldots,H_j))\nonumber\\
&=\sum_{j=d+1}^{d \nu}\binom{N}{j}\sum_{(J_1,\ldots,J_\nu)}\int_{\cH_{\BS^d}^j}\mathbf{1}\{\Phi(\conv\bigcup_{r=1}^\nu\bigcap_{i\in J_r}H_i)\geq (1+\bar\alpha)\tau(a)\}\nonumber\\
&\qquad\qquad\qquad\qquad\times[1-\Phi(\conv\bigcup_{r=1}^\nu\bigcap_{i\in J_r}H_i)]^{N-j}\,  \mu^j(\dint(H_1,\ldots,H_j))\nonumber\\
&\leq\sum_{j=d+1}^{d \nu}\binom{N}{j}\binom{j}{d}^\nu\left[1-(1+\bar\alpha)\tau(a)\right]^{N-j}.
\end{align}
Summation over $N$ gives
\begin{align*}
\BP(Z_0\in\cK_{a,\varepsilon})
&\leq\sum_{N=0}^d\BP(\widetilde{X}(\cH_{\BS^d})=N)\\
&\qquad+\sum_{N=d+1}^\infty\BP(Z_0\in\cK_{a,\varepsilon}  \mid \widetilde{X}(\cH_{\BS^d})=N)\BP(\widetilde{X}(\cH_{\BS^d})=N).
\end{align*}
For the second sum, we deduce from (\ref{grossd}) that
\begin{align*}
&\sum_{N=d+1}^\infty\BP(Z_0\in\cK_{a,\varepsilon} \mid \widetilde{X}(\cH_{\BS^d})=N)\BP(\widetilde{X}(\cH_{\BS^d})=N)\\
&\leq\sum_{N=d+1}^\infty\sum_{j=d+1}^{d \nu}\binom{N}{j}\binom{j}{d}^\nu
\left[1-(1+\bar\alpha)\tau(a)\right]^{N-j}\frac{ \gamma_s  ^N}{N!}\exp(-\gamma_s )\\
&=\sum_{j=d+1}^{d \nu}\binom{j}{d}^\nu\exp(-\gamma_s)\frac{\gamma_s^j}{j!}\sum_{N=j}^\infty
\frac{[1-(1+\bar\alpha)\tau(a)]^{N-j}}{(N-j)!}\gamma_s^{N-j}\\
&=\exp\left(-\gamma_s(1+\bar\alpha)\tau(a)\right)\sum_{j=d+1}^{d \nu}\binom{j}{d}^\nu\frac{\gamma_s^j}{j!}\\
&= \sum_{j=d+1}^{d \nu}\binom{j}{d}^\nu\frac{\gamma_s^j}{j!}\,  \exp\left(-\gamma_s\left(1+(1-\bar\kappa)f_a(\varepsilon)\right)\tau(a)\right).%\\
\end{align*}
For the first sum, we get
\begin{align*}
\sum_{N=0}^d\BP(\widetilde{X}(\cH_{\BS^d})=N)
&\leq\sum_{N=0}^d\frac{\gamma_s^{N}}{N!}\exp\left(-\gamma_s\left(1+\bar\alpha
\right)\tau(a)\right),
\end{align*}
since $1\ge (1+\bar\alpha)\tau(a)= (1+(1-\bar\kappa)f_a(\varepsilon) )\tau(a)$.
Combining both estimates, we obtain
$$
\BP(Z_0\in\cK_{a,\varepsilon})\leq c_1\,\max\{1,\gamma_s\}^{d\nu}\,\exp\left(-\gamma_s(1+(1-\bar\kappa)f_a(\varepsilon))\tau(a)\right),
$$
where
$$
c_1=c_1(\Sigma,\vartheta,d,a,\varepsilon,\bar\kappa)\defeq\sum_{j=d+1}^{d  \nu}\binom{j}{d}^\nu\frac{1}{j!}+\sum_{N=0}^d\frac{1}{N!}
\le e+(d\nu)^{d \nu},
$$
which completes the argument.
\end{proof}

In the next lemma, we specify $\Sigma=\sigma_d$ and $\vartheta=\Delta_2^*$ (see the definition in Remark \ref{rem2}). Using Theorem \ref{stability} instead of Proposition \ref{propneu1},
we arrive at the following more explicit result. Let $\Ca$ denote a spherical cap of volume $a$.

\begin{lemma}\label{zaehler2}
Under the preceding assumptions, let $\Phi=2U_1$.
If $0<a<\omega_{d+1}/2$, $\varepsilon\in (0,1]$  and
$$
\overline{\cK}_{a,\varepsilon}\defeq\{K\in\overline{\cK}_s^{d}: \sigma_{d}(K)\ge a,  \Delta_2^*(K)\geq\varepsilon\},
$$
then
$$
\BP(Z_0\in\overline{\cK}_{a,\varepsilon})\leq c_2\,\max\{1,\gamma_s\}^{d\nu}\,
\exp\left(-\gamma_s(1+\bar \beta\, \varepsilon^{2(d+1)} )\Phi(\Ca)\right),
$$
where $c_2$ and $\nu$ depend on $a,d,\varepsilon$ and the constant $\bar\beta$ depends on $a,d$.
\end{lemma}

\begin{proof}
Let $Z_0\in\overline{\cK}_{a,\varepsilon}\cap\cK_s^d$, hence $\Delta_2(Z_0)=\Delta_2^*(Z_0)$. Suppose that all points of the underlying Poisson point process $X$ are
  in $\interior(B_s(e(Z_0),\varepsilon)\cup B_s(-e(Z_0),\varepsilon))$. Then we immediately get $\Delta_0(Z_0)<\varepsilon$. By (\ref{delta-absch}), it follows that
$\Delta_2(Z_0)<\varepsilon$, a contradiction to $Z_0\in\overline{\cK}_{a,\varepsilon}$.

Therefore, there is a point $x\in X$ such that $x\in B_s(e(Z_0),\pi/2)$ and $d_s(e(Z_0),x)\geq\varepsilon$
or there is a point $x\in B_s(-e(Z_0),\pi/2)$ and
$d_s(-e(Z_0),x)\geq\varepsilon$. In either case, we obtain
$$
\sigma_d(Z_0)\leq\frac{\omega_{d+1}}{2}-\frac{\omega_{d+1}}{2}\frac{\varepsilon}{\pi}.
$$
Now let $C$ be a spherical cap with $\sigma_d(C)\leq\omega_{d+1}/2-(\varepsilon/\pi)\cdot(\omega_{d+1}/2)$
and denote its radius by $\alpha_C$. Since $\alpha_C\le \pi/2$,
\begin{align*}
\sigma_d(C)=\omega_d\int_0^{\alpha_C}\sin^{d-1}(t)\, \dint t&\leq \frac{\omega_{d+1}}{2}-\frac{\varepsilon}{\pi}\frac{\omega_{d+1}}{2}=\omega_d\int_0^{\pi/2}\sin^{d-1}(t)\,  \dint t-\frac{\varepsilon\omega_{d+1}}{2\pi},
\end{align*}
and thus
\begin{equation*}
\frac{\varepsilon\omega_{d+1}}{2\pi\omega_d}\leq\int_{\alpha_C}^{\pi/2}\sin^{d-1}(t)\,  \dint t\leq\frac{\pi}{2}-\alpha_C,
\end{equation*}
which yields
\begin{equation}\label{ungl1a}
 \alpha_C\leq\frac{\pi}{2}-\frac{\varepsilon\omega_{d+1}}{2\pi\omega_d}.
 \end{equation}
Analogously to the proof of Lemma \ref{zaehler}, we consider
$N\in\mathbb{N}$, $N\ge d+1$,  and $H_1,\ldots,H_N\in \mathcal{H}_{\BS^d}$ such that the Crofton cell
$P(H_{(N)})$ of the induced tessellation satisfies $P(H_{(N)})\in\overline{\cK}_{a,\varepsilon}\cap\cK_s^d$.
Let $C$ be a
spherical cap satisfying $\sigma_d(C)=\sigma_d(P(H_{(N)}))\ge a$ and denote its radius by $\alpha_C$. Using Theorem
\ref{stability} instead of Proposition \ref{propneu1} and the monotonicity of $\Phi$, we get
\begin{equation}\label{proto}
\Phi(P(H_{(N)}))\geq(1+{\beta}\varepsilon^2)\Phi(C)\geq(1+{\beta}\varepsilon^2)\Phi(\Ca),
\end{equation}
where
$$
{\beta}=\min\left\{\frac{\frac{2}{d}\binom{d+1}{2}\sin^{d+1}(\alpha_C)
\tan^{-2d}(\alpha_C)}{1+\binom{d+1}{2}\left(\frac{\pi}{2}\right)^2\tan^{-d}
(\alpha_C)}, \frac{8}{\pi^2} D\left(\frac{\pi}{2}-\alpha_C\right)\right\}
$$
with
$$
\alpha_{\Ca}\le \alpha_C\leq\frac{\pi}{2}-\varepsilon\frac{\omega_{d+1}}{2\pi\omega_d}.
$$
Recalling
$$
D(x)=\int_0^x\sin^{d-1}(t)\, \dint t\geq\int_0^x\left(\frac{2}{\pi}\right)^{d-1}t^{d-1}\, \dint t=\frac{2^{d-1}}{d\pi^{d-1}}x^{d},
$$
and using the fact that $\tan$ is increasing on $[0,\pi/2)$ and $\tan(x)\geq x$ for $x\in[0,\pi/2)$, we obtain
\begin{align*}
{\beta}&\geq \min\left\{\frac{\frac{2}{d}\binom{d+1}{2}\sin^{d+1}(\alpha_{\Ca})
\tan^{-2d}\left(\frac{\pi}{2}-\varepsilon\frac{\omega_{d+1}}{2\pi\omega_d}\right)}{1+\binom{d+1}{2}
\left(\frac{\pi}{2}\right)^2\tan^{-d}(\alpha_{\Ca})},\left(\frac{2}{\pi}\right)^2D\left(\varepsilon
\frac{\omega_{d+1}}{2\pi\omega_d}\right)\right\}\\
&\geq 2 \min\left\{\frac{\frac{1}{d}\binom{d+1}{2}\sin^{d+1}(\alpha_{\Ca})\left(\frac{\omega_{d+1}}{2\pi\omega_d}
\right)^{2d}}{1+\binom{d+1}{2}\left(\frac{\pi}{2}\right)^2\tan^{-d}(\alpha_{\Ca})}
\varepsilon^{2d},\left(\frac{2}{\pi}\right)^{d+1}\frac{\omega_{d+1}^d}{d(2\pi\omega_d)^d}\varepsilon^{2d}\right\}
=:2\bar\beta\cdot\varepsilon^{2d},
\end{align*}
where we made use of $\varepsilon\leq 1$ in the second to last line. Note that $\bar \beta>0$ depends only on $a$ and $d$. Combining
this with (\ref{proto}), we get
$$
\Phi(P(H_{(N)}))\geq(1+2\bar \beta\,\varepsilon^{2(d+1)})\Phi(\Ca).
$$
Proceeding as in the proof of Lemma \ref{zaehler}, with $\bar\kappa=1/2$ and $f_a(\varepsilon)=2\bar\beta \varepsilon^{2(d+1)}$, we obtain the required  result.
\end{proof}

Now we are able to prove the following general theorem.

\begin{theorem}\label{thm:genresult}
Let $Z_0$ be the Crofton cell of a  hypersphere tessellation derived from an isotropic Poisson point process $X$ on $\BS^d$ with intensity $\gamma_s$.
Let $\Phi=2U_1$, let $\Sigma$ be an increasing and rotation invariant size functional, and let $a>0$ be such that $\Sigma^{-1}([a,\infty))\neq\emptyset$. Let
$\vartheta$ be a deviation functional for $\Phi,\Sigma,a$, and let $f_a$ be a stability function as in Proposition \ref{propneu1}.
If  $  \varepsilon\in (0,1]$ and $\kappa\in (0,1)$, then there is a constant $c_3>0 $ such that
$$
\BP\left(\vartheta(Z_0)\ge\varepsilon \mid \Sigma(Z_0)\geq a\right)\leq {c}_3\,\exp\left(-\kappa\gamma_sf_a(\varepsilon)\tau(a)\right),
$$
where the constant ${c}_3$ depends  on $\Sigma,a,\vartheta,d,\varepsilon,\kappa$.
\end{theorem}

\begin{proof}
First we note that
\begin{equation}\label{eq1}
\BP\left(\vartheta(Z_0)\geq\varepsilon \mid \Sigma(Z_0)\geq a\right)=
\frac{\BP\left(\vartheta(Z_0)\geq\varepsilon,\ \Sigma(Z_0)\geq a\right)}{\BP\left(\Sigma(Z_0)\geq a\right)}=
\frac{\BP(Z_0\in \cK_{a,\varepsilon})}{\BP\left(\Sigma(Z_0)\geq a\right)}.
\end{equation}
Let $K_a\in\mathcal{E}(\Phi,\Sigma,a)$ be an extremal body, hence $\Sigma(K_a)\ge a$ and
$\Phi(K_a)=\tau(a)$. Let $e_a\in K_a$ be an arbitrary fixed point and let $\rho_a\in\SO(d+1)$ be such that $\rho_ae_a=\sfp$. Since $\Sigma$ and $\Phi$ are rotation invariant, $K_{\sfp} \defeq \rho_a K_a\in\mathcal{E}(\Phi,\Sigma,a)$ with $\sfp\in K_\sfp$, $\Sigma(K_\sfp)\ge a$ and $\Phi(K_\sfp)=\tau(a)$. If $\widetilde X( \cH_{K_\sfp})=0$, then $K_\sfp\subset Z_{0}$. Since $\Sigma$ is increasing, we deduce from $\widetilde X( \cH_{K_\sfp})=0$ that
$\Sigma(Z_{0})\ge \Sigma(K_\sfp)\ge a$. Thus we arrive at
\begin{equation}\label{nenner}
\BP(\Sigma(Z_0)\geq a)\geq\BP(\widetilde{X}(\cH_{K_\sfp})=0)
=\exp(-\gamma_s\tau(a)).
\end{equation}
From (\ref{eq1}), (\ref{nenner}) and Lemma \ref{zaehler} with $\bar\kappa=(1-\kappa)/2$ we obtain
\begin{align*}
\BP\left(\vartheta(Z_0)\geq\varepsilon \mid \Sigma(Z_0)\geq a\right)
&\leq\frac{c_1\, \max\{1,\gamma_s\}^{d\nu}\,\exp\left(-\gamma_s\left(1+
(1-\bar\kappa)f_a(\varepsilon)\right)\tau(a)\right)}{\exp(-\gamma_s\tau(a))}\\
&=c_1\,\max\{1,\gamma_s\}^{d\nu}\,\exp\left(-\gamma_s(1-\bar\kappa)f_a(\varepsilon)\tau(a)\right)\\
&\leq c_3\,\exp\left(-\gamma_s(1-2\bar\kappa)f_a(\varepsilon)\tau(a)\right),
\end{align*}
where the constants $c_3,\nu$ depend on $a,\varepsilon,d,\kappa$ and, of course, on  $\Sigma,\vartheta$.
\end{proof}

Using Lemma \ref{zaehler2} instead of Lemma \ref{zaehler} (and of course $\vartheta=\Delta_2^*$), we obtain a similar result but with a
more explicit constant in the exponent, which was stated as Theorem B in the introduction.

\begin{figure}[ht]
  \centering
  \subfloat[$\gamma_s=1$, the realization contains $17$ great circles]{\label{Bildchen1}%
    \includegraphics[width=0.4\textwidth]
    {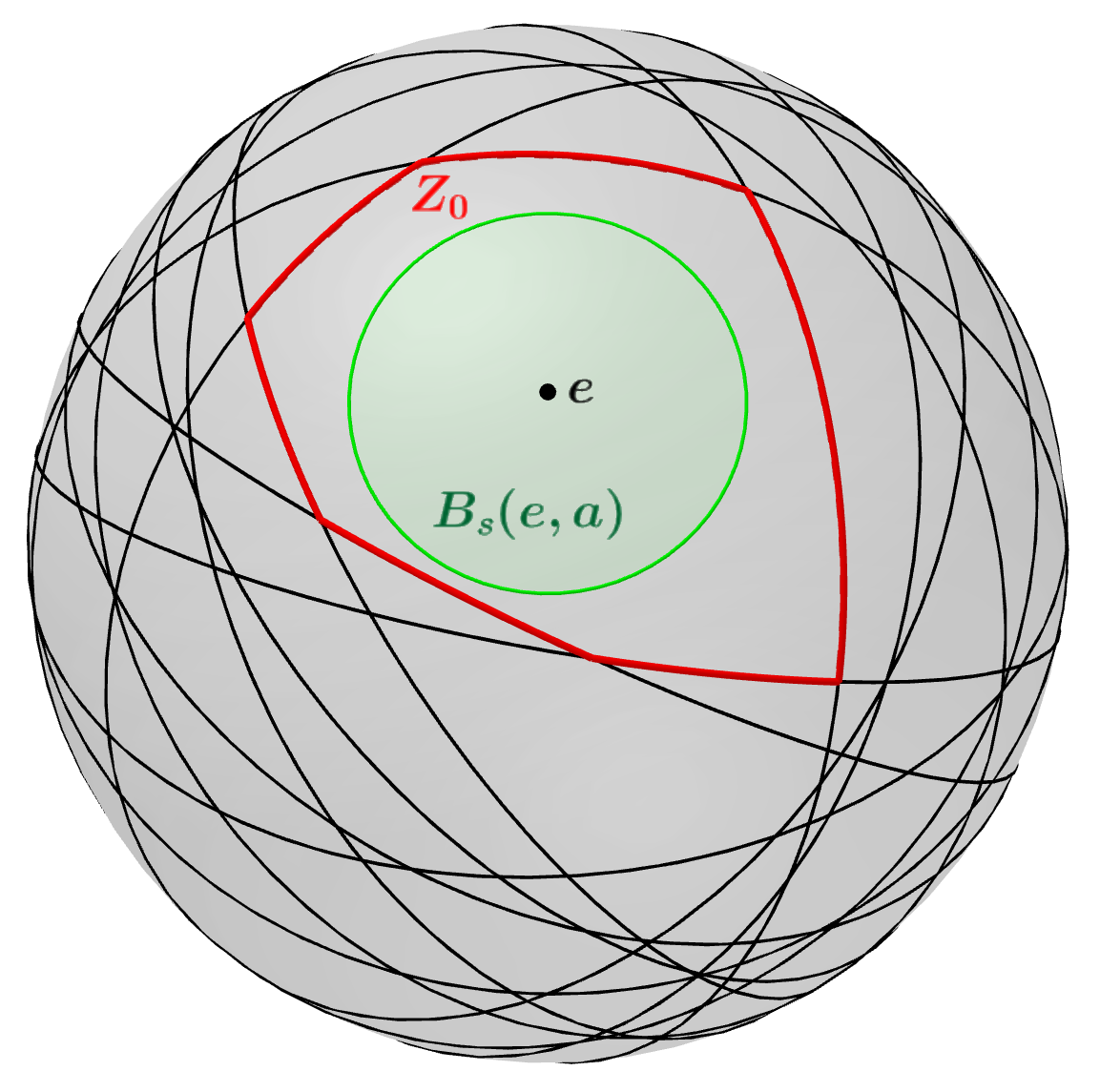}}\quad
  \subfloat[$\gamma_s=2$, the realization contains $31$ great circles]{\label{Bildchen2}
    \includegraphics[width=0.4\textwidth]
    {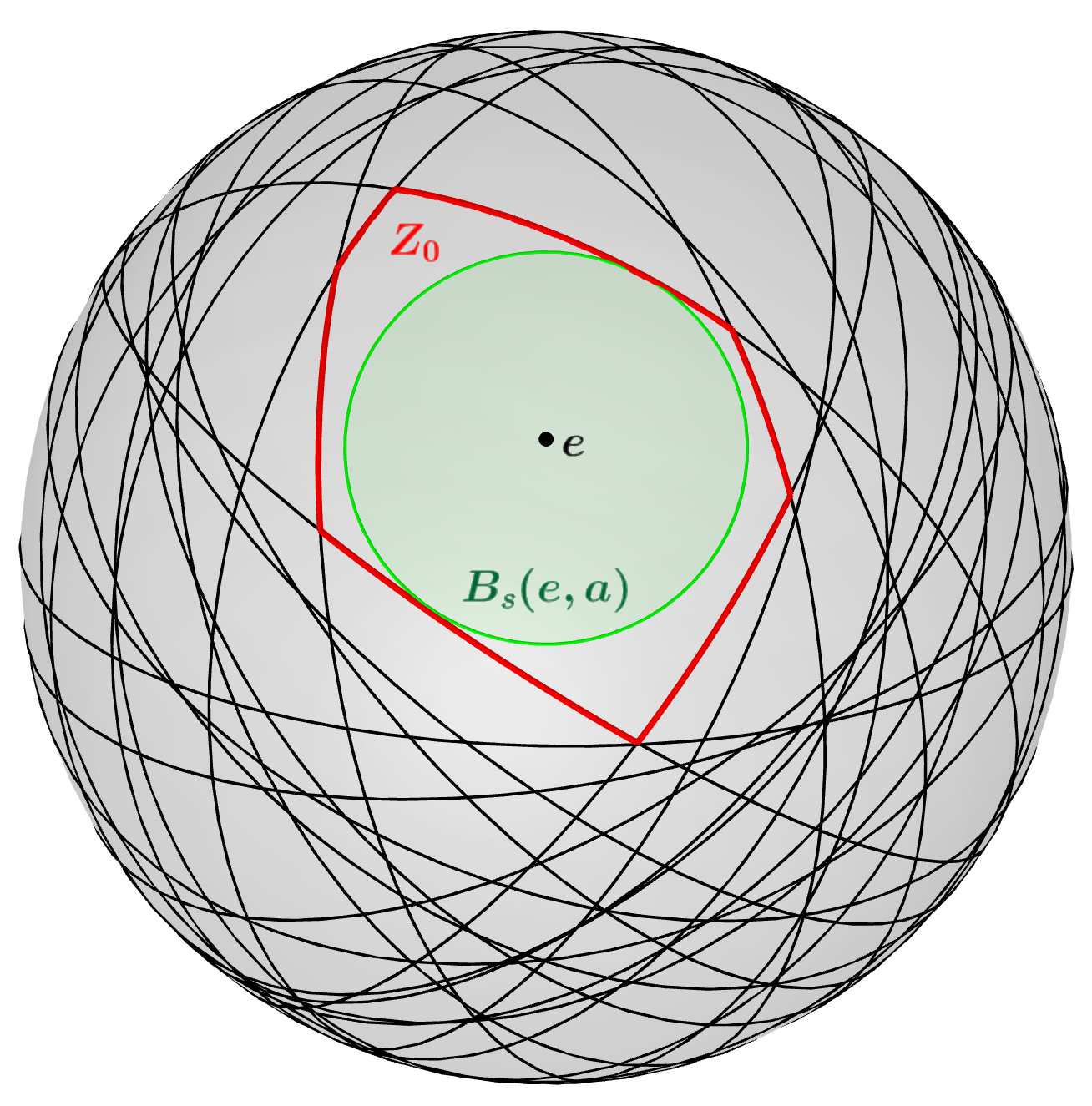}}\quad
  \subfloat[$\gamma_s=4$, the realization contains $61$ great circles]{\label{Bildchen3}%
    \includegraphics[width=0.4\textwidth]
    {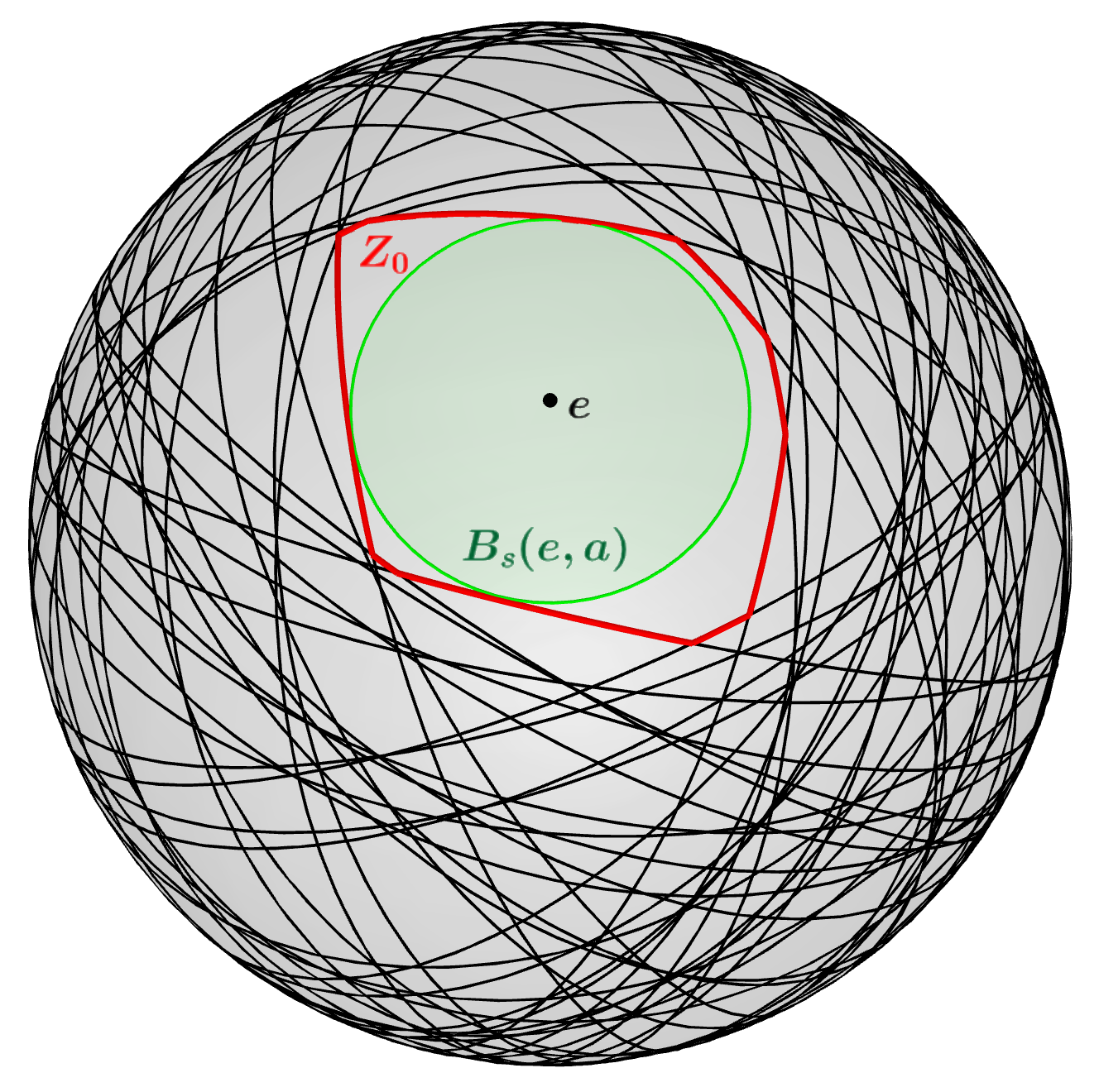}}\quad
  \subfloat[$\gamma_s=10$, the realization contains $118$ great circles]{\label{Bildchen4}
    \includegraphics[width=0.4\textwidth]
    {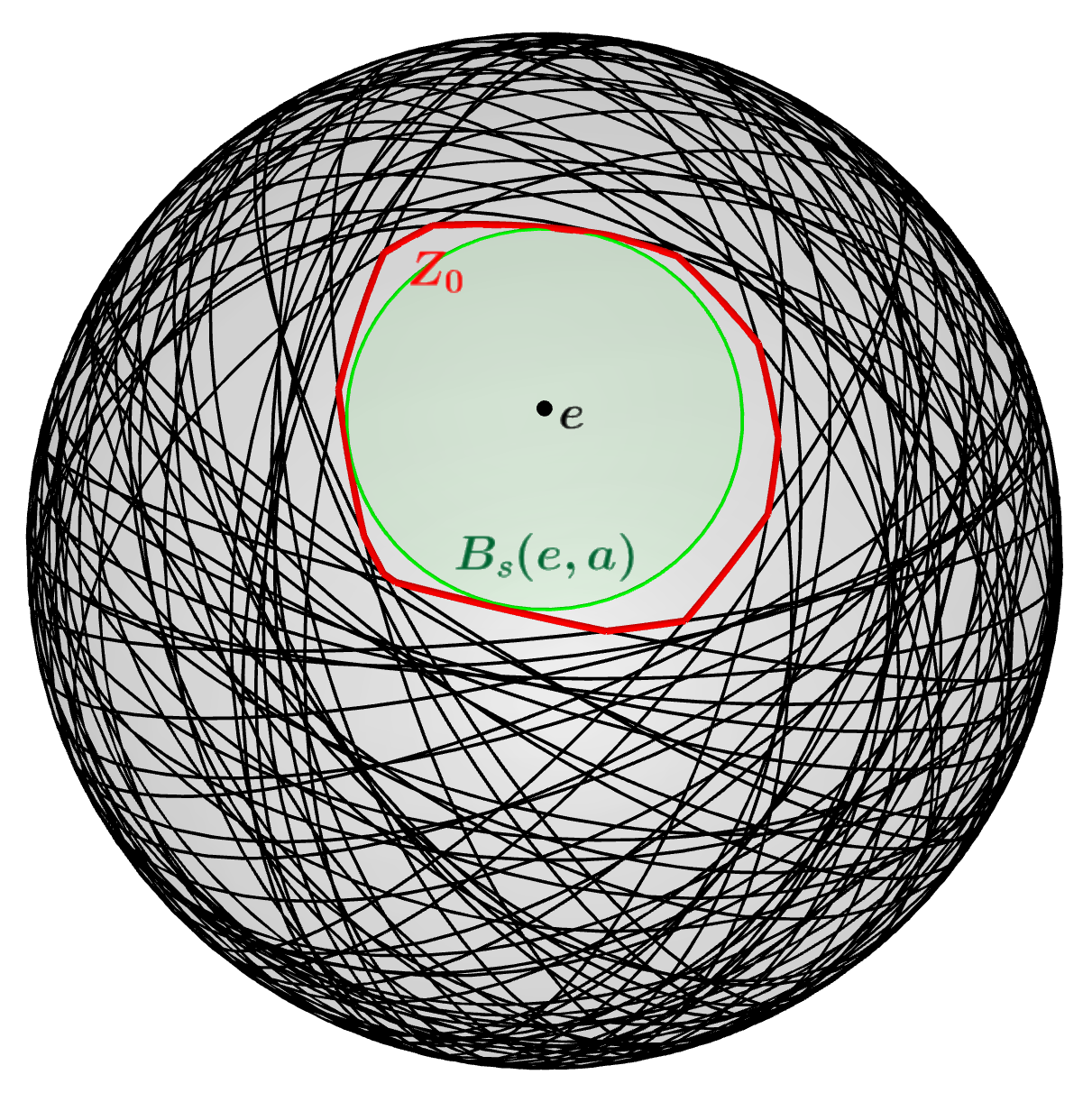}}
  \caption{If the Crofton cell $Z_0$ of $\sfp$ contains a spherical ball $B_s(e,a)$ of radius at least $a$ centred at some point $e\in Z_0$, then $r(Z_0)\ge a$. For easier programming, realizations are   discarded for which $Z_0$ does not contain $B_s(\sfp,a)$. Random points on $\BS^2$ are simulated with R, the figures are produced with GeoGebra.}\label{Fig3}
\end{figure}

By similar arguments, we get the following lemma and the subsequent theorem, which is concerned with the asymptotic shape of Crofton cells
having large spherical inradii in the case where $\Phi=2U_1$. Now the crucial geometric ingredient is Theorem \ref{stab2}. Recall the definition of the size functional $\Sigma_r$ (the inradius functional) and the deviation functional $\vartheta_r$ from Section \ref{subsec:ex2}. An illustration is provided in Fig.~\ref{Fig3}.

\begin{lemma}\label{zaehlerIR}
Let $\Phi=2U_1$. If $a\in (0,\pi/2)$, $ \varepsilon\in(0,1]$ and
$$
\cK_{a,\varepsilon}'\defeq\left\{K\in\overline{\cK}_s^d:\Sigma_r(K)\ge a, \vartheta_r(K)\geq\varepsilon\right\},
$$
then
$$
\BP(Z_0\in\cK_{a,\varepsilon}')\leq c_{4}\,\max\{1,\gamma_s\}^{d\nu}\,\exp\left(-\gamma_s
\left(1+c_5\,\varepsilon^{\frac{d+1}{2}}\right)\Phi(B_s(\sfp,a))\right),
$$
where the constants $c_{4},\nu$ depend on $a,d,\varepsilon$ and $c_5\ge  0.02^da^{d-2}(\pi/2-a)^{d-1}$.
\end{lemma}

\begin{proof}
Using Theorem \ref{stab2} instead of Proposition \ref{propneu1}, the argument is essentially the same as in the
proof of Lemma \ref{zaehler}.
\end{proof}

\begin{theorem}\label{phinball}
If $a\in (0,\pi/2)$ and $\varepsilon\in (0,1]$, then there are constants $c_{6},c_{7}>0$ such that
$$
\BP(\vartheta_r(Z_0)\geq\varepsilon \mid \Sigma_r(Z_0)\geq a)\leq c_{6}\,\exp\left(-c_{7}\,\varepsilon^{\frac{d+1}{2}}\,\gamma_s\right),
$$
where $c_{6}$ depends on $a,d,\varepsilon$ and $c_{7}\ge  0.3\cdot 0.01^da^{d-1}(\pi/2-a)^{d-1}$.
\end{theorem}

\begin{proof}
Combining
$$
\BP(\Sigma_r(Z_0)\geq a)=\BP(\widetilde{X}(\cH_{B_s(\sfp,a)})=0)=\exp\left(-\gamma_s\Phi(B_s(\sfp,a))\right)
$$
with Lemma \ref{zaehlerIR}, the result follows as before by using the lower bound in \eqref{eq:UsBa2} and $2U_1=\Phi$.
\end{proof}

\section{Asymptotic distribution of the size of the Crofton cell}
\label{sec:aymdistr}
Similar to \cite[Theorem 2]{LC}, we determine the asymptotic distribution function of a general increasing and rotation invariant size functional $\Sigma$
of the Crofton cell $Z_0$ of the hypersphere tessellation derived from an isotropic Poisson point process with intensity $\gamma_s$,  as the intensity
$\gamma_s$ tends to infinity. We use the techniques developed in the proof of Lemma \ref{zaehler} to obtain the following theorem.
We write $\tau(a)=\tau(\Phi,\Sigma,a)$ for the associated isoperimetric constant, where $\Phi=2 U_1$.
\begin{theorem}\label{asymp}
Let $\Sigma$ be an increasing and rotation invariant size functional. If $a>0$ is such that $\Sigma^{-1}([a,\infty))\neq\emptyset$, then
$$
\lim_{\gamma_s\rightarrow\infty}\gamma_s^{-1}\, \ln \BP(\Sigma(Z_0)\geq a)=- \tau(\Phi,\Sigma,a).
$$
\end{theorem}

\begin{proof} Since we are only interested in $\gamma_s\to\infty$, we can assume that $\gamma_s \ge 1$.
Let $\kappa\in(0,1)$ and $\cK_{a,0}=\{K\in\overline{\cK}_s^d:\Sigma(K)\geq a\}$.
Let $N\in\mathbb{N}$, $N\ge d+1$, and let $H_1,\ldots,H_N\in \mathcal{H}_{\BS^d}$ be such that $P(H_{(N)})\in\cK_{a,0}\cap
\cK^d_s$,
in particular  $\Phi(P(H_{(N)}))\ge \tau(a)$.  By Lemma \ref{Ecken} we obtain a number
$\nu=\nu(\Sigma,d,a,\kappa)$ and a spherical polytope $Q(P(H_{(N)}))=:Q$ with at most $\nu$ vertices and $\text{ext}(Q)\subset\text{ext}(P(H_{(N)}))$
such that
$$
\Phi(Q)\geq\left(1-\frac{\kappa}{2}\right)\Phi(P(H_{(N)}))\geq\left(1-\frac{\kappa}{2}\right)\tau(a).
$$
Proceeding as in the proof of Lemma \ref{zaehler},  we obtain
$$
\BP(Z_0\in\cK_{a,0} \mid \widetilde{X}(\cH_{\BS^d})=N)\leq\sum_{j=d+1}^{d \nu}\binom{N}{j}
\binom{j}{d}^\nu\left[1-\left(1-\frac{\kappa}{2}\right)\tau(a)\right]^{N-j}
$$
for $N\geq d+1$. After summation over $N$, where we deal with the cases $N\in\{0,\ldots,d\}$ as in the proof of Lemma
\ref{zaehler}, and since $\gamma_s\ge 1$, we get
\begin{align}\label{nenner2}
\BP(Z_0\in\cK_{a,0})&\leq c_8 \gamma_s^{d\nu}\exp\left(-\left(1-\frac{\kappa}{2}\right)\tau(a)
\gamma_s\right)\nonumber\\
&\leq c_{9}\,\exp\left(-\left(1-\kappa\right)\tau(a)\gamma_s\right)
\end{align}
for suitable constants $c_8,c_{9}>0$ which depend only on $\Sigma,a,d$ and $\kappa$. For the last inequality, we used that
$x\mapsto x^{d\nu}\exp(-\kappa/2\, x)$ is bounded.

Combining (\ref{nenner2}) with (\ref{nenner}), we get
$$
\exp\left(-\gamma_s\tau(a)\right)\leq\BP(\Sigma(Z_0)\geq a)\leq c_{9}\,\exp\left(-(1-\kappa)\gamma_s
\tau(a)\right).
$$
This yields
$$
\liminf_{\gamma_s\rightarrow\infty}\gamma_s^{-1}\ln\BP(\Sigma(Z_0)\geq a)\geq -  \tau(a)
$$
and
$$
\limsup_{\gamma_s\rightarrow\infty}\gamma_s^{-1}\ln\BP(\Sigma(Z_0)\geq a)\leq -(1-\kappa)\tau(a).
$$
The left-hand side of the second estimate is independent of $\kappa$ and therefore
$$
\lim_{\gamma_s\rightarrow\infty}\gamma_s^{-1}\ln\BP(\Sigma_d(Z_0)\geq a)=- \tau(a),
$$
which completes this proof.
\end{proof}

\section{Limit shapes}
\label{sec:limshape}

In its original form, Kendall's conjecture was vaguely expressed in terms of a weak convergence statement. Adjusting the approach from \cite{HRS,LC,HS22}, we briefly explain how rigorous results can be obtained in the present setting.

First, we introduce a notion of shape. Let ${\sf G}$ be a subgroup of $\SO(d+1)$. If for given $K,L\in\widehat{\cK}^d_s\defeq\overline{\cK}^d_s\cup\{\BS^d\}$ there is some $g\in {\sf G}$ such that $K=gL$, then we write $K\sim_{\sf G} L$ and say that $K$ and $L$ have the same ${\sf G}$-shape. We endow the quotient space ${\mathcal S}_{\sf G}:= \widehat{\cK}^d_s/\hspace{-4pt}\sim_{\sf G}$ with the induced topology, which is the finest topology such that the canonical projection $s_{\sf G}:\overline{\cK}^d_s\to{\mathcal S}_{\sf G}$ is continuous. For $K\in \overline{\cK}^d_s$, the set $s_{\sf G}(K)=\{gK:g\in{\sf G}\}\subset \overline{\cK}^d_s$  is the class of all spherically convex bodies with the same ${\sf G}$-shape as $K$, and $s_{\sf G}(\BS^d)=\{\BS^d\}$.

Let the situation of Theorem \ref{thm:genresult} be given. Explicitly, suppose that  $Z_0$ is the Crofton cell of the  hypersphere tessellation derived from an isotropic Poisson point process $X$ on $\BS^d$ with intensity $\gamma_s$. Let $\Phi=2U_1$, and let
  $\Sigma$ be an increasing and rotation invariant size functional.
For $a>0$ such that $\Sigma^{-1}([a,\infty))\neq\emptyset$, we define the probability measure
$$
\mu_{a,\gamma_s}\defeq \BP\left(s_{\sf G}(Z_0)\in\cdot\mid \Sigma(Z_0)\ge a\right)
$$
on the Borel $\sigma$-algebra of ${\mathcal S}_{\sf G}$, which is well defined by \eqref{nenner}. We refer to $\mu_{a,\gamma_s}$ as the conditional law of the ${\sf G}$-shape of $Z_0$, given that $\Sigma(Z_0)\ge a$. A shape $s_{\sf G}(B)$, where $B\in \overline{\cK}^d_s $, is said to be the limit shape of $Z_0$ with respect to $\Sigma$ and $a$ as $\gamma_s\to\infty$, if
$$
\lim_{\gamma_s\to\infty}\mu_{a,\gamma_s}=\delta_{s_{\sf G}(B)}\quad\text{weakly},
$$
where $\delta_{s_{\sf G}(B)}$ is the Dirac measure concentrated at the singleton $s_{\sf G}(B)\in {\mathcal S}_{\sf G}$.

The following theorem states that if the extremal bodies in $\mathcal{E}(\Phi,\Sigma,a)$ have a unique ${\sf G}$-shape, then this ${\sf G}$-shape is the limit shape of $Z_0$. Note that $\mathcal{E}(\Phi,\Sigma,a)$ is rotation invariant, if $\Phi$ and $\Sigma$ are rotation invariant.

\begin{theorem}
    \label{thm:limshape}
Let $Z_0$ be the Crofton cell of the hypersphere  tessellation derived from an isotropic Poisson point process $X$ on $\BS^d$ with intensity $\gamma_s$. Let  $\Phi=2U_1$, and let $\Sigma$ be an increasing and rotation invariant size functional. Let  $a>0$ be such that $\Sigma^{-1}([a,\infty))\neq\emptyset$. If  $s_{\sf G}(\mathcal{E}(\Phi,\Sigma,a))=s_{\sf G}(B)$, for some $B\in \overline{\cK}^d_s $,  then $s_{\sf G}(B)$ is the limit shape of $Z_0$ with respect to $\Sigma$ and $a$ as $\gamma_s\to\infty$.
\end{theorem}

\begin{proof}
We will show that if $\mathcal{C}\subset \mathcal{S}_{\sf G}$ is closed, then
\begin{equation}\label{eq:claim}
\limsup_{\gamma_\to\infty}\mu_{a,\gamma_s}(\mathcal{C})\le\delta_{s_{\sf G}(B)}(\mathcal{C}),
\end{equation}
which yields the stated result. For this, let a closed set $\mathcal{C}\subset \mathcal{S}_{\sf G}$ be fixed. Let $\vartheta$ be the canonical deviation functional associated with $\Phi,\Sigma,a$ as in \eqref{devfunc}.

If $s_{\sf G}(B)\in \mathcal{C}$, then \eqref{eq:claim} is clearly true.

Hence suppose that $s_{\sf G}(B)\notin \mathcal{C}$.
If the set $\mathcal{K}_a\defeq\{K\in \overline{\cK}^d_s\cap s_{\sf G}^{-1}(\mathcal{C}):\Sigma(K)\ge a\}$ is nonempty, then $\vartheta$ attains its minimum at some $K_0\in \mathcal{K}_a$, since $\vartheta$ is continuous and $\mathcal{K}_a$ is compact (here we use that $\overline{\cK}^d_s$ is compact,  $s_{\sf G}^{-1}(\mathcal{C})$ is closed and $\Sigma$ is continuous). If $\vartheta(K_0)=0$, then $K_0\in \mathcal{E}(\Phi,\Sigma,a)$ and it follows that $s_{\sf G}(K_0)=s_{\sf G}(B)\notin \mathcal{C}$,
which is in contradiction to $K_0\in s_{\sf G}^{-1}(\mathcal{C})$.

Thus we have shown that there is some $\varepsilon\in (0,1]$ such that $\vartheta(K)\ge \varepsilon$ for all $K\in \mathcal{K}_a$, which will be used at \eqref{eq:lastinequ}.

Since  $\Sigma^{-1}([a,\infty))\neq\emptyset$, there is some $K\in \overline{\cK}^d_s$ such that $\Sigma(K)\ge a$. Using that $\Sigma$ is increasing and $K\subset B_s(e_0,\pi/2)$, for some $e_0\in\BS^d$, we get $a\le \Sigma(K)\le \Sigma(B_s(e_0,\pi/2))\le \Sigma(\BS^d)
$. Since $\Sigma$ is rotation invariant, $\Sigma(B_s(e,\pi/2))\ge a$ for all $e\in\BS^d$. In other words, if $X(\BS^d)=0$, then $\Sigma(Z_0)=\Sigma(\BS^d)\ge a$, and if $X(\BS^d)=1$, then
$\Sigma(Z_0)=\Sigma(B_s(e,\pi/2)\ge a$, where $e$ depends on the realization of $Z_0$.
We thus obtain
\begin{align*}
    \BP(\Sigma(Z_0)\ge a)&\ge\BP(\Sigma(Z_0)\ge a,X(\BS^d)=0)+\BP(\Sigma(Z_0)\ge a,X(\BS^d)=1)\\
    &= \BP(X(\BS^d)=0)+ \BP(X(\BS^d)=1)=(1+\gamma_s)e^{-\gamma_s},
\end{align*}
Now we can conclude that
\begin{align}\label{eq:lastinequ}
    \mu_{a,\gamma_s}(\mathcal{C})&=\BP\left(s_{\sf G}(Z_0)\in\mathcal{C}\mid \Sigma(Z_0)\ge a\right)\nonumber\\
    &\le \BP(Z_0\in \mathcal{K}_a)\mid \Sigma(Z_0\ge a)+\BP(Z_0=\BS^d\mid\Sigma(Z_0)\ge a)\nonumber\\
    &\le \BP(\vartheta(Z_0)\ge\varepsilon \mid \Sigma(Z_0)\ge a)+\frac{e^{-\gamma_s}}{(1+\gamma_s)e^{-\gamma_s}}.
\end{align}
An application of Theorem \ref{thm:genresult} shows that the right-hand side converges to zero as $\gamma_s\to\infty$.
\end{proof}

\begin{examples}\label{ex:secneu8}{\rm
Let $Z_0$ be the Crofton cell of the hypersphere  tessellation derived from an isotropic Poisson point process $X$ on $\BS^d$ with intensity $\gamma_s$. Let
  $\Phi=2U_1$. Then Theorem \ref{thm:limshape} applies, for instance, in the following situations.
\begin{itemize}
\item If $\Sigma=\sigma_d$ (volume), $a\in (0,\omega_{d+1}/2)$, ${\sf G}=\SO(d+1)$, the limit shape of $Z_0$  is the class of spherical caps with spherical Lebesgue measure $a$.
\item If $\Sigma=r$ (inradius), $a\in (0,\pi/2)$, ${\sf G}=\SO(d+1)$, the limit shape of $Z_0$  is the class of spherical caps with radius $a$.
\item Consider $\Sigma=\text{diam}$ (diameter) and $a\in (0,\pi)$. Then one can show that
$$
\min\{\Phi(K):K\in \overline{\cK}^d_s,\text{diam}(K)\ge a\}\ge \frac{a}{\pi}
$$
and $\mathcal{E}(\Phi,\Sigma,a)$ is the set of all geodesic segments of length $a$. With respect to ${\sf G}=\SO(d+1)$, the limit shape of $Z_0$ is the class of geodesic segments of length $a$ if size is measured by the diameter. So far,  we did not provide an explicit stability function in this example. If one chooses
$$
\vartheta_a(K)\defeq\min\{\alpha\ge 0: I\subset K\subset I_\alpha, I\text{ is a geodesic segment of length  }a\},
$$
for $K\in \overline{\cK}^d_s$ with $\text{diam}(K)=a$, then $\vartheta_a(K)\ge \varepsilon$ implies that
$$
\Phi(K)\ge \left(1+c_d\,\varepsilon^2\right)\frac{a}{\pi}.
$$
\end{itemize}
In all these examples, we have ${\sf G}=\SO(d+1)$. For the analysis of typical cells, the subgroup of rotations that keep $\sfp$ fixed is relevant. In this context,
Theorem \ref{thmgenZ} can be used to derive a general result similar to Theorem \ref{thm:limshape}.}
\end{examples}

\section{Typical cells of particle processes in spherical space}
\label{sec:typcells}
After having studied the Crofton cell in the previous section, we now turn to typical cells.
In Euclidean space, there is a very intuitive representation for the distribution of the
typical particle of a stationary particle process. In the special case of a particle process derived from a tessellation, this
leads to the notion of the typical cell of the tessellation. In spherical space, some modifications are required.
Instead of specializing the general framework of \cite{Last, Last2010b}, where
  random measures are studied in a general homogeneous space, we prefer a more direct reasoning which provides additional insights.
  The current approach is similar to the one in hyperbolic space, but some simplifications are possible due to
  the compactness of the unit sphere and its isometry group. However, the compactness of the sphere also requires some additional considerations.

\subsection{Particle processes}\label{subse:PP}
Let $\lambda$ denote the unique Haar probability measure on the compact isometry group $\Ih\defeq\SO(d+1)$ of $\BS^d$ (see \cite[Theorem 9.2.6]{Cohn}), which is left invariant, right invariant, and inversion invariant. The isometry group operates continuously and transitively on $\BS^d$. For each $x\in\BS^d$, the map $\varrho\mapsto \varrho x$ from $\Ih$ to $\BS^d$ is open and proper (inverse images of compact sets are compact), hence $\BS^d$ is a homogeneous $\Ih$-space (see \cite{Nachbin,SW2008}).
%The $d$-dimensional Hausdorff measure (i.e. spherical Lebesgue measure) $\mathcal{H}^{d}=\sigma_d$ on the Borel $\sigma$ field of $\BS^d$ is finite and rotation invariant.
For all $x\in\BS^d$,
\begin{equation}\label{eq:hdlambda}
\sdn=
\int_{\Ih} \1\{\varrho x\in \cdot\}\, \lambda(\dint \varrho).
\end{equation}
Let
$\Ih(\sfp) \defeq \{\varrho \in \Ih \colon \varrho \sfp=\sfp\}$  denote the isotropy group of $\sfp$, which carries a unique $\Ih(\sfp)$ invariant probability measure $\kappa(\sfp, \cdot)$. We
define $\kappa(\sfp,\Ih \setminus \Ih(\sfp))\defeq 0$ to extend
$\kappa(\sfp,\cdot)$ to $\Ih$. In addition, we define $\kappa(o,\cdot)\defeq \lambda$. More generally,
$$\Ih(\sfp,x) \defeq \{ \varrho \in \Ih \colon  \varrho \sfp=x\},\quad x\in\BS^d,$$
is the set of isometries that map $\sfp$ to $x$. Choosing an arbitrary
$\varrho_x \in \Ih(\sfp,x)$, we define
\begin{align}\label{e:kappa}
\kappa(x,B) \defeq \int_{\Ih(\sfp)} \1 \{\varrho_x \circ \varrho \in B\} \, \kappa(\sfp,\dint \varrho), \quad x\in\BS^d, B \in \mathcal{B}(\Ih).
\end{align}
This definition is independent of the choice of $\varrho_x$ (see
\cite{Last}) and the map $x\mapsto \varrho_x$ can be chosen in a
measurable way. Hence $\kappa$ is a stochastic transition kernel from
$\BS^d$ to $\Ih$. Moreover, $\kappa(x,\cdot)$ is concentrated on
$\Ih(\sfp,x)$.
The kernel $\kappa$ disintegrates $\lambda$ in the sense that
\begin{align}\label{disintlambda}
\int_{\BS^d}\int_{\Ih} \I\{\varrho\in\cdot\}\,\kappa(x,\dint \varrho)\,\sdn(\dint x)=\lambda.
\end{align}
A simple direct argument for \eqref{disintlambda} is given in \cite[Sec.~2.1]{HLS24}.

 Let $\cC^d$ denote the space of nonempty compact
  subsets of $\BS^d$, endowed with the Fell topology (that coincides with the topology induced by the Hausdorff metric \cite[Theorem 12.3.2]{SW2008}), which is a locally compact Hausdorff space with countable base of the topology (see the remarks after \cite[Theorem 12.2.1]{SW2008}).
By a particle process on $\BS^d$ we mean a point process $\zeta$ on  $\cC^d$.
We require $\zeta$ to be finite, that is, $\zeta(\cC^d)<\infty$.
In other words, $\zeta$ is a measurable mapping
from $\Omega$ to $\mathbf{N}$, where
$(\Omega,\mathcal{F},\BP)$ is the underlying probability space
and $\mathbf{N}$ is the space of all measures on
$\cC^d$ which take values in $\N_0$. Measurability refers to the smallest
$\sigma$-field on $\mathbf{N}$ such that the mappings $\mu\mapsto\mu(A)$
(from $\mathbf{N}$ to $[0,\infty)$) are measurable, for
each Borel set $A\subset \cC^d$.
Later we will consider a particle process $\zeta$ concentrated on
the space $\widehat{\cK}^d_s=\overline{\cK}_s^d\cup\{\BS^d\}$ of nonempty compact convex subsets of
$\BS^d$, that is, $\BP(\zeta(\cC^d\setminus \widehat{\cK}^d_s)=0)=1$.
The isometry group $\Ih$ acts continuously on $\cC^d$.
For $\mu\in\mathbf{N}$
and $\varrho\in \Ih$, let $\varrho \mu$ denote the image measure of
$\mu$ under $\varrho$. The map $(\varrho,\mu)\mapsto \varrho\mu$
is measurable. A particle process $\zeta$ is called  {\em isotropic}
if $\varrho\zeta\overset{d}{=}\zeta$, for each
$\varrho\in \Ih$.

Suppose that $\zeta$ is an isotropic particle process on $\BS^d$.
Then the {\em intensity measure} $\Lambda\defeq\BE\zeta$
is isometry (rotation) invariant, that is,
$\Lambda(A)=\Lambda(\varrho A)$, for   $\varrho\in \Ih$
and  measurable sets $A\subset \cC^d$,
where $\varrho A\defeq\{\varrho C\colon C\in A\}$. If $\zeta$ is a Poisson process,
then the distribution of $\zeta$ is completely determined by $\Lambda$;
see \cite[Proposition 3.2]{LastPenrose}. In this case,
$\Lambda$ is finite.
In the following, we establish an important disintegration result for a finite isotropic measure $\Lambda$ on $\cC^d$.
For this purpose, we fix a {\em centre function} by which we mean a measurable map
$c_s\colon \cC^{d} \rightarrow \BS^d\cup\{o\}$ that is isometry covariant in the
sense that
\begin{align}\label{eq:c_h_property}
c_s(\varrho C)= \varrho c_s(C), \quad C \in \cC^{d},  \varrho \in \Ih.
\end{align}
If $C=\BS^d$, then \eqref{eq:c_h_property} requires that $c_s(\BS^d)=\varrho c_s(\BS^d)$ for all $\varrho \in \Ih$, hence $c_s(\BS^d)=o$.
To provide an explicit example of a centre function, let $\varphi(C)\defeq \text{cl}\conv(C)$, $C\in \cC^d$, be the closed convex hull of $C$, which is a spherically convex body or the intersection of $\BS^d$ with a linear subspace of $\R^{d+1}$. If $\psi(C)\defeq (\varphi(C)\cap (-\varphi(C)))^\perp\cap \varphi(C)\neq\emptyset$, then $\psi(C)$ is a proper convex body and we define $c_s(C)\defeq R(\psi(C))$, which is the uniquely determined centre of the circumball of $\psi(C)$. Otherwise, we define $c_s(C)=o$. In this example, we have $c_s( C)\in C$ for $C \in \overline{\cK}^{d}_s$.
The number $\gamma=\BE[\zeta(\cC^{d})]=\Lambda(\cC^{d})$ is said to be the {\em intensity} of $\zeta$.
We assume that $\gamma$ is positive and finite.

Let
$$
\cC_{z}^{d} \defeq \{C \in \cC^{d}\colon \ c_s(C)=z\},\quad z\in \{o,\sfp\},
$$
be the set of all compact sets with centre at $z$, and $\cC^d_{o,\sfp}\defeq \cC_{o}^{d} \cup \cC_{\sfp}^{d} $. Clearly, $\cC^d_{o}$ and $\cC^d_*\defeq \cC^d\setminus \cC^d_o$ are $\Ih$ invariant sets and $\cC^d_{\sfp}$ is an $\Ih(\sfp)$ invariant set. We write $\mu\resmes A$ for the restriction of a measure $\mu$ to a set $A$. Theorem \ref{thmdisintegration}  is stated for an isometry invariant, non-zero, finite Borel measure $\Lambda$ on $\cC^d$. Hence it applies in particular to the intensity measure of an isotropic particle process on $\BS^d$.

\begin{theorem}\label{thmdisintegration}
  Let $\Lambda$ be an isometry invariant finite Borel measure on
$\cC^d$ with $\gamma\defeq \Lambda(\cC^d)\in (0,\infty)$. Let $c_s\colon \cC^{d} \rightarrow \BS^d\cup\{o\}$ be a centre
  function.  There exists a unique probability measure
$\BQ$ on $\cC^d$, concentrated
on $\cC^d_{o,\sfp}$, such that $\BQ\resmes \cC^d_{o}$ is rotation invariant, $\BQ\resmes \cC^d_{\sfp}$ is invariant under rotations fixing $\sfp$, and $\BQ$  satisfies
\begin{equation}\label{eqdisintegration}
\Lambda =\gamma \int_{\cC^{d}}\int_{\Ih} \I\{\varrho C\in\cdot\} \,\lambda(\dint \varrho) \,\BQ(\dint C).
   \end{equation}
   This measure is given by
\begin{align}\label{e:BQ}
\BQ = \frac{1}{\gamma} \int_{\cC^{d}} \int_{\Ih}
 \1 \{\varrho^{-1} C \in \cdot\} \, \kappa(c_s(C),\dint \varrho) \, \Lambda( \dint C).
\end{align}
\end{theorem}

The proof of Theorem \ref{thmdisintegration} is similar to the argument for \cite[Theorem 2.1]{HLS24}. Some simplifications are  due to the fact that the measures involved are finite and a constant weight function $w$ can be used, however some additional case distinctions are required. Another variant of the arguments is provided in detail in the proof of Theorem \ref{thm:typVch} below.  This is also true for the following remark.

\begin{remark}{\rm
Let $\overline{\kappa}$ be an arbitrary probability kernel from $\BS^d\cup\{o\}$ to $\Ih$. If $\overline{\kappa}(x,\cdot)$ is concentrated on $\Ih(\sfp,x)$, for $\mathcal{H}^d$-almost all $x\in\BS^d$, then
\begin{align}\label{e:BQ2a}
 {\BQ} = \frac{1}{\gamma} \int_{\cC^{d}} \int_{\Ih}
 \1 \{\varrho^{-1} C \in \cdot\} \, \overline{\kappa}(c_s(C),\dint \varrho) \, \Lambda( \dint C).
\end{align}
This shows that while $\kappa$ as given at \eqref{e:kappa}, together with $\kappa(o,\cdot)=\lambda$, is a natural choice, any other kernel $\overline{\kappa}$ leads to the same measure. For instance, we can choose $\overline{\kappa}(x,\cdot)=\delta_{\varrho_x}$ provided that $\varrho_x\in \Ih(\sfp,x) $ is a measurable function of $x\in \BS^d$ and $\varrho_o\in\Ih$ is arbitrary.
}
\end{remark}

\begin{remark}{\rm
If $\Lambda=\mathbb{E} \zeta$ is the intensity measure of an isotropic particle process $\zeta$ and $\gamma=\Lambda(\cC^{d})\in (0,\infty)$, then the probability measure $\BQ$ given in \eqref{e:BQ} (and \eqref{e:BQ2a})
is called the distribution of the {\em typical particle} (of $\zeta$).
It is convenient to introduce a random element
$\TG$ with this distribution. Roughly speaking,  the typical particle
is obtained by choosing one of the particles $C$ of $\zeta$ at random
and applying to it a rotation $\varrho^{-1}$ with $\varrho\in \Ih(\sfp, c_s(C))$, so that $c_s(\varrho^{-1} C)=\sfp$, if $c_s(C)\in\BS^d$; if $c_s(C)=o$, then we take a uniform random rotation of $C$. While in the Euclidean case a canonical selection of an isometry mapping $\sfp$ to  $c_s(C)$ is available, this is
not the case in spherical space.}
\end{remark}

\subsection{The typical cell of Poisson hypersphere  tessellations}
\label{sec:typcellPois}
As in \cite{AZ}, we can interpret the tessellation generated by the isotropic hypersphere process $\widetilde{X}$ as an isotropic
particle process $X'$ on $\widehat{\cK}^d_s$ with intensity $\gamma_{X'}\in (0,\infty)$. The distribution of the typical particle
$Z$ of $X'$ is then given by
\begin{equation*}%\label{typpart}
\BP_{Z}(\cdot)\defeq\BP(Z\in\cdot)
=\frac{1}{\gamma_{X'}}\mathbb{E}\,\int_{\widehat{\cK}^d_s}\int_{\Ih}\mathbf{1}\{\varphi^{-1}K\in\cdot\}\, \kappa(c_s(K),\dint \varphi)\, X'(\dint K).
\end{equation*}
The following relation
between the typical cell and the Crofton cell of an isotropic tessellation on $\BS^d$ is a special case of a well-known
relationship valid for tessellations in homogeneous spaces (see, e.g., \cite[Corollary 8.4]{Last}).
Its Euclidean counterpart can be found in \cite[Theorem 10.4.1]{SW2008}. We include the simple proof for convenience and add an
explicit expression for the intensity of $X'$ if the tessellation is induced by a Poisson hypersphere
process. In advance, we point out some properties of the functions
$$
h_m: [0,\infty)\rightarrow\R,\qquad
  t\mapsto(-1)^{m+1}e^{-t}+2\sum_{i=0}^{\lfloor\frac{m}{2}\rfloor}\frac{t^{m-2i}}{(m-2i)!},
$$
for $m\in\mathbb{N}_0$, which will occur in the explicit expression of the intensity of $X'$.
\begin{lemma}
The functions $h_m,\ m\in\N_0,$ have the following properties:
\begin{itemize}
\item[$(1)$] $h_m'=h_{m-1},\ m\geq1$;
\item[$(2)$] $h_0(t)=2-e^{-t}\geq1,$ $h_1(t)=e^{-t}+2t\geq1$;
\item[$(3)$] $h_m(0)=1,\ m\geq0$;
\item[$(4)$] $h_m$ is strictly increasing and $h_m\geq1,$ $m\geq0$;
\item[$(5)$] $h_m$ is convex for $m\geq 1$;
\item[$(6)$] $0\leq h_m(t)-\left(1+\frac{t}{1!}+\frac{t^2}{2!}+\ldots+\frac{t^m}{m!}\right)\leq\frac{t^m}{m!},\ m\in\mathbb{N}_0$.
\end{itemize}
\end{lemma}

\begin{proof}
First, note that for $m\in\N_0$
$$
h_m(t)=(-1)^{m+1}e^{-t}+\sum_{i=0}^m\left(1+(-1)^{m-i}\right)\frac{t^i}{i!}.
$$
The properties (1), (2) and (3) now follow directly from this identity, (4) and (5) are an immediate consequence of the
first three properties.
Assertion (6) follows easily by induction over $m$.
\end{proof}

The auxiliary function $h_d$ is used to relate the intensities $\gamma_{X'}$ and $\gamma_s$ in the following theorem.

\begin{theorem} \label{NC}
Let $f:\widehat{\cK}_s^d\rightarrow[0,\infty)$ be measurable and rotation invariant. Let $X'$ be an isotropic tessellation of $\BS^d$
with intensity $\gamma_{X'}\in (0,\infty)$.  If $Z_0$ denotes the spherical Crofton cell and $Z$ the typical cell of $X'$, then
\begin{equation}\label{typ0}
\BE\left[f(Z_0)\right]=\gamma_{X'}\BE\left[f(Z)\,\sdn(Z)\right].
\end{equation}
If $X'$ is the hypersphere tessellation induced by $\widetilde{X}=h(X)$, where $X$ is a Poisson point process on $\BS^d$  whose intensity measure has no point masses, then $\BE[ X'\left(\mathcal{C}^d\right)]=h_d\left(\BE[ X(\BS^d)]\right)$. In particular, if $X$ is isotropic with intensity $\gamma_s$,
then $\gamma_{X'}=h_d(\gamma_s) $.
\end{theorem}

\begin{proof}
From (\ref{eqdisintegration}) with $\Lambda=\BE X'$ and $\BQ=\BP_Z$  and from the rotation invariance of $f$ we get
\begin{align*}
\BE[f(Z_0)]&=\BE\sum_{K\in X'} f(K)\mathbf{1}\{\sfp\in\interior(K)\}\\
&=\gamma_{X'}\int_{\widehat{\cK}_{\sfp}^d}\int_{\Ih}f(\varphi K)\mathbf{1}\{\sfp\in\interior(\varphi K)\}\, \lambda(\dint \varphi)\, \BP_Z(\dint K)\\
&=\gamma_{X'}\int_{\widehat{\cK}_{\sfp}^d}f(K)\int_{\Ih}\mathbf{1}\{\varphi\sfp\in\interior(K)\}\, \lambda(\dint \varphi)\, \BP_Z(\dint K)\\
&=\gamma_{X'}\int_{\widehat{\cK}_{\sfp}^d}f(K)\,\sdn(K)\, \BP_Z(\dint K),
\end{align*}
where \eqref{eq:hdlambda} and $\sdn(\partial K)=0$ for $K\in \widehat{\cK}^d_s$ were used in the last step.

For the second part (where $X$ is a Poisson process), we use Schl\"afli's theorem (see \cite[p.~209--212]{schlaf} or \cite[(1.1)]{CE} in modern language),
which provides an explicit formula for the number of cells $N(k)$ generated by $k\geq1$ hyperspheres in general position,
$$
N(k)=2\sum_{i=0}^{d}\binom{k-1}{i}.
$$
Recall that
$h(x)=\BS^{d}\cap x^\perp=\BS_x$, $x\in\BS^d$. Since the intensity measure of $X$ is diffuse, $X$ is simple and if $x\in X$, then $-x\notin X$, almost surely (by \cite[Theorem 4.4]{LastPenrose} or \cite[Corollary 3.2.3]{SW2008}).
  If $X$ contains no points, we consider the whole of $\BS^d$ as one cell and thus define $N(0)\defeq 1$. Therefore we get
\begin{align}\label{gamma}
\BE\left[ X'(\mathcal{C}^d)\right]&= \BE\left[N(X(\BS^d))\right] \nonumber\\
&= \sum_{k=1}^\infty2\,\sum_{i=0}^d\binom{k-1}{i}\, \BP(X(\BS^d)=k)+ \BP(X(\BS^d)=0)\nonumber\\
&=2\,\sum_{i=0}^d\frac{e^{-\BE[ X(\BS^d)]}}{i!}\,\sum_{k=i+1}^\infty
\frac{\BE[ X(\BS^d)]^k}{(k-i-1)!}\,\frac{1}{k}+ e^{-\BE[ X(\BS^d)]}.
\end{align}
The relation
\begin{equation}\label{eq:relneu}
\frac{1}{i!}\sum_{k=i+1}^\infty\frac{1}{k}\cdot\frac{x^k}{(k-i-1)!}=e^x\sum_{k=0}^i(-1)^k\frac{x^{i-k}}{(i-k)!}+(-1)^{i+1},\quad x\in[0,\infty),
\end{equation}
holds, since both sides of relation \eqref{eq:relneu} are zero for $x=0$ and the differentials are equal. Plugging \eqref{eq:relneu} with $x=\BE[ X(\BS^d)]$ into \eqref{gamma} and observing that $\sum_{i=0}^d(-1)^{i+1}=-\frac{1}{2}(1+(-1)^d)$, we get
\begin{align*}
\BE\left[ X'(\mathcal{C}^d)\right]
&= (-1)^{d+1}e^{-\BE[ X(\BS^d)]} + {2}\sum_{i=0}^d\sum_{k=0}^i\left[(-1)^k\,
 \frac{\BE[ X(\BS^d)]^{i-k}}{(i-k)!}\right].
\end{align*}
Since
$$
\sum_{i=0}^d\sum_{k=0}^i (-1)^k\frac{x^{i-k}}{(i-k)!} =\sum_{i=0}^{\lfloor\frac{d}{2}\rfloor}\frac{x^{d-2i}}{(d-2i)!},
$$
the remaining assertion follows.
\end{proof}

Now we are able to extend our results concerning the asymptotic shape of the spherical Crofton cell to typical cells.
We state  an abstract version of such a result, for specific size functionals the argument can be adjusted as before.
Again we require that the size functional $\Sigma$ be increasing, but now we also need that $\Sigma$ is rotation invariant and simple. The latter means that  $\Sigma(K)=0$
whenever $K$ is not $d$-dimensional. This condition is clearly satisfied by the volume and the inradius.  A simple compactness and continuity argument shows that if $a >0$ with
$\Sigma^{-1}([a,\infty))\neq\emptyset$, then the condition implies that there is
a positive constant $c(a)>0$ such that $\sdn(K)\ge c(a)$ whenever $\Sigma(K)\ge a$. Alternatively, one can assume that $\sdn$ can be bounded from above by an increasing (positive) function of $\Sigma$.

\begin{theorem}\label{thmgenZ}
Let $\widetilde X$ be the hypersphere process on $\BS^d$ derived from an isotropic Poisson point process $X$ on $\BS^d$ with intensity $\gamma_s\in (0,\infty)$.  Let $Z$ be the typical cell
of the induced tessellation. Let $\Phi=2U_1$, let $\Sigma$ be an increasing, rotation invariant, simple size functional, and let $a>0$ be such that $\Sigma^{-1}([a,\infty))\neq\emptyset$. Let $\vartheta$ be a rotation invariant deviation functional for $\Phi,\Sigma,a$, and let $f_a$ be a stability function as in Proposition \ref{propneu1}.
If $\varepsilon\in(0,1]$ and $\kappa\in (0,1)$,
then
$$
\BP(\vartheta(Z)\geq\varepsilon \mid \Sigma(Z)\geq a)\leq c_{10}\,\exp\left(-\kappa\gamma_sf_a(\varepsilon)\tau(a)\right),
$$
where the constant $c_{10}>0$ depends on $a,d,\varepsilon,\kappa,\Sigma,\vartheta$ and $\tau(a)=\tau(\Phi,\Sigma,a)$.
\end{theorem}

\begin{proof}
We first note the trivial upper bound $\sdn(Z_0)\leq 1$. In order to estimate the denominator, we use
(\ref{typ0}) and (\ref{nenner}) to obtain
\begin{align*}
\BP(\Sigma(Z)\geq a)&=\BE\left[\mathbf{1}\{\Sigma(Z)\geq a\}\right]
=\frac{1}{\gamma_{X'}}\,\gamma_{X'}\,\BE\left[\mathbf{1}\{\Sigma(Z)\geq a\}\,\frac{1}{\sdn(Z)}\,\sdn(Z)
\right]\\
&=\frac{1}{\gamma_{X'}}\,\BE\left[\mathbf{1}\{\Sigma(Z_0)\geq a\}\,\frac{1}{\sdn(Z_0)}\right]\\
&\geq\frac{1}{\gamma_{X'}}\,\BE\left[\mathbf{1}\{\Sigma(Z_0)\geq a\}\right]\\
&\geq\frac{1}{\gamma_{X'}}\,\exp\left(-\gamma_s\tau(a)\right),
\end{align*}
where $\tau(a)=\tau(\Phi,\Sigma,a)$. For the numerator, we use (\ref{typ0}), (\ref{phszaehler}), the fact that
 $\Phi$, $\Sigma$ and $\vartheta$ are rotation invariant (by assumption), and proceed as above to get
\begin{align*}
\BP(\Sigma(Z)\geq a,\ \vartheta(Z)\geq\varepsilon)&=\frac{1}{\gamma_{X'}}\,\BE\left[\mathbf{1}\{\Sigma(Z_0)\geq a,
\, \vartheta(Z_0)\geq\varepsilon\}\,\frac{1}{\sdn(Z_0)}\right]\\
&\leq\frac{1}{\gamma_{X'}}\,\BE\left[\mathbf{1}\{\Sigma(Z_0)\geq a,\, \vartheta(Z_0)\geq\varepsilon\}\,
\frac{1}{c(a)}\right]\\
&\leq\frac{c_1}{c(a)\,\gamma_{X'}}\max\{1,\gamma_s\}^{d \nu}\,\exp\left(-\gamma_s
\left(1+(1-\bar\kappa)f_a(\varepsilon)\right)\tau(a)\right),
\end{align*}
where $\bar\kappa\defeq (1-\kappa)/2\in (0,1/2)$.
Combining the preceding two estimates, we obtain the  result as in the proof of Theorem \ref{thm:genresult}.
\end{proof}

A result similar to Theorem \ref{asymp} holds for the asymptotic distribution of the typical cell under the assumptions
of Theorem \ref{thmgenZ}.

\begin{theorem}\label{asymptyp}
Let $a>0$ be such that $\Sigma^{-1}([a,\infty))\neq\emptyset$.  Then
$$
\lim_{\gamma_s\rightarrow\infty}\gamma_s^{-1}\, \ln \BP(\Sigma(Z)\geq a)=-\tau(\Phi,\Sigma,a).
$$
\end{theorem}

\section{Spherical Poisson--Voronoi tessellations}
\label{sec:sphPVT}

After the investigation of the spherical Crofton cell and the typical cell of Poisson hypersphere tessellations,
it is natural to explore  Poisson--Voronoi tessellations in spherical space.

Let $A\subset\BS^d$ be finite. The Voronoi cell associated with $x\in A$ is given by
$$
C(x,A)\defeq \left\{y\in\BS^d:\ d_s(y,x)\leq d_s(y,z)\text{ for all }z\in A\right\}.
$$
If  $A$ is a singleton, then $C(x,A)=\BS^d$ for $x\in A$.
The set of all these cells $C(x,A)$, $x\in A$, forms the Voronoi tessellation generated by $A$. Note that for $x\in A$ and $\varphi\in \Ih$ we have
\begin{equation}\label{rcov}
C(\varphi x,\varphi A)=\varphi C(x,A).
\end{equation}
Let $X$ be an isotropic point process on $\BS^d$ with intensity measure
$
\BE[X(\cdot)]=\gamma_s\,\sdn(\cdot)
$.
As before, we assume that $\gamma_s\in(0,\infty)$.

For an isotropic Voronoi tessellation there is a very natural way to choose a centre function, applied to the cells of the tessellation, namely the nucleus $x\in X$ of the Voronoi cell $C(x,X)$. In analogy to \eqref{e:BQ2a} we therefore define the distribution $\BP_{\ZV}$ of the typical cell $\ZV$ of the Voronoi tessellation associated with $X$ by
\begin{align}\label{eq:deftypV}
\BP_{\ZV}&\defeq \frac{1}{\gamma_s}\BE \int_{\BS^d}\int_{\Ih}\indi\{\varphi^{-1} C(x,X)\in\cdot\}\,\overline{\kappa}(x,\dint \varphi)\, X(\dint x)\nonumber\\
&\;=\frac{1}{\gamma_s}\BE \int_{\BS^d}\int_{\Ih}\indi\{ C(\sfp,\varphi^{-1}X)\in\cdot\}\,\overline{\kappa}(x,\dint \varphi)\, X(\dint x),
\end{align}
where $\overline{\kappa}$ is a probability kernel from $\BS^d$ to $\Ih$ such that $\overline{\kappa}(x,\cdot)$ is concentrated on $\Ih(\sfp,x)$ for $\mathcal{H}^d$-almost all $x\in\BS^d$. The definition is independent of the special choice of the kernel (see Remark \ref{rem:sec10char}), but the kernel $\kappa$ from \eqref{e:kappa} is certainly a natural choice. In the case of an isotropic Poisson process $X$, relation \eqref{eq:indtypVor} in particular confirms that the definition is independent of the particular choice of the kernel $\overline{\kappa}$ and also provides a useful geometric characterization of the typical cell. A characterization of the distribution of the typical cell in the spirit of Theorem \ref{thmdisintegration} is provided in Section \ref{sec:typVor}.
  A different approach to arrive at the preceding definition of the typical cell is described in \cite[Sec.~10]{HR17} and the literature cited there.

If $X$ is an isotropic Poisson point process on $\BS^d$ with $\BE [X(\cdot)]=\gamma_s\,\sdn(\cdot)$, then an application of Mecke's theorem (see \cite[Theorem 4.1]{LastPenrose} or \cite[Theorem 3.2.5]{SW2008}), Fubini's theorem and the isotropy of $X$ yield
\begin{align}\label{eq:indtypVor}
\BP_{\ZV}&=\frac{1}{\gamma_s}\BE \int_{\BS^d}\int_{\Ih}\indi\{ C(\sfp,\varphi^{-1}X)\in\cdot\}\,\overline{\kappa}(x,\dint \varphi)\, X(\dint x)\nonumber\\
&=\int_{\BS^d}\int_{\Ih}
\BE\left[\mathbf{1}\{C( \sfp,\varphi^{-1}(X+\delta_x))\in\cdot\}\right]\, \overline{\kappa}(x,\dint \varphi)\,\sdn(\dint x) \nonumber\\
&=\int_{\BS^d}\int_{\Ih}
\BE\left[\mathbf{1}\{C( \sfp, X+\delta_\sfp)\in\cdot\}\right]\, \overline{\kappa}(x,\dint \varphi)\,\sdn(\dint x) \nonumber\\
&=\BP(C(\sfp,X+\delta_{\sfp} )\in\cdot),
\end{align}
since $\overline{\kappa}$ and $\sdn$ are probability measures.

This relation shows that the typical cell of a spherical Poisson--Voronoi tessellation is the spherical Crofton cell
of a special spherical Poisson hyperplane tessellation. The underlying (a.s.~simple) spherical hyperplane process $Y$ is supported by the set of all great
subspheres, having equal spherical distance to the spherical origin $\sfp$ and to a point $x\in X$. (If $X(\omega)=\emptyset$, then
$C(\sfp,X(\omega)+\delta_{\sfp}) =\BS^d$.)  This suggests to consider the
non-isotropic spherical hyperplane process
$$
Y\defeq\sum_{x\in X}\delta_{(x-\sfp)^\bot\cap\BS^d}=\int_{\BS^d}\mathbf{1}\{(x-\sfp)^\bot\cap\BS^d\in\cdot\}\, X(\dint x),
$$
which is a spherical Poisson hyperplane process with intensity measure
$\BE[Y(\cdot)]=\gamma_s\,\tilde{\mu}(\cdot)$.

In order to derive a deviation result, we start with the proof of an upper bound for the deviation probability.
Recall from Section \ref{subsec:ex3} that $\vartheta_{\sfp}(K)=R_\sfp(K)-r_\sfp(K)$ for $K\in\overline{\cK}^d_\sfp$ and $\vartheta_\sfp(\BS^d)\defeq 0$.

%
%%%%%%%%%%%%%%%%%%%%%%%%%%%%%%%%%%%%%%%%%%%%
%%%%%%%%%%%%%%%%%%%%%%%%%%%%%%%%%%%%%%%%%%%
%
\begin{lemma}\label{pvzaehler}
Let $0<a<\pi/2$ and let $\varepsilon\in (0,1]$. Let $
\widetilde{\cK}_{a,\varepsilon}\defeq\{K\in\overline{\cK}^d_\sfp:r_{\sfp}(K)\ge a, \vartheta_{\sfp}(K)\geq\varepsilon\}$. If $X$ is an isotropic
Poisson process with intensity $\gamma_s$, then
\begin{equation*}%\label{pvzaehlereq}
\BP(\ZV\in\widetilde{\cK}_{a,\varepsilon})\leq c_{11}\,
\max\{1,\gamma_s\}^{d \nu}\,\exp\left(-\gamma_s\left(1+  \tfrac{3}{4}c(a,d)\,\varepsilon^{\frac{d+1}{2}}\right)
\,\sdn(B_s(\sfp,2a))\right),
\end{equation*}
where the constants $\nu,c_{11}$ depend on $a,d,\varepsilon$.
\end{lemma}

\begin{proof}
Let $\widetilde{\Phi}_\sfp(\cdot)\defeq \widetilde{U}_\sfp(\cdot)$ and
$N\in\mathbb{N}$. We can proceed as in the proof of Lemma \ref{zaehler}, if we replace $\mu$ by $\tilde\mu_\sfp$ and $\Phi$ by $\widetilde{\Phi}_\sfp$.

For $H_1,\ldots,H_N\in \mathcal{H}_{\BS^d},$ we define $H_{(N)}\defeq(H_1,\ldots,H_N)$ and let $P(H_{(N)})$ denote the spherical
Crofton cell of the tessellation induced by $H_1,\ldots,H_N$. For $N\ge d+1$, let
 $H_1,\ldots,H_N$ be such that
$P(H_{(N)})\in\widetilde{\cK}_{a,\varepsilon}\cap\cK^d_s$ and define $\alpha\defeq\frac{1}{4}c(a,d)\varepsilon^{\frac{d+1}{2}}/(1+c(a,d)\varepsilon^{\frac{d+1}{2}})$ (that is, we choose $\bar\kappa=1/4$), where the constant $c(a,d)$ is taken from Theorem \ref{stab3}. Then we obtain
  $(1-\alpha)(1+c(a,d)\varepsilon^{\frac{d+1}{2}})=1+(1-\bar\kappa)c(a,d)\varepsilon^{\frac{d+1}{2}}$.
By a result analogous to Lemma \ref{Ecken} and by Theorem \ref{stab3}, there are at most
$\nu=\nu(a,d,\varepsilon)$ vertices of $P(H_{(N)})$ such that the spherical convex hull $Q(H_{(N)})$ of these vertices
satisfies
\begin{align*}
1\ge \widetilde{\Phi}_\sfp(Q(H_{(N)}))&\geq(1-\alpha)\widetilde{\Phi}_\sfp(P(H_{(N)}))\\
&\geq(1-\alpha)(1+c(a,d)\varepsilon^{\frac{d+1}{2}})\widetilde{\Phi}_\sfp(B_s(\sfp,a))
=\left(1+\tfrac{3}{4}c(a,d)\varepsilon^{\frac{d+1}{2}}\right)\widetilde{\Phi}_\sfp(B_s(\sfp,a)) .
\end{align*}
It follows from \eqref{starast} that $\tilde{\mu}^N_\sfp$-almost surely any
$N$ hyperspheres are in general position and therefore $P(H_{(N)})\in \cK^d_s$  holds almost surely for $N\ge d+1$.

 Thus, proceeding as in the proof of Lemma \ref{zaehler}, we obtain
$$
\BP(\ZV\in\widetilde{\cK}_{a,\varepsilon})\leq c_{11}\,\max\{1,\gamma_s\}^{d\nu}\,\exp
\left(-\gamma_s\left(1+\tfrac{3}{4}c(a,d) \,\varepsilon^{\frac{d+1}{2}}\right)\,\widetilde{\Phi}(B_s(\sfp,a))\right).
$$
Since $\widetilde{\Phi}_\sfp(B_s(\sfp,a))=
\widetilde{U}_\sfp(B_s(\sfp,a))= \sdn(B_s(\sfp,2a))$,
this completes the proof.
\end{proof}

For the spherical inball radius, centred at $\sfp$, as our size functional, and for $0<a<\pi/2$, we get
$$
\BP(r_{\sfp}(\ZV)\geq a)=\exp\left(-\gamma_s\tilde{\mu}_\sfp(\cH_{B_s(\sfp,a)})\right)
=\exp\left(-\gamma_s\,\sdn(B_s(\sfp,2a))\right).
$$
In combination with Lemma \ref{pvzaehler} and using again \eqref{lubounds},
we obtain the following theorem for the asymptotic shape of the typical cell
of a sphericall Poisson--Voronoi tessellation. In addition, we describe the asymptotic distribution of the centred
inball radius.

\begin{theorem}\label{pvziel}
Let $0<a<\pi/2$ and $\varepsilon\in(0,1]$. If $X$ is an isotropic Poisson
process on $\BS^d$ with intensity $\gamma_s\in (0,\infty)$, then the typical cell $\ZV$ of the induced spherical
Poisson--Voronoi tessellation satisfies
\begin{equation*}
\BP(\vartheta_{\sfp}(\ZV)\geq\varepsilon \mid r_{\sfp}(\ZV)\ge a)\leq
c_{12}\,\exp\left(- c_{13}\,\varepsilon^{\frac{d+1}{2}}\,\gamma_s\right),
\end{equation*}
where  the constant $c_{12}>0$ depends  on $a,d,\varepsilon$ and
$$c_{13}\ge (1/2)c(a,d)\sdn(B_s(\sfp,2a))\ge \pi^{-4d} a^{d-1}
(\min\left\{a,\tfrac{\pi}{2}-a\right\})^{d-1}.
$$
Moreover,
$$
\lim_{\gamma_s\rightarrow\infty}\gamma_s^{-1}\, \ln \BP(r_\sfp(\ZV)\geq a)=-\sdn(B_s(\sfp,2a)).
$$
\end{theorem}

\section{The typical Voronoi cell revisited}\label{sec:typVor}

In this section we establish a counterpart to Theorem \ref{thmdisintegration} for Voronoi tessellations induced by
an isotropic point process $X$ on $\BS^d$ with intensity $\gamma_s\in (0,\infty)$. The finite Borel measure
$$
\widehat\Lambda\defeq\BE\int_{\BS^d}\indi\{(x,C(x,X))\in\cdot\}\, X(\dint x)
$$
on $\BS^d\times \widehat{\cK}^d_s$ will be considered instead of the intensity measure $\Lambda$ of a general particle process $\zeta$
in Section \ref{sec:typcells}. We define $\sigma(x,C)\defeq(\sigma x,\sigma C)$ for $\sigma\in \Ih$ and $(x,C)\in \BS^d\times\widehat{\cK}^d_s$ and denote by $\sigma\widehat\Lambda$ the image measure of $\widehat\Lambda$ under $\sigma$. It follows from \eqref{rcov} and the isotropy of $X$ that $\sigma\widehat{\Lambda}=\widehat\Lambda$.

\begin{theorem}\label{thm:typVch}
Let $X$ be an isotropic point process on $\BS^d$ with intensity $\gamma_s\in (0,\infty)$.
There exists a unique Borel probability measure $\widehat Q$ on $\BS^d\times \widehat{\cK}^d_s$ which is concentrated on $\{\sfp\}\times\widehat{\cK}^d_s$, invariant under rotations of \ $\BS^d$ fixing $\sfp$ and such that
\begin{equation}\label{eq:chtypVor}
\widehat\Lambda=\gamma_s\int_{\BS^d\times \widehat{\cK}^d_s}\int_{\Ih}\indi\{(\varrho x,\varrho K)\in\cdot\}\,
\lambda(\dint\varrho)\, \widehat Q(\dint(x,K)).
\end{equation}
If $\overline\kappa$ is an arbitrary probability kernel from $\BS^d$ to $\Ih$  such that $\overline\kappa(x,\cdot)$ is concentrated on $\Ih(\sfp,x)$, for $\sigma_d$-a.e.~$x\in\BS^d$, then $\widehat Q$ is given by
\begin{align*}
\widehat Q_{\overline \kappa}&\defeq \frac{1}{\gamma_s}\BE\int_{\BS^d}\int_{\Ih}\indi\{(\sfp,C(\sfp,\varphi^{-1}X))\in\cdot\}
\,\overline\kappa(x,\dint\varphi)\, X(\dint x)\\
&\,\,=\frac{1}{\gamma_s}\int_{\BS^d\times \widehat{\cK}^d_s}\int_{\Ih}\indi\{(\sfp,\varphi^{-1} L)\in\cdot\}\, \overline\kappa(x,\dint\varphi)\, \widehat\Lambda(\dint(x,L)).
\end{align*}
In particular, $\widehat Q_{\overline \kappa}$ is independent of the particular choice of the kernel $\overline\kappa$.
\end{theorem}

\begin{remark}\label{rem:sec10char}{\rm
Except for the uniqueness statement, as a consequence of Theorem \ref{thm:typVch} we obtain corresponding results for the image measures of $\widehat\Lambda$ and $\widehat Q$ under the projection map from $\BS^d\times \widehat{\cK}^d_s$ to $\widehat{\cK}^d_s$. Note that the projection $\widehat Q$ is precisely the distribution $\BP_Z$ of the typical cell of the Voronoi tessellation induced by $X$, introduced in \eqref{eq:deftypV}.
}
\end{remark}

\begin{proof}
It is easy to check that $\widehat Q_{\overline \kappa}$ is well defined and the second representation in terms of $\widehat\Lambda$ is equal to the definition,  by noting that $C(\sfp,\varphi^{-1}X)=\varphi^{-1}C(x,X)$ for $\overline{\kappa}(x,\cdot)$-a.e.~$\varphi\in \Ih$. Moreover it is clear that  $\widehat Q_{ \overline \kappa}$ is concentrated on $\{\sfp\}\times \widehat{\cK}^d_s$.  The proof of the theorem will be accomplished by showing that (1)
$\widehat Q_{  \kappa}$ is invariant under rotations fixing $\sfp$, (2) $\widehat Q_{\overline \kappa}$ satisfies relation \eqref{eq:chtypVor}, (3)  the uniqueness assertion holds, and (4)  $ \widehat Q_{\overline \kappa}=\widehat Q_{  \kappa}$.

For the proof of (1), let $\sigma\in\Ih(\sfp)$ be given. Using the rotation invariance of $\widehat\Lambda$ and basic invariance properties of the kernel $\kappa$, we get
\begin{align*}
\gamma_s\widehat Q_{ \kappa}(\sigma^{-1}(\cdot))&=
\int_{\BS^d\times \widehat{\cK}^d_s}\int_{\Ih}\indi\{(\sfp,\sigma\circ \varphi^{-1} L)\in\cdot\}\,  \kappa(x,\dint\varphi)\, \widehat\Lambda(\dint(x,L))\\
&=
\int_{\BS^d\times \widehat{\cK}^d_s}\int_{\Ih}\indi\{(\sfp,\sigma\circ  \varphi^{-1}\circ\sigma L)\in\cdot\}\,  \kappa(\sigma  x,\dint\varphi)\, \widehat\Lambda(\dint(x,L))\\
&=\int_{\BS^d\times \widehat{\cK}^d_s}\int_{\Ih}\indi\{(\sfp, \sigma\circ ( \sigma\circ\varphi)^{-1}\circ \sigma L) \in\cdot\}\,  \kappa(  x,\dint\varphi)\, \widehat\Lambda(\dint(x,L))\\
&=\int_{\BS^d\times \widehat{\cK}^d_s}\int_{\Ih}\indi\{(\sfp,   ( \varphi\circ\sigma^{-1})^{-1}  L) \in\cdot\}\,  \kappa(  x,\dint\varphi)\, \widehat\Lambda(\dint(x,L))\\
&=\int_{\BS^d\times \widehat{\cK}^d_s}\int_{\Ih}\indi\{(\sfp, \varrho^{-1}  L)\in\cdot\}\,  \kappa(  x,\dint\varrho)\, \widehat\Lambda(\dint(x,L))\\
&=\gamma_s\widehat Q_{ \kappa}( \cdot).
\end{align*}

For the proof of (2), we use again the second representation of $\widehat Q_{\overline \kappa}$,   the right invariance of $\lambda$, the defining properties of $\overline\kappa$, Fubini's theorem and the rotation invariance of $\widehat\Lambda$  to get
\begin{align*}
&\gamma_s\int_{\BS^d\times \widehat{\cK}^d_s}\int_{\Ih}\indi\{(\varrho x,\varrho K)\in\cdot\}\,
\lambda(\dint\varrho)\, \widehat Q_{\overline \kappa}(\dint(x,K))\\
&=\int_{\BS^d\times \widehat{\cK}^d_s}\int_{\Ih}\int_{\Ih}\indi\{(\varrho \sfp,\varrho\circ\varphi^{-1} L)\in\cdot\}\, \lambda(\dint\varrho)\,\overline\kappa(x,\dint\varphi)\,\widehat\Lambda(\dint(x,L))\\
&=\int_{\BS^d\times \widehat{\cK}^d_s}\int_{\Ih}\int_{\Ih}\indi\{(\varrho\circ\varphi \sfp,\varrho L)\in\cdot\}\, \lambda(\dint\varrho)\,\overline\kappa(x,\dint\varphi)\,\widehat\Lambda(\dint(x,L))\\
&=\int_{\BS^d\times \widehat{\cK}^d_s}\int_{\Ih}\indi\{\varrho(x,L)\in\cdot\}\, \lambda(\dint\varrho) \,\widehat\Lambda(\dint(x,L))\\
&=\int_{\Ih}\int_{\BS^d\times \widehat{\cK}^d_s}\indi\{\varrho(x,L)\in\cdot\}\,\widehat\Lambda(\dint(x,L))\, \lambda(\dint\varrho) =\widehat\Lambda,
\end{align*}
since $\lambda(\Ih)=1$.

Now we turn to (3). Let $Q^*$ be any Borel measure on $\BS^d\times \widehat{\cK}^d_s$ which satisfies relation \eqref{eq:chtypVor}, is concentrated on $\{\sfp\}\times \widehat{\cK}^d_s$, and is invariant under rotations of \ $\BS^d$ fixing $\sfp$.  Using these properties (in this order, the first property is used again in the final step), basic properties of the kernel $\kappa$ and Fubini's theorem, we obtain
\begin{align*}
&\frac{1}{\gamma_s}\int_{\BS^d\times \widehat{\cK}^d_s}\int_{\Ih}\indi\{(\sfp,\varphi^{-1} L)\in\cdot\}\,  \kappa(x,\dint\varphi)\, \widehat\Lambda(\dint(x,L))\\
&=\int_{\BS^d\times \widehat{\cK}^d_s}\int_{\Ih}\int_{\Ih}\indi\{  (\sfp,\varphi^{-1}\circ\varrho K)\in\cdot\}\,\kappa(\varrho x,\dint\varphi)\, \lambda(\dint\varrho)\, Q^*(\dint(x,K))\\
&=\int_{\BS^d\times \widehat{\cK}^d_s}\int_{\Ih}\int_{\Ih}\indi\{  (\sfp,\varphi^{-1}\circ\varrho K)\in\cdot\}\,\kappa(\varrho \sfp,\dint\varphi)\, \lambda(\dint\varrho)\, Q^*(\dint(x,K))\\
&=\int_{\BS^d\times \widehat{\cK}^d_s}\int_{\Ih}\int_{\Ih}\indi\{  (\sfp,(\varrho\circ\varphi)^{-1}\circ\varrho K)\in\cdot\}\,\kappa( \sfp,\dint\varphi)\, \lambda(\dint\varrho)\, Q^*(\dint(x,K))\\
&=\int_{\BS^d\times \widehat{\cK}^d_s}\int_{\Ih}\int_{\Ih}\indi\{  (\sfp, \varphi^{-1} K)\in\cdot\}\,\kappa( \sfp,\dint\varphi)\, \lambda(\dint\varrho)\, Q^*(\dint(x,K))\\
&=\int_{\BS^d\times \widehat{\cK}^d_s}\int_{\Ih}\indi\{  (\sfp, \varphi^{-1} K)\in\cdot\}\,\kappa( \sfp,\dint\varphi)\, Q^*(\dint(x,K))\\
&=\int_{\Ih}\int_{\BS^d\times \widehat{\cK}^d_s}\indi\{  (\sfp, \varphi^{-1} K)\in\cdot\}\, Q^*(\dint(x,K))\,\kappa( \sfp,\dint\varphi)
\\
&=\int_{\BS^d\times \widehat{\cK}^d_s}\indi\{  (\sfp,  K)\in\cdot\}
\, Q^*(\dint(x,K))
=Q^*.
\end{align*}

Finally, we prove statement (4). Here we use \eqref{eq:chtypVor} for the measure $\widehat{Q}_{\kappa}$, \eqref{disintlambda} and the fact that
$\widehat{Q}_\kappa$ is concentrated on $\{\sfp\}\times \widehat{\cK}^d_s$ (in the third and the final step), in the fifth and the sixth step we use that for $\varrho_z,\varphi\in\Ih(\sfp,z)$ and $\varrho\in\Ih(\sfp)$  we have $\varrho_z\circ\varrho\sfp=z$ and $\varphi^{-1}\circ\varrho_z\in \Ih(\sfp)$. Moreover we apply Fubini's theorem, the invariance properties of the kernel $\kappa(\sfp,\cdot)$, the invariance property of the measure $\widehat Q_\kappa$, and the fact that $\kappa,\overline\kappa$ are probability kernels. Thus we get
\begin{align*}
\widehat{Q}_{\overline\kappa}
&=\frac{1}{\gamma_s}\int_{\BS^d\times \widehat{\cK}^d_s}\int_{\Ih}\indi\{(\sfp,\varphi^{-1} L)\in\cdot\}\, \overline\kappa(x,\dint\varphi)\, \widehat\Lambda(\dint(x,L))\\
&=\int_{\BS^d\times \widehat{\cK}^d_s}\int_{\Ih}\int_{\Ih}\indi\{(\sfp,\varphi^{-1}\circ\varrho K)\in\cdot\}\, \overline\kappa(\varrho x,\dint\varphi)\,\lambda(\dint\varrho)\, \widehat{Q}_\kappa(\dint(x,K))\\
&=\int_{\BS^d\times \widehat{\cK}^d_s}\int_{\BS^d}\int_{\Ih}\int_{\Ih}\indi\{(\sfp,\varphi^{-1}\circ\varrho K)\in\cdot\}\, \overline\kappa(\varrho \sfp,\dint\varphi)\,\kappa(z,\dint\varrho)\, \sdn(\dint z)\, \widehat{Q}_\kappa(\dint(x,K))\\
&=\int_{\BS^d\times \widehat{\cK}^d_s}\int_{\BS^d}\int_{\Ih}\int_{\Ih}\indi\{(\sfp,(\varphi^{-1}\circ\varrho_z)\circ\varrho K)\in\cdot\}\, \overline\kappa(\varrho_z\circ \varrho \sfp,\dint\varphi)\,\kappa(\sfp,\dint\varrho)\, \sdn(\dint z)\, \widehat{Q}_\kappa(\dint(x,K))\\
&=\int_{\BS^d\times \widehat{\cK}^d_s}\int_{\BS^d}\int_{\Ih}\int_{\Ih}\indi\{(\sfp,(\varphi^{-1}\circ\varrho_z)\circ\varrho K)\in\cdot\}\,\kappa(\sfp,\dint\varrho)\, \overline\kappa(z,\dint\varphi)\, \sdn(\dint z)\, \widehat{Q}_\kappa(\dint(x,K))\\
&=\int_{\BS^d\times \widehat{\cK}^d_s}\int_{\BS^d}\int_{\Ih}\int_{\Ih}\indi\{(\sfp, \varrho K)\in\cdot\}\,\kappa(\sfp,\dint\varrho)\, \overline\kappa(z,\dint\varphi)\, \sdn(\dint z)\, \widehat{Q}_\kappa(\dint(x,K))\\
&=\int_{\BS^d\times \widehat{\cK}^d_s}\int_{\BS^d}\int_{\Ih} \indi\{(\sfp, \varrho K)\in\cdot\}\,\kappa(\sfp,\dint\varrho)\, \sdn(\dint z)\, \widehat{Q}_\kappa(\dint(x,K))\\
&=\int_{\BS^d\times \widehat{\cK}^d_s} \int_{\Ih} \indi\{(\sfp, \varrho K)\in\cdot\}\,\kappa(\sfp,\dint\varrho) \, \widehat{Q}_\kappa(\dint(x,K))\\
&= \int_{\Ih}\int_{\BS^d\times \widehat{\cK}^d_s} \indi\{(\sfp, \varrho K)\in\cdot\} \, \widehat{Q}_\kappa(\dint(x,K))\,\kappa(\sfp,\dint\varrho)\\
&=  \int_{\BS^d\times \widehat{\cK}^d_s} \indi\{(\sfp,   K)\in\cdot\} \, \widehat{Q}_\kappa(\dint(x,K))
=\widehat Q_\kappa.
\end{align*}
Thus we have completed the proofs of (1)--(4), which establishes the theorem.
\end{proof}

\end{document}